\documentclass[english]{article}

\usepackage{epsfig,euscript}

\usepackage{multirow}
\usepackage[table,xcdraw]{xcolor}
\usepackage{geometry}
\usepackage{setspace}
\usepackage{authblk}
\usepackage[utf8]{inputenc}
\usepackage{color,soul,yhmath}
\usepackage{xcolor}
\usepackage{algpseudocode}
\usepackage{bm}
\usepackage{graphicx}
\usepackage{hyperref}
\usepackage{natbib}
\usepackage{bbm}
\usepackage{amsmath}
\usepackage{amssymb}
\usepackage{xr}
\usepackage{booktabs}
\newcommand{\ra}[1]{\renewcommand{\arraystretch}{#1}}

\usepackage[makeroom]{cancel}

\hypersetup{
    colorlinks=false,
    linkcolor=black,
    filecolor=magenta,
    urlcolor=cyan,
    }
 \usepackage{xr}
\makeatletter

\newtheorem{theorem}{Theorem}

\newtheorem{corollary}{Corollary}

\newtheorem{lemma}{Lemma}

\newtheorem{proposition}{Proposition}
\newtheorem{remark}{Remark}

\numberwithin{equation}{section}
\def\c#1{\mathcal{#1}}
\def\b#1{\mathbbm{#1}}
\def\bf#1{\mathbf{#1}}

\def\mat#1#2{\text{mat}_{#1}({#2})}
\def\red#1{\textcolor{red}{#1}}

\def\wt#1{\widetilde{#1}}
\def\wh#1{\widehat{#1}}
\def\vec#1{\textbf{vec}\big(#1\big)}

\def\cY{\mathcal{Y}}
\def\cM{\mathcal{M}}
\def\cC{\mathcal{C}}
\def\cF{\mathcal{F}}
\def\cE{\mathcal{E}}

\def\cN{\mathcal{N}}

\def\A{\mathbf{A}}
\def\B{\mathbf{B}}
\def\C{\mathbf{C}}

\def\Q{\mathbf{Q}}
\def\S{\mathbf{S}}
\def\D{\mathbf{D}}
\def\Y{\mathbf{Y}}

\def\F{\mathbf{F}}
\def\E{\mathbf{E}}
\def\H{\mathbf{H}}
\def\v{\mathbf{v}}
\def\R{\mathbf{R}}

\def\I{\mathbf{I}}
\def\X{\mathbf{X}}
\def\Z{\mathbf{Z}}
\def\0{\mathbf{0}}
\def\1{\mathbf{1}}

\def\diag{\textnormal{diag}}
\def\Amk{\A_{\text{-}k}}
\def\bLambda{\boldsymbol{\Lambda}}
\def\bSigma{\boldsymbol{\Sigma}}
\def\bgamma{\boldsymbol{\gamma}}
\def\bepsilon{\boldsymbol{\epsilon}}
\def\tr{\textnormal{tr}}
\DeclareMathOperator*{\plim}{plim}
\def\dmk{d_{\text{-}k}}
\def\rmk{r_{\text{-}k}}
\def\Amk{\A_{\text{-}k}}

\def\Aemk{\A_{e,\text{-}k}}
\newcommand{\norm}[1]{\big\| #1 \big\|}

\renewcommand{\tilde}{\widetilde}

\def\6bullets{\bullet\bullet\bullet\bullet\bullet\bullet}

\begin{document}
\setlength{\parindent}{18pt}
	\doublespacing
	\begin{titlepage}
		
		\title{Tensor Time Series Imputation through Tensor Factor Modelling}
		        \author{Zetai Cen\thanks{Zetai Cen is PhD student, Department of Statistics, London School of Economics. Email: Z.Cen@lse.ac.uk}}
				\author{Clifford Lam\thanks{Clifford Lam is Professor, Department of Statistics, London School of Economics. Email: C.Lam2@lse.ac.uk}}
		
		\affil{Department of Statistics, London School of Economics and Political Science}
		
		\date{}
		
		\maketitle

\begin{abstract}
We propose tensor time series imputation when the missing pattern in the tensor data can be general, as long as any two data positions along a tensor fibre are both observed for enough time points. The method is based on a tensor time series factor model with Tucker decomposition of the common component. One distinguished feature of the tensor time series factor model used is that there can be weak factors in the factor loadings matrix for each mode. This reflects reality better when real data can have weak factors which drive only groups of observed variables, for instance, a sector factor in financial market driving only stocks in a particular sector. Using the data with missing entries, asymptotic normality is derived for rows of estimated factor loadings, while consistent covariance matrix estimation enables us to carry out inferences. As a first in the literature, we also propose a ratio-based estimator for the rank of the core tensor under general missing patterns. Rates of convergence are spelt out for the imputations from the estimated tensor factor models. Simulation results show that our imputation procedure works well, with asymptotic normality and corresponding inferences also demonstrated. Re-imputation performances are also gauged when we demonstrate that using slightly larger rank then estimated gives superior re-imputation performances. A Fama-French portfolio example with matrix returns and an OECD data example with matrix of Economic indicators are presented and analyzed, showing the efficacy of our imputation approach compared to direct vector imputation.
\end{abstract}
		
		\bigskip
		\bigskip


\noindent

	\end{titlepage}
	
	\setcounter{page}{2}

\maketitle


\newpage
\section{Introduction}\label{sec:introduction}
Large dimensional panel data is easier to obtain than ever thanks to a quickly evolving internet speed and more diverse download platforms. Together with the advancement of statistical analyses for these data over the past decade, researchers also open up more to time series data with higher order, namely, tensor time series data. A prominent example would be order 2 tensor time series, i.e., matrix-valued time series. \cite{Wangetal2019} proposes a factor model using a Tucker decomposition of the common component in the modelling.
An example on monthly import-export volume of products among different countries is given in \cite{Chenetal2022}, where factor modelling using Tucker decomposition is explored, and generalized to higher order tensors. Focusing on matrix-valued time series, \cite{Changetal2023} proposes a tensor-CP decomposition for modelling the data. \cite{Zhang2024} and \cite{Chenetal2021} propose autoregressive models for matrix-valued time series.
For a more comprehensive review on matrix-valued time series analysis, please refer to \cite{Tsay2023}.

A less addressed topic in large time series analysis is the treatment of missing data, in particular,  imputation of missing data and the corresponding inferences. While there are numerous data-centric methods in various scientific fields for imputing multivariate time series data (see \cite{Chaponetal2023} for environmental time series, \cite{KazijevsSamad2023} for health time series, \cite{Zhaoetal2023} and  \cite{Zhangetal2021}
for using deep-learning related architectures for imputations, to name but a few),
almost none of them address statistically how accurate their methods are, and all of them are not for higher order tensor time series. We certainly can line up the variables in a tensor time series to make it a longitudinal panel, but in doing so we lose special structures and insights that can be utilized for forecasting and interpretation of the data. More importantly, transforming a moderate sized tensor to a vector means the length of the vector can be much larger than the sample size, creating curse of dimensionality.

For imputing large panel of time series with statistical analyses, \cite{Bai_Ng} defines the concept of TALL and WIDE blocks of data and proposes an iterative TW algorithm in imputing missing values in a large panel, while \cite{Cahanetel2023} improves the TW algorithm to a Tall-Project (TP) algorithm so that there is no iterations needed. Both papers use factor modelling for the imputations, and derive rates of convergence when all factors are pervasive and the number of factors known. Asymptotic normality for rows of estimated factor loadings and the corresponding practical inferences are also developed as well. \cite{Xiong_Pelger} also bases their imputations on a factor model for a large panel of time series with pervasive factors and number of factors known, and build a method for imputing missing values under very general missing patterns, with asymptotic normality and inferences also developed.

To the best of our knowledge, for tensor time series with order larger than 1 (i.e, at least matrix-valued), there are no theoretical analyses on imputation performances. Imputation methodologies developed on tensor time series are also scattered around very different applications. See \cite{Chenetal2022b} on traffic tensor data and  \cite{Panetal2021} for RNA-sequence tensor data for instance.

In view of all the above, as a first in the literature, we aim to develop a tensor imputation method accompanied by theoretical analyses in this paper.
Like \cite{Cahanetel2023}, we use factor modelling for tensor time series as a basis for our imputation method.  Unlike \cite{Cahanetel2023}, \cite{Bai_Ng} or \cite{Xiong_Pelger} though, we develop a method that can consistently estimate the number of factors, or the core tensor rank, in a Tucker decomposition-based factor model for the tensor time series with missing values. Our method can be considered a combination of \cite{Heetal2022a} for the tensor factor model, and \cite{Xiong_Pelger} for the imputation methodology with general missingness. In Section \ref{sec:factormodelwithmissingness}, we introduce two motivating examples and  our methodology at the same time. One is the Fama-French portfolio return data with missing  entries, to be analyzed in Section \ref{subsec:famafrench}. The other is a set of monthly and quarterly OECD Economic indicators, with missingness naturally occurring for the quarterly recorded indicators relative to the monthly ones. We analyze this set of OECD data in Section  \ref{subsec:OECD}.

As a further contribution, we also allow factors to be weak. A weak factor corresponds to a column in a factor loading matrix being sparse, or approximately sparse. This implies that not all units in a tensor has dynamics contributed by all the factors inside the core tensor. In \cite{Chen_Lam}, they allow for weak factors in their analyses, and discovers that there are potentially weak factors in the NYC taxi traffic data, which we are going to analyze in the supplementary materials of this paper.
We prove consistency of our imputations under general missingness, and develop asymptotic normality and practical inferences for rows of factor loading matrix estimators, with rates of convergence in all consistency results spelt out. Our method is available in the R package \texttt{tensorMiss}, which has used the \texttt{Rcpp} package to greatly boost computational speed.

The rest of the paper is organised as follows. Section \ref{sec:notations} introduces the notations used in this paper. Section \ref{sec:factormodelwithmissingness} presents the Fama-French portfolio returns data and the OECD data as two motivating examples, before describing the tensor factor model and the imputation methodology we employ.  Section \ref{sec:Assumption_theories} lays down the main assumptions for the paper, with consistent estimation and rates of convergence of all factor loading matrix estimators and imputed values presented. Asymptotic normality and the estimators of the corresponding covariance matrices for practical inferences are introduced as well in Section \ref{subsec:asymptotic_notmality}, before our proposed ratio-based estimators for the number of factors in Section \ref{subsec:estimate_number_of_factors}. Section \ref{sec:numerical} presents extensive simulation results for our paper, together with an analysis for the Fama-French portfolio return data in Section \ref{subsec:famafrench} and an analysis for the OECD Economic data in Section \ref{subsec:OECD}.
The NYC taxi traffic data, together with extra simulations are presented in the supplementary materials for the paper.
All proofs are in the supplementary materials associated with this paper.

\section{Notations}\label{sec:notations}
Throughout this paper, we use the lower-case letter, bold lower-case letter, bold capital letter, and calligraphic letter, i.e., $x,\bf{x},\bf{X},\c{X}$, to denote a scalar, a vector, a matrix, and a tensor respectivel;y.
We also use $x_i, X_{ij}, \bf{X}_{i\cdot}, \bf{X}_{\cdot i}$ to denote, respectively, the $i$-th element of $\bf{x}$, the $(i,j)$-th element of $\bf{X}$, the $i$-th row vector (as a column vector) of $\bf{X}$, and the $i$-th column vector of $\bf{X}$. We use $\otimes$ to represent the Kronecker product, and $\circ$ the Hadamard product. We use $a\asymp b$ to denote $a=O(b)$ and $b=O(a)$. Hereafter, given a positive integer $m$, define $[m]:=\{1,2,\dots,m\}$. The $i$-th largest eigenvalue of a matrix $\X$ is denoted by $\lambda_i(\bf{X})$. The notation $\bf{X}\succcurlyeq 0$ (resp. $\bf{X} \succ 0$) means that $\bf{X}$ is positive semi-definite (resp. positive definite). We use $\bf{X}'$ to denote the transpose of $\bf{X}$, and $\diag(\X)$ to denote a diagonal matrix with the diagonal elements of $\X$, while $\diag(\{x_1, \dots, x_n\})$ represents the diagonal matrix with $\{x_1, \dots, x_n\}$ on the diagonal.

\textbf{Norm notations}:
For a given set, we denote by $|\cdot|$ its cardinality. We use $\|\bf{\cdot}\|$ to denote the spectral norm of a matrix or the $L_2$ norm of a vector, and $\|\bf{\cdot}\|_F$ to denote the Frobenius norm of a matrix. We use $\|\cdot\|_{\max}$ to denote the maximum absolute value of the elements in a vector, a matrix or a tensor. The notations $\|\cdot\|_1$ and $\|\cdot\|_{\infty}$ denote the $L_1$ and $L_{\infty}$-norm of a matrix respectively, defined by $\|\X\|_{1} := \max_{j}\sum_{i}|X_{ij}|$ and $\|\X\|_{\infty} := \max_{i}\sum_{j}|X_{ij}|$. WLOG, we always assume the eigenvalues of a matrix are arranged by descending orders, and so are their corresponding eigenvectors.

\begin{figure}[!htp]
	\includegraphics[width=15cm]{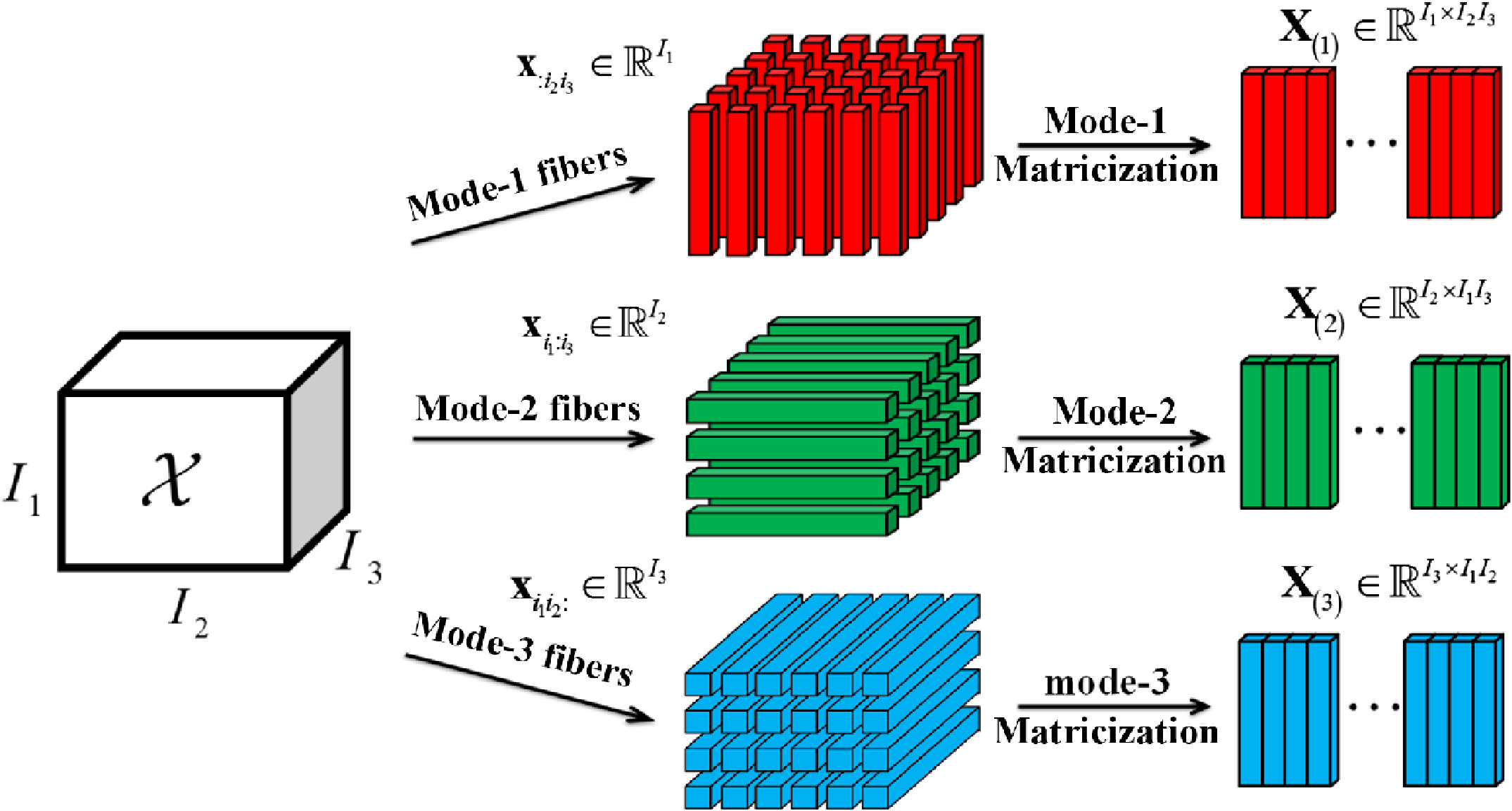}\\
	\caption{Illustration of the mode-$k$ fibers and its corresponding unfolding matrix.}
	\label{Fig: tensor_example}
\end{figure}

\textbf{Tensor-related notations}: For the rest of this section, we briefly introduce the notations and operations for tensor data, which will be sufficient for this paper. For more details on tensor manipulations, readers are referred to \cite{Kolda_Bader}. A multidimensional array with $K$ dimensions is an \textit{order}-$K$ tensor, with its $k$-th dimension termed as \textit{mode}-$k$. For an order-$K$ tensor $\c{X} = (X_{i_1,\ldots,i_K}) \in \mathbb{R}^{I_1\times\cdots\times I_K}$, a column vector $(X_{i_1, \ldots, i_{k-1}, i, i_{k+1}, \ldots, i_K})_{i\in[I_k]}$ represents a \textit{mode-$k$ fibre} for the tensor $\c{X}$. We denote by $\mat{k}{\c{X}} \in \mathbb{R}^{I_k\times I_{\text{-}k}}$ (or sometimes $\bf{X}_{(k)}$, with $I_{\text{-}k} := (\prod_{j=1}^K I_j)/I_k$) the \textit{mode-$k$ unfolding/matricization} of a tensor, defined by placing all mode-$k$ fibres into a matrix. See Figure \ref{Fig: tensor_example} for an illustration (figure from \cite{Taoetal2019}).

We denote by $\c{X} \times_k \A$ the \textit{mode-$k$ product} of a tensor $\c{X}$ with a matrix $\A$, defined by
\[\mat{k}{\c{X}\times_k \A} := \A\, \mat{k}{\c{X}}.\]
Finally, we use the notation $\Vec{\cdot}$ to denote the vectorization of a matrix or the vectorization of the mode-1 unfolding of a tensor.

\section{Two Motivating Examples and The Imputation Procedure}\label{sec:factormodelwithmissingness}
We first describe two motivating data examples in Section \ref{subsec:fama-french_data} and \ref{subsec:OECD_data}, before presenting our imputation procedure for a general order-$K$ mean zero tensor $\c{Y}_t = (\cY_{t,i_1, \dots,i_K})\in\b{R}^{d_1\times d_2\times \dots \times d_K}$ for each $t\in[T]$. The two data examples will be analyzed in details in Section \ref{subsec:famafrench} and \ref{subsec:OECD} respectively.

\subsection{Example: The Fama-French portfolio returns}\label{subsec:fama-french_data}
This is a set of Fama-French portfolio returns data with missingness. Stocks are categorized into ten levels of market equity (ME) and ten levels of book-to-equity ratio (BE) which is the book equity for the last fiscal year divided by the end-of-year ME. At the end of June each year, both ME and BE use NYSE deciles as breakpoints, with stocks allocated accordingly. Moreover, the stocks in each of the $10\times 10$ categories form exactly two portfolios, one being value weighted, and the other of equal weight. Hence, there are two sets of 10 by 10 portfolios with their time series to be studied. We use monthly data from January 1974 to June 2021, so that $T =570$, and for both value weighted and equal weighted portfolios we have each of our data set as an order-2 tensor $\c{X}_t\in \b{R}^{10\times 10}$ for $t\in[570]$.
For more details, please visit

\url{https://mba.tuck.dartmouth.edu/pages/faculty/ken.french/Data_Library/det_100_port_sz.html}.

If no stocks are allocated to a category at a timestamp, the corresponding return data is unavailable and hence missing. It is reasonable to argue the missingness might depend on the rows of the loading matrix, i.e., extreme categories tend to contain fewer stocks, but independent of latent factors and noise. The total number of missing entries is $161$ and hence the percentage of missing is $161/(10\times 10\times 570) =0.28\%$ for both the value weighted and equal weighted series. However, the irregular missing pattern here can be harmful if we are after a complete case analysis. For full observation after a timestamp, we may only start from July 2009 and hence $74.7\%$ of the data would be ditched. On the other hand, we might ditch four categories to obtain a complete data set but lose the potential insights on the return series of the four categories.

\subsection{Example: OECD Economic indicators}\label{subsec:OECD_data}
In this example, we study a group of economic indicators for a selection of countries obtained from the Organization for Economic Co-operation and Development (OECD). The data consists of monthly/quarterly observations of 11 economic indicators: current account balance as percentage of GDP (CA-GDP), consumer price index (CP), merchandise exports (EX), merchandise imports (IM), short-term interest rates (IR3TIB), long-term interest rates (IRLT), interbank rates (IRSTCI), producer price index (PP), production volume (PRVM), retail trade volume (TOVM) and unit labour cost (ULC). They are observed for 17 countries: Belgium (BEL), Canada (CAN), Denmark (DNK), Finland (FIN), France (FRA), Germany (DEU), Greece (GRC), Italy (ITA), Luxembourg (LUX), Netherlands (NLD), Norway (NOR), Portugal (PRT), Spain (ESP), Sweden (SWE), Switzerland (CHE), United Kingdom (GBR) and United States (USA), with data spanning from January 1971 to December 2023. We correspond respectively rows and columns to countries and indicators, so that we have our data as an order-2 tensor $\c{Y}_t \in \b{R}^{17\times 11}$ for $t\in[636]$. For more details, see key short-term economic indicators available at \url{https://data.oecd.org/}.

The data is naturally missing for three reasons: unavailable indicator records for some countries at early time periods, quarterly indicators are only available at the end of each quarter, and are sometimes unrecorded. Similar to the Fama-French data, we suppose the missing pattern is dependent on the loading matrices by arguing that relatively less important indicators are only available quarterly. The percentage of missing data is $26.2\%$, which leads to significantly inefficient use of data if we hope to analyse a set of complete data. The fact that the data is observed at least quarterly in the long run ensures the existence of a lower bound on the proportion of available data, which in turn satisfies Assumption (O1) in Section \ref{sec:Assumption_theories}.

\subsection{The model and the imputation procedure}\label{subsec:model_imputation}
\noindent\textbf{The Model}: Suppose the order-$k$ mean zero tensor $\cY_t$ is modelled by
\begin{equation}
    \label{eqn: tenfac-all}
    \c{Y}_t = \c{C}_t+\c{E}_t= \c{F}_t \times_1 \bf{A}_{1}\times_2 \bf{A}_{2}\times_3 \dots\times_{K} \bf{A}_K
    + \c{E}_t, \;\;\; t\in[T],
\end{equation}
where $\c{C}_t$ is the common component and $\c{E}_t$ the error tensor. The core tensor is $\c{F}_t\in\b{R}^{r_1\times r_2\times \dots \times r_K}$, and each mode-$k$ factor loading matrix $\bf{A}_{k}$ has dimension $d_k \times r_k$. See \cite{Heetal2022a} amongst others using the same tensor factor model. Using the QR decomposition, if we can decompose $\A_k = \Q_k\Z_k^{1/2}$ (see Assumption (L1) in Section \ref{subsec:Assumptions} for details), then (\ref{eqn: tenfac-all}) can be written as
\begin{equation}
\begin{split}\label{eqn: tenfac-all2}
  \c{Y}_t &= \c{F}_{Z,t}\times_1 \Q_1 \times_2 \cdots \times_K \Q_K + \cE_t, \;\;\; t\in[T], \; \text{ where }\\
  \c{F}_{Z,t} &:=  \c{F}_t\times_1 \bf{Z}_1^{1/2}\times_2 \dots\times_{K} \bf{Z}_K^{1/2}.
\end{split}
\end{equation}
Model (\ref{eqn: tenfac-all}) is an extension to the usual time series factor model ($K=1$):
\[\c{Y}_t = \mat{1}{\c{Y}_t} = \mat{1}{\c{F}_t\times_1 \A_1} + \mat{1}{\cE_t} = \A_1\mat{1}{\cF_t} + \mat{1}{\cE_t} = \A_1\cF_t + \cE_t,\]
and also for a matrix-valued time series factor model ($K=2$):
\begin{align*}
  \cY_t &= \mat{1}{\cY_t} = \A_1\mat{1}{\cF_t}\A_2' + \mat{1}{\cE_t} = \A_1\cF_t\A_2' + \cE_t.
\end{align*}

\noindent\textbf{The Imputation Procedure}: We only observe partial data. Define the missingness tensor $\c{M}_t = (\c{M}_{t,i_1, \dots,i_K})\in\b{R}^{d_1\times d_2\times \dots \times d_K}$ with
\begin{align*}
  \c{M}_{t,i_1, \dots,i_K} = \left\{
                           \begin{array}{ll}
                             1, & \hbox{if $\cY_{t,i_1, \dots,i_K}$ is observed;} \\
                             0, & \hbox{otherwise.}
                           \end{array}
                         \right.
\end{align*}
Our aim is to recover the value for the common component $\c{C}_{t,i_1,\dots,i_K}$ if $\c{M}_{t,i_1,\dots,i_K}=0$. Assuming first the number of factors $r_k$ is known for all modes, we want to obtain the estimators of the factor loading matrices, $\wh{\bf{Q}}_k$ for $k\in[K]$, and then the estimated core tensor series $\wh{\c{F}}_{Z,t}$ for $t\in[T]$. See (\ref{eqn: tenfac-all2}) for the definition of $\Q_k$ and $\c{F}_{Z,t}$.
We can then estimate the common components at time $t$ by
\begin{align}
  \wh{\cC}_t &= \wh{\cF}_{Z,t} \times_1 \wh{\Q}_{1}\times_2 \cdots\times_{K}\wh{\Q}_K. \label{eqn:est_common_components}
\end{align}
With (\ref{eqn:est_common_components}), we can impute $\cY_t$ using
\begin{align*}
\wt\cY_{t,i_1,\ldots,i_K} = \left\{
  \begin{array}{ll}
    \cY_{t,i_1,\ldots,i_K}, & \hbox{if $\cM_{t,i_1,\ldots,i_K} = 1$;} \\
    \wh{\cC}_{t,i_1,\ldots,i_K}, & \hbox{if $\cM_{t,i_1,\ldots,i_K} = 0$.}
  \end{array}
\right.
\end{align*}
We leave the discussion of estimating $r_k$ to Section \ref{subsec:estimate_number_of_factors}. See Section \ref{subsec:est_Q} in how to obtain $\wh\Q_k$ and Section \ref{subsec:imputation_procedure} in how to obtain $\wh\cF_{Z,t}$.

\subsection{Estimation of factor loading matrices}\label{subsec:est_Q}
In this paper, we make use of the following notation:
\begin{equation}\label{eqn:set_psi}
    \psi_{k,ij,h} := \Big\{t\in[T]\mid \mat{k}{\c{M}_t}_{ih}\mat{k}{\c{M}_t}_{jh}=1\Big\}.
\end{equation}
Hence $\psi_{k,ij,h}$ is the set of time periods where both the $i$-th and $j$-th entries of the $h$-th mode-$k$ fibre are observed, $i,j\in[d_k]$, $h\in[\dmk]$ with $\dmk := d_1\cdots d_K/d_k$.

Inspired by \cite{Xiong_Pelger} for a vector time series panel, our method relies on the reconstruction of the mode-$k$ sample covariance matrix $\S_k$, defined for $i,j\in[d_k]$,
\begin{equation}
    \label{eqn: second-order-moment-entry}
    \begin{split}
        (S_k)_{ij} &:= \frac{1}{T}\sum_{t=1}^T\mat{k}{\c{Y}_t}_{i\cdot}'
        \mat{k}{\c{Y}_t }_{j\cdot}
        =\sum_{h=1}^{d_{\text{-}k}}\frac{1}{T}\sum_{t=1}^T
        \mat{k}{\c{Y}_t }_{ih}
        \mat{k}{\c{Y}_t }_{jh}.
    \end{split}
\end{equation}
With missing entries characterized by $\cM_t$ and $\psi_{k,ij,h}$ in (\ref{eqn:set_psi}),
we can generalize the above to
\begin{equation}
    \label{eqn: weighted-Sk-entry}
    \begin{split}
        (\wh{S}_k)_{ij} =&
        \sum_{h=1}^{d_{\text{-}k}}
        \Bigg\{\frac{1}{|\psi_{k,ij,h}|}\sum_{t\in\psi_{k,ij,h}}
        \mat{k}{\c{Y}_t}_{ih}
        \mat{k}{\c{Y}_t}_{jh}
         \Bigg\}.
    \end{split}
\end{equation}
Intuitively, the cross-covariance between unit $i$ and $j$ at the $h$-th mode-$k$ fibre is estimated inside the curly bracket in (\ref{eqn: weighted-Sk-entry}) using only the corresponding available data. PCA can now be performed on $\wh{\bf{S}}_k$, and $\wh{\bf{Q}}_k$ is obtained as the first $r_k$ eigenvectors of $\wh{\bf{S}}_k$.

\subsection{Estimation of the core tensor series}\label{subsec:imputation_procedure}
With $\wh\Q_k$ available (which is estimating the factor loading space of $\Q_k$, with $\wh{\Q}_k$ having orthonormal columns), we can estimate $\cF_{Z,t}$ (equivalently $\vec{\cF_{Z,t}}$) by observing from model (\ref{eqn: tenfac-all2}) that
\begin{align*}
  \vec{\cY_t} &= \Q_\otimes \vec{\cF_{Z,t}} + \vec{\cE_t}, \;\;\; \text{where } \Q_\otimes := \Q_K\otimes\cdots\otimes\Q_1.
\end{align*}
If $\Q_\otimes$ is known, then the least squares estimator of $\vec{\cF_{Z,t}}$ is given by
\begin{align*}
  \vec{\cF_{Z,t}} &= (\Q_\otimes'\Q_\otimes)^{-1}\Q_\otimes'\vec{\cY_t}
= \Bigg(\sum_{j=1}^d \Q_{\otimes,j\cdot}\Q_{\otimes,j\cdot}'\Bigg)^{-1}
\Bigg(\sum_{j=1}^d\Q_{\otimes,j\cdot}[\vec{\cY_t}]_j\Bigg).
\end{align*}
With missing data, using the missingness tensor $\cM_t$, the above can be generalized to
\begin{equation}
    \label{eqn: Ft-regression}
     \vec{\wh\cF_{Z,t}} =
     \Bigg(\sum_{j=1}^d [\vec{\cM_t}]_j
     \wh\Q_{\otimes, j\cdot}\wh\Q_{\otimes, j\cdot}' \Bigg)^{-1}
     \Bigg(\sum_{j=1}^d [\vec{\cM_t}]_j
     \wh\Q_{\otimes, j\cdot}
     [\vec{\cY_t}]_j \Bigg).
\end{equation}

\section{Assumptions and Theoretical Results}\label{sec:Assumption_theories}
We present our assumptions for consistent imputation and estimation of factor loading matrices, with the corresponding theoretical results presented afterwards.
\subsection{Assumptions}\label{subsec:Assumptions}
\begin{itemize}
    \item[(O1)] (Observation patterns)
\textit{\\
1. $\c{M}_t$ is independent of $\c{F}_s$ and $\c{E}_s$ for any $t,s\in[T]$.\\
2. Given $\c{M}_t$ with $t\in[T]$, for any $k\in[K], i,j\in[d_k],h\in[d_{\text{-}k}]$, there exists a constant $\psi_0$ such that
$$
\frac{|\psi_{k,ij,h}|}{T}\geq \psi_0>0.
$$
}
\end{itemize}
\begin{itemize}
    \item[(M1)] (Alpha mixing)
\textit{The elements in $\c{F}_t$ and $\c{E}_t$ are $\alpha$-mixing. A vector process $\{\bf{x}_t: t=0,\pm 1, \pm2,\dots\}$ is $\alpha$-mixing if, for some $\gamma>2$, the mixing coefficients satisfy the condition that
\begin{equation*}
    \sum_{h=1}^\infty
    \alpha(h)^{1-2/\gamma}<\infty,
\end{equation*}
where $\alpha(h)=\sup_\tau\sup_{A\in\c{H}_{-\infty}^\tau,
B\in\c{H}_{\tau+h}^\infty}|\b{P}(A\cap B)-\b{P}(A)\b{P}(B)|$ and $\c{H}_\tau^s$ is the $\sigma$-field generated by $\{\bf{x}_t: \tau\leq t\leq s\}$.
}
\end{itemize}
\begin{itemize}
    \item[(F1)] (Time series in $\c{F}_t$)
\textit{There is $\c{X}_{f,t}$ the same dimension as $\c{F}_{t}$, such that $\c{F}_{t}=\sum_{q\geq 0}a_{f,q}\c{X}_{f,t-q}$. The time series $\{\c{X}_{f,t}\}$ has i.i.d. elements with mean $0$ and variance $1$, with uniformly bounded fourth order moments. The coefficients $a_{f,q}$ are such that $\sum_{q\geq 0}a_{f,q}^2=1$ and $\sum_{q\geq 0}|a_{f,q}|\leq c$ for some constant $c$.
}
\item[(L1)] (Factor strength)
\textit{We assume for $k\in[K]$, $\bf{A}_k$ is of full column rank and independent of factors and errors series. Furthermore, as $d_k\to\infty$,
\begin{equation}
\label{eqn: L1}
    \bf{Z}_k^{-1/2}\bf{A}_k'\bf{A}_k\bf{Z}_k^{-1/2}
    \to \bf{\Sigma}_{A,k},
\end{equation}
where $\bf{Z}_k=\textnormal{diag}(\bf{A}_k'\bf{A}_k)$ and $\bf{\Sigma}_{A,k}$ is positive definite with all eigenvalues bounded away from 0 and infinity. We assume $(\bf{Z}_k)_{jj}\asymp d_k^{\alpha_{k,j}}$ for $j\in[r_k]$, and $1/2<\alpha_{k,r_k}\leq \dots\leq \alpha_{k,2}\leq \alpha_{k,1}\leq 1$.
}
\end{itemize}
With Assumption (L1), we can denote $\bf{Q}_k:=\bf{A}_k\bf{Z}_k^{-1/2}$ and hence $\bf{Q}_k'\bf{Q}_k\to\bf{\Sigma}_{A,k}$. We need $\alpha_{k,j} > 1/2$ in order that the ratio-based estimator of the number of factors in Section \ref{subsec:estimate_number_of_factors} works.

\begin{itemize}
    \item[(E1)] (Decomposition of $\c{E}_t$)
\textit{We assume $K$ is constant, and
\begin{equation}
\label{eqn: NE1}
    \c{E}_t = \c{F}_{e,t}\times_1 \A_{e,1}\times_2 \cdots \times_K\A_{e,K}
    + \bSigma_{\epsilon}\circ\bm\epsilon_t,
\end{equation}
where $\c{F}_{e,t}$ is an order-$K$ tensor with dimension $r_{e,1}\times\cdots\times r_{e,K}$, containing independent elements with mean $0$ and variance $1$. The order-$K$ tensor $\bm\epsilon_t\in\b{R}^{d_1\times \cdots\times d_K}$ contains independent mean zero elements with unit variance, with the two time series $\{\bm\epsilon_t\}$ and $\{\c{F}_{e,t}\}$ being independent. The order-$K$ tensor $\bSigma_{\epsilon}$ contains the standard deviations of the corresponding elements in $\bepsilon_t$, and has elements uniformly bounded.
}

\textit{
Moreover, for each $k\in[K]$, $A_{e,k} \in \mathbb{R}^{d_k\times r_{e,k}}$ is such that
$\norm{\A_{e,k}}_1=O(1)$. That is, $\A_{e,k}$ is (approximately) sparse.
}



\end{itemize}
\begin{itemize}
    \item[(E2)] (Time series in $\c{E}_t$)
\textit{There is $\c{X}_{e,t}$ the same dimension as $\c{F}_{e,t}$, and $\c{X}_{\epsilon,t}$ the same dimension as $\bm\epsilon_t$, such that $\c{F}_{e,t}=\sum_{q\geq 0}a_{e,q}\c{X}_{e,t-q}$ and $\bm\epsilon_t=\sum_{q\geq 0}a_{\epsilon,q}\c{X}_{\epsilon, t-q}$, with $\{\c{X}_{e,t}\}$ and $\{\c{X}_{\epsilon, t}\}$ independent of each other, and each time series has independent elements with mean $0$ and variance $1$ with uniformly bounded fourth order moments. Both $\{\c{X}_{e,t}\}$ and $\{\c{X}_{\epsilon, t}\}$ are independent of $\{\c{X}_{f,t}\}$ from (F1).}

\textit{The coefficients $a_{e,q}$ and $a_{\epsilon, t}$ are such that $\sum_{q\geq 0}a_{e,q}^2=\sum_{q\geq 0}a_{\epsilon,q}^2=1$ and $\sum_{q\geq 0}|a_{e,q}|, \sum_{q\geq 0}|a_{\epsilon,q}|\leq c$ for some constant $c$.
}


\item[(R1)] (Further rate assumptions) {\em
We assume that, with $d := d_1\cdots d_K$ and $g_s := \prod_{k=1}^Kd_k^{\alpha_{k,1}}$,
\begin{align*}
dg_s^{-2}T^{-1}d_k^{2(\alpha_{k,1} - \alpha_{k,r_k})+1} &= o(1), \;\;\;
dg_s^{-1}T^{-1}d_k^{2(\alpha_{k,1} - \alpha_{k,r_k})}=o(1),\;\;\; 
dg_s^{-1}d_k^{\alpha_{k,1}-\alpha_{k,r_k}-1/2} = o(1). 
\end{align*}
}
\end{itemize}
Assumption (O1) means that the missing mechanism is independent of the factors and the noise series, which is also assumed in \cite{Xiong_Pelger} for the purpose of identification. It also means that the missing pattern can depend on the $K$ factor loading matrices, allowing for a wide variety of missing patterns that can vary over time and units in different dimensions. Condition 2 of (O1) implies that the number of time periods that any two individual units are both observed are at least proportional to $T$, which simplifies proofs and presentations, and is also used in \cite{Xiong_Pelger}. Assumption (M1) is a standard assumption in vector time series factor models, which facilitates proofs using central limit theorem for time series without losing too much generality. Assumption (F1), (E1) and (E2) are exactly the corresponding assumptions in \cite{Chen_Lam}, allowing for serial correlations in the factor series, and serial and cross-sectional dependence within and among the error tensor fibres. These three assumptions facilitate the proof of asymptotic normality in Section \ref{subsec:asymptotic_notmality}, and boil down to similar assumptions in \cite{ChenFan2023} for matrix time series and in \cite{Barigozzietal2023b} for general tensor time series (see Proposition \ref{Prop:assumption_implications} in the supplementary materials for the technical details). Together with Assumption (M1), we implicitly restrict the general linear processes in (F1) and (E2) to be, for instance, of short rather than long dependence.


Assumption (L1) is quite different from assumptions in other papers on factor models, in the sense that we allow for the existence of weak factors alongside the pervasive ones. \cite{Chen_Lam} adapted the same assumption, which allows each column of $\A_k$ to be completely dense (i.e., a pervasive factor) or sparse to a certain extent. A diagonal entry in $\Z_k$ then records how dense a column really is, and the corresponding strength of factors defined.
{
Assumption (L1) is similar to, yet more general than, Assumption 1(iii) in \cite{Onatski2012} which requires $\bSigma_{A,k}$ to be diagonal. If all factors are pervasive, (\ref{eqn: L1}) can be read as $d_k^{-1} \A_k' \A_k \to \bf{\Sigma}_{A,k}$ which is akin to Assumption 3 of \cite{ChenFan2023} for $K=2$.
}
Modelling with weak factors is closer to reality, and empirical evidence can be found in economics and finance, etc. For instance, apart from a pervasive market factor, there can be weaker sector factors in a large selection of stock returns \citep{Trzcinka1986}. More recent work on factor models specifically focuses on weak factors with real data examples confirming the existence of weak factors, such as \cite{Freyaldenhoven2022} and \cite{ChenLam2024a}.


%

Finally, Assumption (R1) gives the technical rates needed for the proof of various theorems in the paper because of the existence of weak factors. If all factors are pervasive (i.e., $\alpha_{k,j}=1$), then the conditions are automatically satisfied. Suppose $K=2$, $T\asymp d_1 \asymp d_2$ and the strongest factors are all pervasive (i.e., $\alpha_{k,1}=1$), then we need $\alpha_{k,r_k} > 1/2$ for (R1) to be satisfied. This condition is the same as the one remarked right after we stated Assumption (L1). A factor with $\alpha_{k,j}$ close to 0.5 presents a significantly weak factor with only more than $d_k^{1/2}$ of elements are non-zero in the corresponding column of $\A_k$.

\begin{remark}\label{remark:1}
With the missing entries imputed using the estimated common components $\wh\cC_{t,i_1,\dots,i_K}$, we have a completed data set which could be used for re-estimation and hence re-imputation. The convergence could be shown empirically to be accelerated by such a procedure. The rate improvement would be from the difference between $T$ and $\psi_0 T$, where $\psi_0$ is the lowest proportion of observation among all entries from Assumption (O1). We omit the lengthy proofs as eventually the rates only differ by a constant, but we note here that re-imputation can indeed improve our imputation, which is essentially credited to the more observations used when we have an initially good imputation.
\end{remark}

\subsection{Consistency: factor loadings and imputed values}\label{subsec:consistency}
We present consistency results in this section. For $k\in[K],j\in[d_k]$, define
\begin{align}
\H_{k,j} &:=\wh\D_k^{-1}\sum_{i=1}^{d_k} \wh\Q_{k,i\cdot}
    \sum_{h=1}^{d_{\text{-}k}}
        \frac{1}{|\psi_{k,ij,h}|}
        \sum_{t\in\psi_{k,ij,h}}
        \Big(
    \sum_{m=1}^{\rmk}\Lambda_{k,hm} \mat{k}{\cF_{Z,t}}_{\cdot m}
    \Big)'\Q_{k,i\cdot}\Big(
    \sum_{m=1}^{\rmk}\Lambda_{k,hm} \mat{k}{\cF_{Z,t}}_{\cdot m}
    \Big)', \label{eqn:H_j} \\
\H_k^{a} &:= \frac{1}{T}\sum_{t=1}^T\wh\D_k^{-1}\wh\Q_k'\Q_k\mat{k}{\cF_{Z,t}}\bLambda_k'\bLambda_k\mat{k}{\cF_{Z,t}}',
\label{eqn:H_k(all)}
\end{align}
where $\wh\D_k := \wh\Q_k'\wh\S_k\wh\Q_k$ is a diagonal matrix of eigenvalues of $\wh\S_k$ defined in (\ref{eqn: weighted-Sk-entry}). Hence $\H_{k,j} = \H_k^a$ if there are no missing entries, i.e., $|\psi_{k,ij,h}| = T$ for each $k\in[K], i,j\in[d_k]$ and $h\in[\dmk]$. Furthermore, each $\H_{k,j}$ and $\H_k^a$ can be shown asymptotically invertible (see Lemma \ref{lemma:3} and \ref{lemma:limit} in the supplementary materials).

{
We first present a consistency result for the factor loading matrix estimator $\wh\Q_k$ of $\Q_k$. In particular, our theoretical rates are shown in the presence of potential weak factors. To compare with results in similar literature, we will end this section with a simplified result. Readers interested in the rates under only pervasive factors can go straight to Corollary \ref{corollary:simplified_consistency}.
}

\begin{theorem}\label{thm:factor_loading_consistency}
Under Assumptions (O1), (M1), (F1), (L1), (E1), (E2) and (R1), for any $k\in[K]$, we have
\begin{equation*}
    \frac{1}{d_k}
    \sum_{j=1}^{d_k}\Big\|
    \wh\Q_{k,j\cdot}-\H_{k,j}
    \Q_{k,j\cdot}
    \Big\|^2 =
    O_P\Bigg(d_k^{2(\alpha_{k,1} - \alpha_{k,r_k})-1}
    \bigg(\frac{1}{T\dmk}+\frac{1}{d_k}\bigg)
    \frac{d^2}{g_s^2} \Bigg) = o_P(1),
\end{equation*}
where $g_s$ is defined in Assumption (R1). Furthermore,
for each $i,j\in[d_k]$ and $h\in[\dmk]$, if the set $\psi_{k,ij,h}$ defined in (\ref{eqn:set_psi}) has cardinality satisfying
\[\frac{|\psi_{k,ij,h}|}{T} =1-\eta_{k,ij,h} \geq 1-\eta = \psi_0,\]
then we have
\begin{align*}
    \frac{1}{d_k}
    \sum_{j=1}^{d_k}\Big\|
    \wh{\bf{Q}}_{k,j\cdot}-\bf{H}_k^a
    \bf{Q}_{k,j\cdot}
    \Big\|^2 &= \frac{1}{d_k}\norm{\wh\Q_k - \Q_k\H_k^{a'}}_F^2\\
    &= O_P\Bigg(d_k^{2(\alpha_{k,1} - \alpha_{k,r_k})-1}
    \Bigg[
    \bigg(\frac{1}{T\dmk}+\frac{1}{d_k}
    \bigg)\frac{d^2}{g_s^2}
    + \min\bigg(\frac{1}{T}, \frac{\eta^2}{(1-\eta)^2}\bigg)\Bigg]\Bigg) = o_P(1).
\end{align*}
\end{theorem}
The proof of the theorem can be found in the supplementary materials of this paper. The two results in Theorem \ref{thm:factor_loading_consistency} coincide with each other if $\eta=0$, i.e., there are no missing values.

We present the two results in the theorem to highlight the difficulty of obtaining consistency when there are missing values. Since a factor loading matrix is not uniquely defined, in the second result in Theorem 1 we are estimating how close $\wh\Q_k$ is to a version of $\Q_k$ in Frobenius norm, namely $\Q_k\H_k^a$, which is still defining the same factor loading space as $\Q_k$ does. With missing data, such a feat is complicated, in the sense that for the $j$-th row of $\Q_k$, $\wh\Q_{k,j\cdot}$, there corresponds an $\H_{k,j}$  different from $\H_k^a$ in general, so that $\wh\Q_{k,j\cdot}$ is close to $\H_{k,j}\Q_{k,j\cdot}$. The extra rate $\min(1/T, \eta^2/(1-\eta)^2)$ in the second result is essentially measuring how close each  $\H_{k,j}$ is to $\H_k^a$. See Lemma \ref{lemma:3} in the supplementary materials as well.

\begin{theorem}\label{thm:imputation_consistency}
Under the Assumptions in Theorem \ref{thm:factor_loading_consistency}, suppose we further have
$d_k^{2\alpha_{k,1} - 3\alpha_{k,r_k}} = o(\dmk)$. Define
\[g_{\eta} := \min\bigg(\frac{1}{T}, \frac{\eta^2}{(1-\eta)^2}\bigg), \;\;\; g_w := \prod_{k=1}^Kd_k^{\alpha_{k,r_k}}.\]
Then we have the following.

\noindent 1. The error of the estimated factor series has rate
\begin{equation*}
\begin{split}
    &\Big\|\textnormal{\textbf{vec}}
    (\wh{\cF}_{Z,t}) -
    \big(\H_K^{a'}\otimes\dots\otimes\H_1^{a'}\big)^{-1}
    \textnormal{\textbf{vec}}(\cF_{Z,t})
    \Big\|^2\\
    &=
    O_P\Bigg(\max_{k\in[K]}\Bigg\{
     T^{-1}d
    d_k^{3\alpha_{k,1} - 2\alpha_{k,r_k}} g_s^{-1} + d^2 g_s^{-1}d_k^{2\alpha_{k,1}- 3\alpha_{k,r_k}-1} + g_{\eta}g_s d_k^{2\alpha_{k,1}- 3\alpha_{k,r_k}+1}
    \Bigg\} + \frac{d}{g_w}\Bigg).
\end{split}
\end{equation*}
2. For any $k\in[K],i_k\in[d_k],t\in[T]$, the squared individual imputation error is
\begin{equation*}
\begin{split}
    &(\wh{\c{C}}_{t,i_1,\dots,i_K}-
    \c{C}_{t,i_1,\dots,i_K})^2 \\
    &=  O_P\Bigg(\max_{k\in[K]}\Bigg\{
     T^{-1}d
    d_k^{3\alpha_{k,1} - 2\alpha_{k,r_k}} g_s^{-1}g_w^{-1} + d^2 g_s^{-1}g_w^{-1}d_k^{2\alpha_{k,1}- 3\alpha_{k,r_k}-1} + g_{\eta}g_sg_w^{-1} d_k^{2\alpha_{k,1}- 3\alpha_{k,r_k}+1}
    \Bigg\} + \frac{d}{g_w^2}\Bigg).
\end{split}
\end{equation*}
3. The average imputation error is given by
\begin{equation*}
\begin{split}
    &\hspace{5mm}
    \frac{1}{Td} \sum_{t=1}^T
    \sum_{i_1,\dots,i_K=1}^{d_1,\dots,d_K}
    (\wh{\c{C}}_{t,i_1,\dots,i_K}-
    \c{C}_{t,i_1,\dots,i_K})^2 \\
    &= O_P\Bigg(\max_{k\in[K]}\Bigg\{
     T^{-1}
    d_k^{3\alpha_{k,1} - 2\alpha_{k,r_k}} g_s^{-1} + dg_s^{-1}d_k^{2\alpha_{k,1}- 3\alpha_{k,r_k}-1} + d^{-1}g_{\eta}g_s d_k^{2\alpha_{k,1}- 3\alpha_{k,r_k}+1}
    \Bigg\} + \frac{1}{g_w}\Bigg).
\end{split}
\end{equation*}
\end{theorem}
The proof can be found in the supplementary materials, which utilizes some rates from the proof of Theorem \ref{thm:asymp_normality_loadings} in the supplementary materials (without the need for extra rate restrictions like Theorem \ref{thm:asymp_normality_loadings} though). The complication of missing data comes explicitly from the rate $g_{\eta}$. The average squared imputation error in result 3 improves upon individual squared error in result 2 when weak factors exist, with degree of improvements larger when the difference in strength of factors is larger.

{
Our rate can be considered a generalization to a general order tensor, with general factor strengths and missing data, see the comparison of our results with others' below Corollary \ref{corollary:simplified_consistency}. Such generalizations have certainly revealed that when there are weak factors, especially when the strongest and weakest factor strengths are quite different, those rates of convergence greatly suffer.
}

\begin{corollary}
\label{corollary:simplified_consistency}
(Simplified Theorem \ref{thm:factor_loading_consistency} and \ref{thm:imputation_consistency} under pervasive factors.)
Let Assumption (O1), (M1), (F1), (L1), (E1) and (E2) hold. If all factors are pervasive such that $\alpha_{k,j}=1$ for all $k\in[K], j\in[r_k]$, then with the renormalised loading and core factor estimators defined as $\wh\A_k = \sqrt{d_k} \, \wh\Q_k$ and $\wh\cF_t = \wh\cF_{Z,t} /\sqrt{d}$, we have the following:

\noindent 1. The (renormalized) loading estimator is consistent such that for any $k\in[K]$,
\begin{align*}
    \frac{1}{d_k} \sum_{j=1}^{d_k} \Big\| \wh\A_{k,j\cdot} -\H_{k,j} \A_{k,j\cdot} \Big\|^2
    &= O_P\Big( \frac{1}{T\dmk}+\frac{1}{d_k}\Big) = o_P(1) ,\\
    \frac{1}{d_k} \sum_{j=1}^{d_k}\Big\| \wh{\A}_{k,j\cdot} -\bf{H}_k^a \A_{k,j\cdot} \Big\|^2
    &= O_P\Big\{ \frac{1}{T\dmk}+\frac{1}{d_k}
    + \min\Big(\frac{1}{T}, \frac{\eta^2}{(1-\eta)^2}\Big) \Big\} = o_P(1).
\end{align*}
2. The (renormalized) core factor estimator is consistent such that for any $t\in[T]$,
\[
\Big\|\textnormal{\textbf{vec}} (\wh{\cF}_{t}) - \big(\H_K^{a'} \otimes \dots \otimes \H_1^{a'}\big)^{-1} \textnormal{\textbf{vec}}(\cF_{t}) \Big\|^2
=
O_P\Big\{ \max_{k\in[K]}\Big( \frac{1}{T \dmk} + \frac{1}{d_k^2} \Big) + \min\Big(\frac{1}{T}, \frac{\eta^2}{(1-\eta)^2}\Big) \Big\}.
\]
3. The imputation is consistent both for each entry and on average (with the same rate), such that for any $k\in[K],i_k\in[d_k],t\in[T]$,
\[
(\wh{\c{C}}_{t, i_1, \dots, i_K} -\c{C}_{t, i_1, \dots, i_K})^2
=
O_P\Big\{ \max_{k\in[K]}\Big( \frac{1}{T \dmk} + \frac{1}{d_k^2} \Big) + \min\Big(\frac{1}{T}, \frac{\eta^2}{(1-\eta)^2}\Big) + \frac{1}{d} \Big\}.
\]
\end{corollary}

{
When $K=1$ with missing data, result 1 has rate $1/\min(d_1,T)$, which is the same as the rate in Theorem 1 of \cite{Xiong_Pelger}. If $K=2$ and $\eta=0$ (i.e, no missing values), result 1 has rate  $1/\min(d_k,T\dmk)$, which is consistent with Theorem 1 of \cite{ChenFan2023}, for example. For a general order-$K$ tensor without missing data (i.e., $\eta=0$), our Lemma \ref{lemma:5_additional} in the supplementary materials states that
\[ \|\wh\Q_{k,j\cdot} - \H_k^a\Q_{k,j\cdot}\|^2 = O_P\Big(\frac{1}{Td} + \frac{1}{d_k^3}\Big), \;
\text{ implying } \; \frac{1}{d_k}\|\wh\A_k - \A_k\H_k^a\|_F^2 = O_P\Big(\frac{1}{T\dmk} + \frac{1}{d_k^2}\Big), \]
which aligns with Theorem 3.1 of \cite{Heetal2022a} or \cite{Barigozzietal2023b}.
}

{
If $K\geq 2$ and $\eta=0$, result 3 has rate
\[\max_{k\in[K]}\bigg(\frac{1}{T\dmk} + \frac{1}{d_k^2}\bigg) + \frac{1}{d} \asymp \frac{1}{\min(Td_{\text{-}1}, \ldots, Td_{\text{-}K}, d_1^2,\ldots,d_K^2)}.\]
This rate is the same as the result in Theorem 4 of \cite{ChenFan2023} for $K=2$, which is a rate for estimating the common component. On the other hand, if $\eta$ is a constant and $K=1$, then result 3 becomes $d_1^{-1} + T^{-1}\asymp 1/\min(d_1,T)$, which is the same rate as result 3 of Theorem 2 in \cite{Xiong_Pelger}.
}

\subsection{Inference on the factor loadings}\label{subsec:asymptotic_notmality}
We establish asymptotic normality of the factor loadings for inference purpose. In Section \ref{subsec:covariance_estimation} we present the covariance matrix estimator for practical use of our asymptotic normality result. First, we define
\begin{equation}\label{eqn:Hka*}
\H_k^{a,\ast} :=  \textnormal{tr}(\Amk'\Amk))^{1/2}\cdot\D_k^{-1/2}\Upsilon_k'\Z_k^{1/2},
\end{equation}
where $\D_k := \tr(\Amk'\Amk)\, \diag\{\lambda_1(\A_k'\A_k),\ldots, \lambda_{r_k}(\A_k'\A_k)\}$, and $\Upsilon_k$ is the eigenvector matrix of $\textnormal{tr}(\Amk'\Amk) \cdot g_s^{-1}d_k^{\alpha_{k,1}-\alpha_{k,r_k}}\, \Z_k^{1/2} \bSigma_{A,k}\Z_k^{1/2}$. It turns out $\H_k^{a,\ast}$ is the probability limit of $\H_k^a$ defined in (\ref{eqn:H_k(all)}). Before presenting our results, we need three additional assumptions.

\begin{itemize}
\item[(L2)] (Eigenvalues)
\textit{For any $k\in[K]$, the eigenvalues of the $r_k\times r_k$ matrix $\bf{\Sigma}_{A,k}\Z_k$ from Assumption (L1) are distinct.
}

\item[(AD1)]
\textit{Define $\omega_B := \dmk^{-1} d_k^{2\alpha_{k,r_k}-3\alpha_{k,1}}g_s^2$ and the following,
\[
\bf\Xi_{k,j} :=
    \plim_{T,d_1,\dots,d_K\to\infty}
    \textnormal{Var}\bigg(
    \sum_{i=1}^{d_k} \Q_{k,i\cdot}
    \sum_{h=1}^{d_{\text{-}k}}
    \frac{1}{|\psi_{k,ij,h}|}
    \sum_{t\in\psi_{k,ij,h}}
    \textnormal{mat}_k(\cE_t)_{jh}
    (\Amk)_{h\cdot}'
    \textnormal{mat}_k(\cF_t)'
    \A_{k,i\cdot}\bigg),
\]
then we assume $T\omega_B\cdot \big\|\D_k^{-1} \H_k^{a,\ast} \bf\Xi_{k,j} (\H_k^{a,\ast})' \D_k^{-1} \big\|_F$ is of constant order.
}

\item[(AD2)]
\textit{Define the filtration $\c{G}^T:=\sigma(\cup_{s=1}^T\c{G}_s)$ with $\c{G}_s:=\sigma(\{\c{M}_{t,i_1, \dots,i_K}\mid t\leq s\}, \A_1, \dots,\A_K)$, and
\begin{align*}
    \Delta_{F,k,ij,h} :=
    \frac{1}{|\psi_{k,ij,h}|}
    \sum_{t\in\psi_{k,ij,h}}
    \textnormal{mat}_k(\cF_t)
    \v_{k,h}\v_{k,h}'\textnormal{mat}_k(\cF_t)'
    - \frac{1}{T}\sum_{t=1}^T \textnormal{mat}_k(\cF_t)
    \v_{k,h}\v_{k,h}'\textnormal{mat}_k(\cF_t)' ,
\end{align*}
where $\v_{k,h} := [\otimes_{l\in[K]\setminus \{k\}}\A_l]_{h\cdot}$. With $\Q_k$ being the normalised mode-$k$ factor loading defined below Assumption (L1), we have for every $k\in[K], j\in[d_k]$, for a function $h_{k,j}: \mathbb{R}^{r_k} \rightarrow \mathbb{R}^{r_k\times r_k}$,
\begin{align*}
\sqrt{T d_k^{\alpha_{k,r_k}} }&\cdot \D_k^{-1}\H_k^{a,\ast}\sum_{i=1}^{d_k} \Q_{k,i\cdot} \A_{k,i\cdot}' \sum_{h=1}^{\dmk}
\Delta_{F,k,ij,h} \A_{k,j\cdot}\\
&\to \c{N}(\0, \D_k^{-1}\H_k^{a,\ast} h_{k,j}(\A_{k,j\cdot})(\H_k^{a,\ast})' \D_k^{-1} )
\;\;\;
\text{$\c{G}^T$-stably}.
\end{align*}
}
\end{itemize}
Assumption (AD1) guarantees a part of the covariance matrix of the asymptotic normality in Theorem \ref{thm:asymp_normality_loadings} is of constant order. It can be regarded as a lower bound condition which is necessary for the dominance of a certain term involved in the asymptotic normality. Since we show the upper bound of $T\omega_B\cdot \big\|\D_k^{-1} \H_k^{a,\ast} \bf\Xi_{k,j} (\H_k^{a,\ast})' \D_k^{-1} \big\|_F$ is of constant order in the proof of Theorem \ref{thm:asymp_normality_loadings} in the supplementary materials, this assumption is not particularly strong.

Assumption (AD2) is required since the missing data creates a discrepancy term $\Delta_{F,k,ij,h}$ as defined in the assumption. This assumption is also parallel to Assumption G3.5 in \cite{Xiong_Pelger}. We demonstrate how this assumption is satisfied with Assumption (O1), (F1), (L1) and two additional but simpler assumptions in Proposition \ref{Prop:example_AD2} in Section \ref{subsec:AD2_simplify}.

\begin{theorem}\label{thm:asymp_normality_loadings}
Let all the assumptions under Theorem \ref{thm:imputation_consistency} hold, in addition to (L2), (AD1) and (AD2) above. With $r_k$ fixed and $d_k, T \rightarrow \infty$ for $k\in[K]$, suppose also $T\dmk = o(d_k^{\alpha_{k,1} + \alpha_{k,r_k}})$. We have
\[
\sqrt{T d_k^{\alpha_{k,r_k}} } \cdot
(\wh\Q_{k,j\cdot}-\H_k^a\Q_{k,j\cdot})
\xrightarrow{\c{D}}
\cN(\0, \D_k^{-1}\H_k^{a,\ast}
(Td_k^{\alpha_{k,r_k}} \cdot\bf\Xi_{k,j} + h_{k,j}(\A_{j\cdot}))(\H_k^{a,\ast})'\D_k^{-1}).
\]
Furthermore, if $Td^{-1}g_s^2g_\eta d_k^{1 + \alpha_{k,1}-3\alpha_{k,r_k}} = o(1)$ is also satisfied, then
\begin{align*}
\sqrt{T\omega_B}&\cdot
(\wh\Q_{k,j\cdot}-\H_k^a \Q_{k,j\cdot})
\xrightarrow{\c{D}}
\cN(\0, T\omega_B \cdot
\D_k^{-1}\H_k^{a,\ast}
\bf\Xi_{k,j}(\H_k^{a,\ast})'\D_k^{-1}).
\end{align*}
\end{theorem}
If all factors are pervasive, the rate condition $T\dmk = o(d_k^{\alpha_{k,1} + \alpha_{k,r_k}})$ reduces to $T\dmk = o(d_k^2)$, which is equivalent to the condition needed for asymptotic normality in \cite{Bai2003} for $K=1$ and \cite{ChenFan2023} for $K=2$. The first asymptotic normality result is compatible to Theorem 2.1 of \cite{Xiong_Pelger} when all factors are pervasive. In their Theorem 2.1, the $\Gamma_{\Lambda, j}^{obs}$ is in fact of rate $N^{-1}$, so that the normalizing rate is $\sqrt{TN}$, which is exactly $\sqrt{Td_1}$ in our first result when $K=1$.

Suppose all factors are pervasive. The rate condition $Td^{-1}g_s^2g_\eta d_k^{1 + \alpha_{k,1}-3\alpha_{k,r_k}} = o(1)$ is automatically satisfied when there is no missing data, i.e., $\eta=0$ so that $g_\eta = 0$. If so, the rate of convergence is $\sqrt{T\omega_B} = \sqrt{Td}$,
{
which is compatible to Theorem 2.1, Theorem 2.2 of \cite{ChenFan2023} and Theorem 3.2 of \cite{Barigozzietal2023b} (after our normalization to their factor loading matrices).
}
The condition is also satisfied when there is only finite number of missing data, so that $\eta \asymp T^{-1}$ and $g_\eta \asymp T^{-2}$, and $d_1, d_2 = o(T)$ for $K=2$.

\begin{remark}\label{remark:2}
We do not establish asymptotic normality for the estimated factor series and common components. The reason is that for tensor with $K>1$, the decomposition in the estimated factor series and the common components cannot be dominated by terms that are asymptotically normal. This is also the reason why \cite{ChenFan2023} does not include asymptotic normality for the estimated factor series and common components.
{
\cite{Barigozzietal2023b} constructs asymptotic normality
for the core factor built upon their projection estimator $\wt{\mathcal{F}}_t$, which is sensible as the projecting loading estimator already has an improved rate. In comparisons, the rate of any least-square-type estimators, such as the one in \cite{ChenFan2023} for matrix data and the one in our case for general tensors, is insufficient for a potentially asymptotically Gaussian term to be dominating.
The main goal of this work is to impute missing entries, and existing methods on tensor factor models using Tucker decomposition should be able to be applied with all missing entries replaced by the consistent imputations.
}
\end{remark}

\subsection{Estimation of covariance matrix}\label{subsec:covariance_estimation}
In order to carry out inferences for the factor loadings using Theorem \ref{thm:asymp_normality_loadings}, we need to estimate the asymptotic covariance matrix for $\wh\Q_{k,j\cdot} - \H_k^a\Q_{k,j\cdot}$. To this end, we use the heteroscedasticity and autocorrelation consistent (HAC) estimators \citep{HAC_1987} based on $\{\wh\Q_k, \mat{k}{\wh\cC_t}, \mat{k}{\wh\cE_t}\}_{t\in[T]}$, where
\begin{align*}
    \mat{k}{\wh\cC_t} &=(\wh\Q_k)\mat{k}{\wh\cF_{Z,t}}(\wh\Q_K\otimes\cdots\otimes\wh\Q_{k+1}\otimes\wh\Q_{k-1}\otimes\cdots\otimes\wh\Q_1)'
    ,\;\; \mat{k}{\wh\cE_t} :=
    \mat{k}{\cY_t}-\mat{k}{\wh\cC_t}.
\end{align*}
With a tuning parameter $\beta$ such that $\beta \rightarrow \infty$ and $\beta/(Td_k^{\alpha_{k,r_k}})^{1/4} \rightarrow 0$, we define two HAC estimators
\begin{align*}
    \wh\bSigma_{HAC} &:=
    \D_{k,0,j} + \sum_{\nu=1}^\beta
    \bigg(1- \frac{\nu}{1+\beta}\bigg)
    \Big( \D_{k,\nu,j}+
    \D_{k,\nu,j}' \Big),\\
   \wh\bSigma_{HAC}^\Delta &:=
    \D_{k,0,j}^\Delta + \sum_{\nu=1}^\beta
    \bigg(1- \frac{\nu}{1+\beta}\bigg)
    \Big( \D_{k,\nu,j}^\Delta+
    (\D_{k,\nu,j}^\Delta)' \Big), \;\;\; \text{where}\\
        \D_{k,\nu,j} &:=
    \sum_{t=1+\nu}^T\bigg(
    \sum_{i=1}^{d_k}
    \Big(\frac{1}{T}
    \sum_{s=1}^T \wh\D_k^{-1}
    \wh\Q_k'
    \wh\C_{(k),s}
    \wh\C_{(k),s,i\cdot}\Big)
    \sum_{h=1}^{\dmk}
    \frac{1}{|\psi_{k,ij,h}|}
    \wh\E_{(k),t,jh} \wh\C_{(k),t,ih} \cdot \b{1}\{t\in\psi_{k,ij,h}\} \bigg) \\
    &\cdot \bigg(
    \sum_{i=1}^{d_k}
    \Big(\frac{1}{T}
    \sum_{s=1}^T \wh\D_k^{-1}
    \wh\Q_k'
    \wh\C_{(k),s}
    \wh\C_{(k),s,i\cdot}\Big)
    \sum_{h=1}^{\dmk}
    \frac{1}{|\psi_{k,ij,h}|}
    \wh\E_{(k),t-\nu,jh} \wh\C_{(k),t-\nu,ih} \cdot \b{1}\{t-\nu\in\psi_{k,ij,h}\} \bigg)',\\
        \D_{k,\nu,j}^\Delta &:=
    \sum_{t=1+\nu}^T\Bigg[
    \sum_{i=1}^{d_k}
    \Big(\frac{1}{T}
    \sum_{s=1}^T \wh\D_k^{-1}
    \wh\Q_k'
    \wh\C_{(k),s}
    \wh\C_{(k),s,i\cdot}\Big)
    \sum_{h=1}^{\dmk}\bigg(
    \frac{1}{|\psi_{k,ij,h}|}
    \wh\C_{(k),t,ih} \wh\C_{(k),t,jh} \cdot \b{1}\{t\in\psi_{k,ij,h}\} \\
    &-
    \frac{1}{T}\wh\C_{(k),t,ih} \wh\C_{(k),t,jh} \bigg) \Bigg] \cdot \Bigg[
    \sum_{i=1}^{d_k} \Big(\frac{1}{T}
    \sum_{s=1}^T \wh\D_k^{-1}
    \wh\Q_k'
    \wh\C_{(k),s}
    \wh\C_{(k),s,i\cdot}\Big) \\
    &\cdot
    \sum_{h=1}^{\dmk} \bigg(
    \frac{1}{|\psi_{k,ij,h}|}
    \wh\C_{(k),t-\nu,ih} \wh\C_{(k),t-\nu,jh} \cdot \b{1}\{t-\nu\in\psi_{k,ij,h}\}
    - \frac{1}{T} \wh\C_{(k),t-\nu,ih} \wh\C_{(k),t-\nu,jh}
    \bigg) \Bigg]',
\end{align*}
where $\wh\C_{(k),s} := \mat{k}{\wh\cC_s}$ and $\wh\E_{(k),s} := \mat{k}{\wh\cE_s}$.

\begin{theorem}\label{thm:covariance_estimator}
  Let all the assumptions under Theorem \ref{thm:imputation_consistency} hold, in addition to (L2), (AD1) and (AD2) above. With $r_k$ fixed and $d_k, T \rightarrow \infty$ for $k\in[K]$, suppose also the rate for the individual common component imputation error in result 2 of Theorem \ref{thm:imputation_consistency} is $o(1)$, together with $T\dmk = o(d_k^{\alpha_{k,1} + \alpha_{k,r_k}})$ and $d_k^{2(\alpha_{k,1} - \alpha_{k,r_k})}[(T\dmk)^{-1} + d_k^{-1}]d^2g_s^{-2} = o(1)$. Then
  \begin{itemize}
    \item[1.] $\wh\D_k^{-1}\wh\bSigma_{HAC}\wh\D_k^{-1}$ is consistent for $\D_k^{-1}\H_k^{a,*}\bf{\Xi}_{k,j}(\H_k^{a,*})'\D_k^{-1}$;

    \item[2.] $\wh\D_k^{-1}\wh\bSigma_{HAC}^{\Delta}\wh\D_k^{-1}$ is consistent for $(Td_k^{\alpha_{k,r_k}})^{-1}\D_k^{-1}\H_k^{a,*}h_{k,j}(\A_{k,j\cdot})(\H_{k}^{a,*})'\D_k^{-1}$;

    \item[3.] $(\wh\bSigma_{HAC} + \wh\bSigma_{HAC}^{\Delta})^{-1/2}\wh\D_k(\wh\Q_{k,j\cdot} - \H_k^a\Q_{k,j\cdot}) \xrightarrow{\c{D}}
\cN(\0, \I_{r_k})$.
  \end{itemize}
\end{theorem}
The extra rate assumption $d_k^{2(\alpha_{k,1} - \alpha_{k,r_k})}[(T\dmk)^{-1} + d_k^{-1}]d^2g_s^{-2} =o(1)$ makes sure that we have Frobenius norm consistency for $\wh\Q_k$ from Theorem \ref{thm:factor_loading_consistency}.
The imputation error from result 2 of Theorem \ref{thm:imputation_consistency} also has rate going to 0 when
are all factors are pervasive, for instance.
With result 3 in particular, we can perform inferences on any rows of $\wh\Q_{k}$. Practical performances of result 3 is demonstrated in Section \ref{subsec:empirical_asymp_normality}. The reason that we need two HAC estimators is that similar to Theorem \ref{thm:factor_loading_consistency}, there is a component for missing data, arising from the fact that $\H_{k,j}$ is different from $\H_k^a$ for each $j\in[d_k]$ in general.

\subsection{Estimation of number of factors}\label{subsec:estimate_number_of_factors}
The reconstructed mode-$k$ sample covariance matrix $\wh\S_k$ is in fact estimating a complete-sample version of a matrix $\R_k^*$, where
\begin{equation}\label{eqn:R*}
\R_k^* := \frac{1}{T}\sum_{t=1}^T \Q_k\mat{k}{\cF_{Z,t}}\bLambda_k'\bLambda_k\mat{k}{\cF_{Z,t}}'\Q_k',
\end{equation}
and $\cF_{Z,t}$ and $\bLambda_k$ are defined in (\ref{eqn: tenfac-all2}). It turns out that we have $\lambda_j(\wh{\S}_k) \asymp_P \lambda_j(\R_k^*)$ for $j\in[r_k]$, and
\[\lambda_j(\R_k^*) \asymp_P d_k^{\alpha_{k,j} - \alpha_{k,1}}g_s, \;\;\; g_s := \prod_{k=1}^Kd_k^{\alpha_{k,1}} \text{ as defined in (R1).}\]
We have the following theorem.
\begin{theorem}\label{thm:number_of_factors}
Let Assumption (O1), (M1), (F1), (L1), (E1), (E2) and (R1) hold. Moreover, assume
\[  \left\{
      \begin{array}{ll}
        dg_s^{-1}d_k^{\alpha_{k,1}-\alpha_{k,r_k}}[(T\dmk)^{-1/2} + d_k^{-1/2}]  = o(d_k^{\alpha_{k,j+1} - \alpha_{k,j}}), & \hbox{$j\in[r_k-1]$ with $r_k\geq 2$;} \\
        dg_s^{-1}[(T\dmk)^{-1/2} + d_k^{-1/2}] = o(1), & \hbox{$r_k=1$.}
      \end{array}
    \right.
\]
Then $\wh{r}_k$ is a consistent estimator of $r_k$, where
\begin{equation}\label{eqn:r_k_hat}
\wh{r}_k := \arg\min_\ell \bigg\{\frac{\lambda_{\ell+1}(\wh\S_k) + \xi}{\lambda_{\ell}(\wh\S_k) + \xi}, \;\ell\in[\lfloor d_k/2 \rfloor]\bigg\}, \;\;\; \xi \asymp d[(T\dmk)^{-1/2} + d_k^{-1/2}].
\end{equation}
\end{theorem}
The extra rate assumption is satisfied, for instance, when all factors corresponding to $\A_k$ are pervasive. An eigenvalue-ratio estimator is considered in \cite{LamYao2012} and \cite{AhnHorenstein2013}, while a perturbed eigenvalue ratio estimator is considered in \cite{Pelger2019}. However, all of these estimators are for a vector time series factor model. Our estimator $\wh{r}_k$ in (\ref{eqn:r_k_hat}) extracts eigenvalues from $\wh{\S}_k$, which is not necessarily positive semi-definite. The addition of $\xi$ can make $\wh\S_k + \xi\I_{d_k}$ positive semi-definite, while stabilizing the estimator. We naturally assume that $r_k < d_k/2$, which is a very reasonable assumption for all applications of factor models. In fact, our recommended choice of $\xi$ is
\[\xi = \frac{1}{5}d[(T\dmk)^{-1/2} + d_k^{-1/2}].\]
The requirement that $\xi \asymp d[(T\dmk)^{-1/2} + d_k^{-1/2}]$ ensures that $\xi = o_P(\lambda_{r_k}(\wh\S_k))$ from our rate assumption in the theorem.
Our simulations in Section \ref{subsec:performance_of_number_of_factors} suggest that this proposal works very well.

\subsection{*How Assumption (AD2) can be implied}\label{subsec:AD2_simplify}
This section details how Assumption (AD2) can be implied from simpler assumptions. Readers can skip this part and go straight to the next section for a more integral reading experience.
We begin by presenting a proposition.
\begin{proposition}\label{Prop:example_AD2}
Let Assumption (O1), (F1), (L1) hold. For a given $k\in[K], j\in[d_k]$, assume also the following:
\begin{itemize}
    \item [1.] The mode-$k$ factor is strong enough such that $\alpha_{k,r_k}>4/5$, and $d_k^{\alpha_{k,1}-\alpha_{k,r_k}} T^{-\epsilon/2} = o(1)$ with some $\epsilon\in(0,1)$.
    \item [2.] There exists some $\psi_{k,ij}$ such that $\psi_{k,ij}=\psi_{k,ij,h}$ for any $i\in[d_k], h\in[\dmk]$. Furthermore, there exists $\omega_{\psi,k,j}$ such that
    \[
    d_k^{-2}\sum_{i=1}^{d_k}\sum_{l=1}^{d_k}
    \Big(\frac{T \cdot \b{1}\{t\in\psi_{k,ij}\} }{|\psi_{k,ij}|} -1 \Big) \Big(\frac{T \cdot \b{1}\{t\in\psi_{k,lj}\} }{|\psi_{k,lj}|} -1 \Big)
    \xrightarrow{p} \omega_{\psi,k,j} .
    \]
\end{itemize}
With the above, Assumption (AD2) is satisfied.
\end{proposition}
Condition 1 and 2 in Proposition \ref{Prop:example_AD2} are on the factor strength and missingness pattern, respectively. Condition 1 is trivially satisfied if all factors are pervasive. Condition 2 can be easily satisfied by assuming that in $\mat{k}{\cY_t}$, all the elements in each row are simultaneously missing with probability $1-p_0$. We then have
\begin{equation*}
\begin{split}
    &\hspace{5mm}
    d_k^{-2} \sum_{i=1}^{d_k}\sum_{l=1}^{d_k}
    \Big(\frac{T \cdot \b{1}\{t\in\psi_{k,ij}\} }{|\psi_{k,ij}|} -1 \Big) \Big(\frac{T \cdot \b{1}\{t\in\psi_{k,lj}\} }{|\psi_{k,lj}|} -1 \Big) \\
    &=
    d_k^{-2} \sum_{i=1}^{d_k}\sum_{l=1}^{d_k}
    \Big(\frac{T^2 \cdot \b{1}\{t\in\psi_{k,ij}\} \cdot
    \b{1}\{t\in\psi_{k,lj}\}}{|\psi_{k,ij}|
    \cdot |\psi_{k,lj}|} -
    \frac{T \cdot \b{1}\{t\in\psi_{k,ij}\} }{|\psi_{k,ij}|}
    - \frac{T \cdot \b{1}\{t\in\psi_{k,lj}\} }{|\psi_{k,lj}|} + 1 \Big) \\
    &\xrightarrow{p}
    p_0^{-1} - 1,
\end{split}
\end{equation*}
which is $\omega_{\psi, k,j}$. Similar to Assumption S3.2 in \cite{Xiong_Pelger}, the value of $\omega_{\psi,k,j}$ can be regarded as a measure of missingness complexity. It is a parameter related to the variance of the stable convergence, and tends to increase when there is a larger portion of data missing.

\section{Numerical Results}\label{sec:numerical}
\subsection{Simulation}\label{subsec:simulation}
We demonstrate the empirical performance of our estimators in this section. Note that we do not have comparisons to other imputation methods since to the best of our knowledge, there are no other general imputation methods available for $K>1$ apart from tensor completion methods for very specific applications as mentioned in the introduction.
{
However, we will make comparison with an alternative approach to impute tensor time series combining \cite{Xiong_Pelger} and \cite{Chen_Lam}, as demonstrated in Section \ref{subsec:compare_iter_vec}.
}
Under different missing patterns which will be described later, we investigate the performance of the factor loading matrix estimators, the imputation, and the estimator of the number of factors. We also demonstrate asymptotic normality as described in Theorem \ref{thm:asymp_normality_loadings}, followed by an example plot of a statistical power function using result 3 of Theorem \ref{thm:covariance_estimator}. Throughout this section, each simulation experiment of a particular setting is repeated 1000 times.

For the data generating process, we use model (\ref{eqn: tenfac-all}) together with  Assumption (E1), (E2) and (F1). More precisely, the elements in $\cF_t$ are independent standardised AR(5) with AR coefficients 0.7, 0.3, -0.4, 0.2, and -0.1. The elements in $\cF_{e,t}$ and $\bm\epsilon_t$ are generated similary, but their AR coefficients are (-0.7, -0.3, -0.4, 0.2, 0.1) and (0.8, 0.4, -0.4, 0.2, -0.1) respectively. The standard deviation of each element in $\bm\epsilon_t$ is generated by i.i.d. $|\cN(0,1)|$.

For each $k\in[K]$, each factor loading matrix $\A_k$ is generated independently with $\A_k=\bf{U}_k\B_k$, where each entry of $\bf{U}_k\in\b{R}^{d_k\times r_k}$ is i.i.d. $\cN(0,1)$, and $\B_k\in\b{R}^{r_k\times r_k}$ is diagonal with the $j$-th diagonal entry being $d_k^{-\zeta_{k,j}}$, $0\leq \zeta_{k,j}\leq 0.5$. Pervasive (strong) factors have $\zeta_{k,j}=0$, while weak factors have $0<\zeta_{k,j}\leq 0.5$. Each entry of $\A_{e,k}\in\b{R}^{d_k\times r_{e,k}}$ is i.i.d. $\cN(0,1)$, but has independent probability of 0.95 being set exactly to 0. We set $r_{e,k}=2$ for all $k\in[K]$ throughout the section.

To investigate the performance with missing data, we consider four missing patterns:
\begin{itemize}
    \item (M-i) Random missing with probability 0.05.
    \item (M-ii) Random missing with probability 0.3.
    \item (M-iii) The missing entries have index $(t,i_1,\ldots,i_K)$, where
    \[0.5T \leq t\leq T, \;\;\; 1\leq i_k\leq 0.5d_k \text{ for all } k\in[K].\]
    \item (M-iv) Conditional random missing such that the unit with index $j$ along mode-1 is missing with probability 0.2 if $(\A_1)_{j,1} \geq 0$, and with probability 0.5 if $(\A_1)_{j,1} < 0$.
\end{itemize}

To test how robust our imputation is under heavy-tailed distribution, we consider two distributions for the innovation process in generating $\cF_t$, $\cF_{e,t}$ and $\bm\epsilon_t$: 1) i.i.d. $\cN(0,1)$; 2) i.i.d. $t_3$.

\subsubsection{Accuracy in the factor loading matrix estimators and imputations}
For both the factor loading matrix estimators and the imputations, since our procedure for vector time series ($K=1$) is essentially the same as that in \cite{Xiong_Pelger}, we show here only the performance for $K=2,3$. We use the column space distance
\[
\mathcal{D}(\Q, \wh\Q) = \Big\|\Q(\Q'\Q)^{-1}\Q' - \wh\Q(\wh\Q'\wh\Q)^{-1}\wh\Q'\Big\|
\]
for any given $\Q,\wh\Q$, which is a common used measure in the literature. For measuring the imputation accuracy, we report the relative mean squared errors (MSE) defined by
\begin{equation}\label{eqn:relative_MSE}
\text{relative MSE}_{\mathcal{S}}=
\frac{\sum_{j\in\mathcal{S}}(\wh{C}_j - C_j)^2}{\sum_{j\in\mathcal{S}}C_j^2},
\end{equation}
where $\mathcal{S}$ either denotes the set of all missing, all observed, or all available units.

We consider the following simulation settings:
\begin{itemize}
\item[(Ia)] $K=2, T=100, d_1=d_2=40, r_1=1, r_2=2$. All factors are pervasive with $\zeta_{k,j}=0$ for all $k,j$. All innovation processes in constructing $\cF_t, \cF_{e,t}$ and $\bepsilon_t$ are i.i.d. standard normal, and missing pattern is (M-i).
\item[(Ib)] Same as (Ia), but one factor is weak with $\zeta_{k,1}=0.2$ for all $k\in[K]$.
\item[(Ic)] Same as (Ia), but all innovation processes are i.i.d. $t_3$, and all factors are weak with $\zeta_{k,j}=0.2$ for all $k,j$.
\item[(Id)] Same as (Ic), but $T=200, d_1=d_2=80$.
\item[(Ie)] $K=3, T=80, d_1=d_2=d_3=20, r_1=r_2=r_3=2$. All factors are pervasive with $\zeta_{k,j}=0$ for all $k,j$. All innovation processes in constructing $\cF_t, \cF_{e,t}$ and $\bepsilon_t$ are i.i.d. standard normal, and missing pattern is (M-i).
\item[(If)] Same as (Ie), but all factors are weak with $\zeta_{k,j}=0.2$ for all $k,j$.
\item[(Ig)] Same as (If), but $T=200, d_1=d_2=d_3=40$.
\end{itemize}

\begin{figure}[t!]
\begin{center}
\centerline{\includegraphics[width=\columnwidth,scale=1]{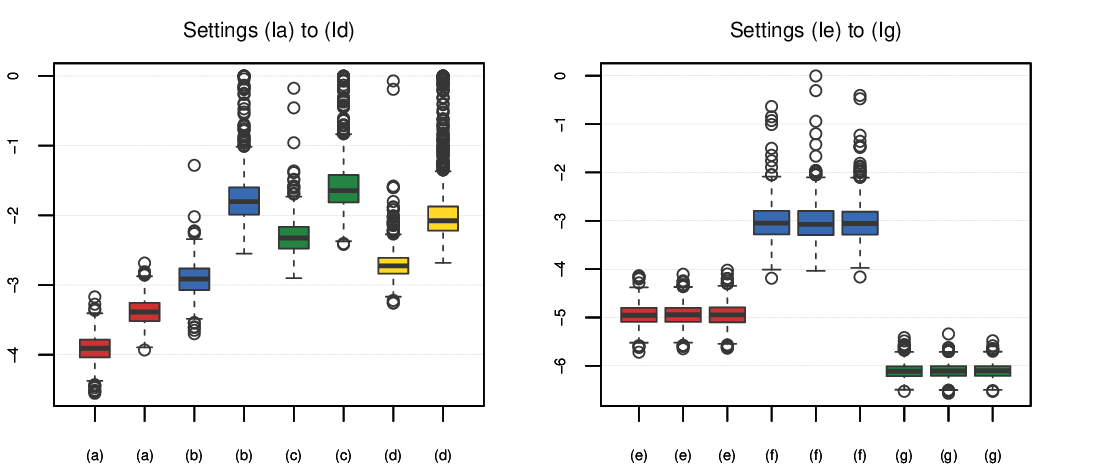}}
\caption{Plot of the column space distance $\mathcal{D}(\Q_k, \wh\Q_k)$ (in log-scale) for $k\in[K]$ for missing pattern (M-i), with $K=2$ on the left panel and $K=3$ on the right. The horizontal axis is indexed by sub-settings (a) to (g), and the $k$-th boxplot within each sub-setting corresponds to the $k$-th factor loading matrix $\Q_k$. Performance on other missing patterns are very similar, and are omitted.}
\label{Fig: sim_loading_consistency}
\end{center}
\end{figure}

Settings (Ia) to (Id) have $K=2$, and settings (Ie) to (Ig) have $K=3$. They all have missing pattern (M-i), but we have considered all settings with missing patterns (M-ii) to (M-iv), with performance of the factor loading matrix estimators very similar to those with missing pattern (M-i). Hence we are only presenting the results for settings (Ia) to (Ig) in Figure \ref{Fig: sim_loading_consistency} for the missing pattern (M-i). The imputation results for the above settings are collected in Table \ref{tab: sim_MSE_imputation}, together with those under different missing patterns.

We can see from Figure \ref{Fig: sim_loading_consistency} that the factor loading matrix estimators perform worse when there are weak factors or when the distribution of the innovation processes is fat-tailed. However, larger dimensions ameliorate the worsen performance. The increase in the loading space distance from $k=1$ to $k=2$ in settings (Ia) to (Id) is due to more factors along mode-2, which naturally incurs more errors compared to smaller $r_k$. In comparison, the loading space error shown in the right panel of Figure \ref{Fig: sim_loading_consistency} are in line for all modes due to the same number of factors along each mode.

From Table \ref{tab: sim_MSE_imputation}, we can see that missing pattern (M-iii) is uniformly more difficult in all settings for imputation. This is understandable as there is a large block of data missing in setting (M-iii), so that we obtain less information towards the ``center'' of the missing block. This is also the reason why under (M-iii), the imputation performance for the missing set is worse than the observed set, unlike for other missing patterns where all imputation performances are close.

Random missing in (M-i) and (M-ii) are relatively easier for our imputation procedure to handle.
Note that if the TALL-WIDE algorithm in \cite{Bai_Ng} were to be extended to the case for $K>1$, it can handle missing pattern (M-iii), but not (M-i) and (M-ii). The design of our method allows us to handle a wider variety of missing patterns, including random missingness.
We want to stress that we have made attempts to generalise the TALL-WIDE algorithm to impute high-order time series data for comparisons, yet the method is almost impossible to use in tensor data. The generalization is also too complicated, and hence is not shown here.


\begin{table}[t!]
\begin{center}
\begin{tabular}{|cc|cccc|ccc|}
\hline
\multicolumn{2}{|c|}{Setting}                         & \multicolumn{4}{c|}{K=2}                                                                 & \multicolumn{3}{c|}{K=3}                                     \\ \hline
\multicolumn{1}{|c|}{Missing Pattern}          & $\c{S}$    & \multicolumn{1}{c|}{(Ia)} & \multicolumn{1}{c|}{(Ib)} & \multicolumn{1}{c|}{(Ic)} & (Id) & \multicolumn{1}{c|}{(Ie)} & \multicolumn{1}{c|}{(If)} & (Ig) \\ \hline
\multicolumn{1}{|c|}{\multirow{3}{*}{(M-i)}}   & obs  & \multicolumn{1}{c|}{.002}     & \multicolumn{1}{c|}{.020}     & \multicolumn{1}{c|}{.066}     & .039      & \multicolumn{1}{c|}{2.61}     & \multicolumn{1}{c|}{120}     & .293     \\ \cline{2-9}
\multicolumn{1}{|c|}{}                         & miss & \multicolumn{1}{c|}{.002}     & \multicolumn{1}{c|}{.020}     & \multicolumn{1}{c|}{.066}     & .039     & \multicolumn{1}{c|}{2.63}     & \multicolumn{1}{c|}{121}     & .294     \\ \cline{2-9}
\multicolumn{1}{|c|}{}                         & all  & \multicolumn{1}{c|}{.002}     & \multicolumn{1}{c|}{.020}     & \multicolumn{1}{c|}{.066}     & .039     & \multicolumn{1}{c|}{2.61}     & \multicolumn{1}{c|}{120}     & .293     \\ \hline
\multicolumn{1}{|c|}{\multirow{3}{*}{(M-ii)}}  & obs  & \multicolumn{1}{c|}{.003}     & \multicolumn{1}{c|}{.025}     & \multicolumn{1}{c|}{.079}     & .045     & \multicolumn{1}{c|}{5.97}     & \multicolumn{1}{c|}{154}     &  .702    \\ \cline{2-9}
\multicolumn{1}{|c|}{}                         & miss & \multicolumn{1}{c|}{.003}     & \multicolumn{1}{c|}{.025}     & \multicolumn{1}{c|}{.079}     & .045     & \multicolumn{1}{c|}{6.06}     & \multicolumn{1}{c|}{155}     &  .703    \\ \cline{2-9}
\multicolumn{1}{|c|}{}                         & all  & \multicolumn{1}{c|}{.003}     & \multicolumn{1}{c|}{.025}     & \multicolumn{1}{c|}{.079}     & .045     & \multicolumn{1}{c|}{6.00}     & \multicolumn{1}{c|}{154}     &  .702    \\ \hline
\multicolumn{1}{|c|}{\multirow{3}{*}{(M-iii)}} & obs  & \multicolumn{1}{c|}{.004}     & \multicolumn{1}{c|}{.025}     & \multicolumn{1}{c|}{.079}     & .048      & \multicolumn{1}{c|}{6.64}     & \multicolumn{1}{c|}{136}     &  1.75    \\ \cline{2-9}
\multicolumn{1}{|c|}{}                         & miss & \multicolumn{1}{c|}{.009}     & \multicolumn{1}{c|}{.036}     & \multicolumn{1}{c|}{.107}     &  .061    & \multicolumn{1}{c|}{14.7}     & \multicolumn{1}{c|}{164}     &  4.02    \\ \cline{2-9}
\multicolumn{1}{|c|}{}                         & all  & \multicolumn{1}{c|}{.005}     & \multicolumn{1}{c|}{.026}     & \multicolumn{1}{c|}{.083}     &  .050    & \multicolumn{1}{c|}{7.19}     & \multicolumn{1}{c|}{138}     &  1.89    \\ \hline
\multicolumn{1}{|c|}{\multirow{3}{*}{(M-iv)}}  & obs  & \multicolumn{1}{c|}{.004}     & \multicolumn{1}{c|}{.027}     & \multicolumn{1}{c|}{.086}     & .047     & \multicolumn{1}{c|}{7.75}     & \multicolumn{1}{c|}{173}     & .888     \\ \cline{2-9}
\multicolumn{1}{|c|}{}                         & miss & \multicolumn{1}{c|}{.004}     & \multicolumn{1}{c|}{.028}     & \multicolumn{1}{c|}{.088}     & .047     & \multicolumn{1}{c|}{8.49}     & \multicolumn{1}{c|}{179}     & .964     \\ \cline{2-9}
\multicolumn{1}{|c|}{}                         & all  & \multicolumn{1}{c|}{.004}     & \multicolumn{1}{c|}{.027}     & \multicolumn{1}{c|}{.086}     & .047     & \multicolumn{1}{c|}{8.00}     & \multicolumn{1}{c|}{175}     & .914     \\ \hline
\end{tabular}
\end{center}
\caption {Relative MSE for settings (Ia) to (Ig), reported for $\c{S}$ as the set containing respectively observed (obs), missing (miss), and all (all) units. For $K=3$, all results presented are multiplied by $10^4$.}
\label{tab: sim_MSE_imputation}
\end{table}

\subsubsection{Performance for the estimation of the number of factors}\label{subsec:performance_of_number_of_factors}
In this section, we demonstrate the performance of our ratio estimator $\wh r_k$ in (\ref{eqn:r_k_hat}) for estimating $r_k$ for $K=1,2,3$. For each $k\in[K]$, we set the value of $\xi$ in Theorem \ref{thm:number_of_factors} as
$\xi = d[(T\dmk)^{-1/2}+d_k^{-1/2}]/5$. We have tried a wide range of values other than $1/5$ for $\xi$ in all settings, but $1/5$ is working the best in vast majority of settings,
{
see simulation results on the sensitivity of different $\xi$ in the supplementary materials.
}
Hence we do not recommend treating it as a tuning parameter in this section for saving computational time.

We present the results under a fully observed scenario and a missing data scenario for each of the following setting:
\begin{itemize}
\item[(IIa)] $K=1, T=d_1=80, r_1=2$. All factors are pervasive with $\zeta_{1,j}=0$ for all $j$. All innovation processes involved are i.i.d. standard normal. We try missing patterns (M-ii), (M-iii), and (M-iv).

\item[(IIb)] Same as (IIa), but one factor is weak with $\zeta_{1,1}=0.1$ and the missing pattern is only (M-ii).

\item[(IIc)] Same as (IIb), but factors are weak with $\zeta_{1,1}=0.1$ and $\zeta_{1,2}=0.15$.

\item[(IId)] Same as (IIc), but $T=160$.

\item[(IIIa)] $K=2, T=d_1=d_2=40, r_1=2,r_2=3$. All factors are pervasive with $\zeta_{k,j}=0$ for all $k,j$. All innovation processes involved are i.i.d. standard normal, and the missing pattern is (M-ii).

\item[(IIIb)] Same as (IIIa), but all factors are weak with $\zeta_{k,j}=0.1$ for all $k,j$.

\item[(IIIc)] Same as (IIIb), but $T=d_1=d_2=80$.

\item[(IVa)] $K=3, T=d_1=d_2=d_3=20, r_1=2, r_2=3, r_3=4$. All factors are pervasive with $\zeta_{k,j}=0$ for all $k,j$. All innovation processes involved are i.i.d. standard normal, and the missing pattern is (M-ii).

\item[(IVb)] Same as (IVa), but all innovation processes are i.i.d. $t_3$.

\item[(IVc)] Same as (IVa), but $T=40$.
\end{itemize}

Since estimating the number of factors with missing data is new to the literature, it is of interest to explore the accuracy of the estimator under different missing patterns. Hence we explore different missing patterns in setting (IIa). Extensive experiments (not shown here) on the imputation accuracy using misspecified number of factors show that underestimation is harmful, while slight overestimation hardly worsen the performance of the imputations. Thus, for each of the above settings, we also compare the performance using re-imputation and iTIP-ER by \cite{Hanetal2022}, where the re-imputation is done by using both $\wh{r}_k$ and $\wh{r}_k +1$ to avoid information loss due to underestimating the number of factors, see Table \ref{tab: sim_num_fac_IIa} and Table \ref{tab: sim_num_fac_II_to_IV}.

\begin{table}[ht!]
\begin{center}
\begin{tabular}{|cccccccc|}
\hline
\multicolumn{8}{|c|}{Setting (IIa)  (True $r_1=2$)} \\ \hline
\multicolumn{1}{|c|}{Missing Pattern} & \multicolumn{1}{c|}{$\wh{r}$} & \multicolumn{1}{c|}{$\wh{r}_{\text{re,0}}$} & \multicolumn{1}{c|}{$\wh{r}_{\text{re,1}}$} & \multicolumn{1}{c|}{$\wh{r}_{\text{iTIP,re,0}}$} & \multicolumn{1}{c|}{$\wh{r}_{\text{iTIP,re,1}}$} & \multicolumn{1}{c|}{\cellcolor[HTML]{9AFF99}$\wh{r}_{\text{full}}$}              & \cellcolor[HTML]{9AFF99}$\wh{r}_{\text{iTIP,full}}$        \\ \hline
\multicolumn{1}{|c|}{}                & \multicolumn{7}{c|}{$\text{Mean}_{(\text{SD})}$}                            \\ \hline
\multicolumn{1}{|c|}{(M-ii)}          & \multicolumn{1}{c|}{$1.98_{(.13)}$}  & \multicolumn{1}{c|}{$1.98_{(.13)}$}     & \multicolumn{1}{c|}{$2.00_{(.06)}$}      & \multicolumn{1}{c|}{$1.97_{(.18)}$}         & \multicolumn{1}{c|}{$1.97_{(.22)}$}          & \multicolumn{1}{c|}{\cellcolor[HTML]{9AFF99}}                     & \cellcolor[HTML]{9AFF99}                     \\ \cline{1-6}
\multicolumn{1}{|c|}{(M-iii)}         & \multicolumn{1}{c|}{$1.92_{(.27)}$}  & \multicolumn{1}{c|}{$1.93_{(.26)}$}     & \multicolumn{1}{c|}{$1.97_{(.20)}$}     & \multicolumn{1}{c|}{$1.90_{(.30)}$}        & \multicolumn{1}{c|}{$1.92_{(.31)}$}         & \multicolumn{1}{c|}{\cellcolor[HTML]{9AFF99}}                     & \cellcolor[HTML]{9AFF99}                     \\ \cline{1-6}
\multicolumn{1}{|c|}{(M-iv)}          & \multicolumn{1}{c|}{$1.98_{(.14)}$} & \multicolumn{1}{c|}{$1.98_{(.14)}$}    & \multicolumn{1}{c|}{$2.01_{(.08)}$}     & \multicolumn{1}{c|}{$1.97_{(.17)}$}        & \multicolumn{1}{c|}{$1.98_{(.24)}$}         & \multicolumn{1}{c|}{\multirow{-3}{*}{\cellcolor[HTML]{9AFF99}$1.99_{(.10)}$}}  & \multirow{-3}{*}{\cellcolor[HTML]{9AFF99}$1.92_{(.28)}$}  \\ \hline
\multicolumn{1}{|c|}{}                & \multicolumn{7}{c|}{Correct Proportion}                 \\ \hline
\multicolumn{1}{|c|}{(M-ii)}          & \multicolumn{1}{c|}{.982} & \multicolumn{1}{c|}{.982}    & \multicolumn{1}{c|}{.996}     & \multicolumn{1}{c|}{.967}        & \multicolumn{1}{c|}{.949}         & \multicolumn{1}{c|}{\cellcolor[HTML]{9AFF99}}                     & \cellcolor[HTML]{9AFF99}                     \\ \cline{1-6}
\multicolumn{1}{|c|}{(M-iii)}         & \multicolumn{1}{c|}{.921} & \multicolumn{1}{c|}{.93}    & \multicolumn{1}{c|}{.96}     & \multicolumn{1}{c|}{.901}        & \multicolumn{1}{c|}{.898}         & \multicolumn{1}{c|}{\cellcolor[HTML]{9AFF99}}                     & \cellcolor[HTML]{9AFF99}                     \\ \cline{1-6}
\multicolumn{1}{|c|}{(M-iv)}          & \multicolumn{1}{c|}{.979} & \multicolumn{1}{c|}{.979}    & \multicolumn{1}{c|}{.993}     & \multicolumn{1}{c|}{.97}        & \multicolumn{1}{c|}{.943}         & \multicolumn{1}{c|}{\multirow{-3}{*}{\cellcolor[HTML]{9AFF99}.99}} & \multirow{-3}{*}{\cellcolor[HTML]{9AFF99}.917} \\ \hline
\end{tabular}
\end{center}
\caption {Results for setting (IIa). Each column reports the mean and SD (subscripted, in bracket) of the estimated number of factors over 1000 replications, followed by the correct proportion  of the estimates. The estimator $\wh{r}$ is our proposed estimator; $\wh{r}_{\text{re,0}}$ and $\wh{r}_{\text{re,1}}$ are similar but used imputed data where the imputation is done using the number of factors as $\wh{r}$ and $\wh{r} +1$, respectively; $\wh{r}_{\text{iTIP,re,0}}$ and $\wh{r}_{\text{iTIP,re,1}}$ are iTIP-ER on imputed data (using $\wh r$ and $\wh r + 1$ respectively); $\wh{r}_{\text{full}}$ and $\wh{r}_{\text{iITP,full}}$ are our estimator and iTIP-ER on fully observed data (in green), respectively.}
\label{tab: sim_num_fac_IIa}
\end{table}

\begin{table}[ht!]
\begin{center}
\begin{tabular}{|cccllccc|}
\hline
\multicolumn{8}{|c|}{Correct Proportion}              \\ \hline
\multicolumn{1}{|c|}{Setting} & \multicolumn{1}{c|}{$\wh{r}$} & \multicolumn{1}{c|}{$\wh{r}_{\text{re,0}}$} & \multicolumn{1}{c|}{$\wh{r}_{\text{re,1}}$} & \multicolumn{1}{c|}{$\wh{r}_{\text{iTIP,re,0}}$} & \multicolumn{1}{c|}{$\wh{r}_{\text{iTIP,re,1}}$} & \multicolumn{1}{c|}{\cellcolor[HTML]{9AFF99}$\wh{r}_{\text{full}}$} & \cellcolor[HTML]{9AFF99}$\wh{r}_{\text{iITP,full}}$ \\ \hline
\multicolumn{8}{|c|}{$K=1$ (True $r_1=2$)}                             \\ \hline
\multicolumn{1}{|c|}{(IIb)}   & \multicolumn{1}{c|}{.556}  & \multicolumn{1}{c|}{.556}      & \multicolumn{1}{c|}{.886}      & \multicolumn{1}{c|}{.526}           & \multicolumn{1}{c|}{.765}             & \multicolumn{1}{c|}{\cellcolor[HTML]{9AFF99}.633}        & \cellcolor[HTML]{9AFF99}.53              \\ \hline
\multicolumn{1}{|c|}{(IIc)}   & \multicolumn{1}{c|}{.626}  & \multicolumn{1}{c|}{.626}      & \multicolumn{1}{c|}{.762}      & \multicolumn{1}{c|}{.594}           & \multicolumn{1}{c|}{.668}             & \multicolumn{1}{c|}{\cellcolor[HTML]{9AFF99}.67}        & \cellcolor[HTML]{9AFF99}.539             \\ \hline
\multicolumn{1}{|c|}{(IId)}   & \multicolumn{1}{c|}{.791} & \multicolumn{1}{c|}{.791}     & \multicolumn{1}{c|}{.817}      & \multicolumn{1}{c|}{.794}           & \multicolumn{1}{c|}{.837}            & \multicolumn{1}{c|}{\cellcolor[HTML]{9AFF99}.812}       & \cellcolor[HTML]{9AFF99}.767             \\ \hline
\multicolumn{8}{|c|}{$K=2$ (True $(r_1,r_2)=(2,3)$)}                             \\ \hline
\multicolumn{1}{|c|}{(IIIa)}  & \multicolumn{1}{c|}{1} & \multicolumn{1}{c|}{1}     & \multicolumn{1}{c|}{1}      & \multicolumn{1}{c|}{.995}           & \multicolumn{1}{c|}{.995}            & \multicolumn{1}{c|}{\cellcolor[HTML]{9AFF99}1}       & \cellcolor[HTML]{9AFF99}.994             \\ \hline
\multicolumn{1}{|c|}{(IIIb)}  & \multicolumn{1}{c|}{.978} & \multicolumn{1}{c|}{.978}     & \multicolumn{1}{c|}{.987}      & \multicolumn{1}{c|}{.985}          & \multicolumn{1}{c|}{.989}            & \multicolumn{1}{c|}{\cellcolor[HTML]{9AFF99}.981}       & \cellcolor[HTML]{9AFF99}.986             \\ \hline
\multicolumn{1}{|c|}{(IIIc)}  & \multicolumn{1}{c|}{.999} & \multicolumn{1}{c|}{.999}     & \multicolumn{1}{c|}{1}     & \multicolumn{1}{c|}{1}          & \multicolumn{1}{c|}{.996}            & \multicolumn{1}{c|}{\cellcolor[HTML]{9AFF99}.999}       & \cellcolor[HTML]{9AFF99}1             \\ \hline
\multicolumn{8}{|c|}{$K=3$ (True $(r_1,r_2,r_3) = (2,3,4)$)}                              \\ \hline
\multicolumn{1}{|c|}{(IVa)}   & \multicolumn{1}{c|}{1} & \multicolumn{1}{c|}{1}     & \multicolumn{1}{c|}{1}     & \multicolumn{1}{c|}{.987}          & \multicolumn{1}{c|}{.987}            & \multicolumn{1}{c|}{\cellcolor[HTML]{9AFF99}1}       & \cellcolor[HTML]{9AFF99}.988             \\ \hline
\multicolumn{1}{|c|}{(IVb)}   & \multicolumn{1}{c|}{.996} & \multicolumn{1}{c|}{.996}     & \multicolumn{1}{c|}{.999}     & \multicolumn{1}{c|}{.991}          & \multicolumn{1}{c|}{.991}            & \multicolumn{1}{c|}{\cellcolor[HTML]{9AFF99}1}       & \cellcolor[HTML]{9AFF99}.991             \\ \hline
\multicolumn{1}{|c|}{(IVc)}   & \multicolumn{1}{c|}{1} & \multicolumn{1}{c|}{1}     & \multicolumn{1}{c|}{1}     & \multicolumn{1}{c|}{.999}          & \multicolumn{1}{c|}{1}            & \multicolumn{1}{c|}{\cellcolor[HTML]{9AFF99}1}       & \cellcolor[HTML]{9AFF99}1             \\ \hline
\end{tabular}
\end{center}
\caption {Results for settings (II), (III), and (IV), excluding (IIa).  Refer to to Table \ref{tab: sim_num_fac_IIa} for the definitions of different estimators. The missing pattern concerned in all settings is (M-ii).}
\label{tab: sim_num_fac_II_to_IV}
\end{table}

From both Table \ref{tab: sim_num_fac_IIa} and \ref{tab: sim_num_fac_II_to_IV}, it is easy to see that our proposed method generally gives more accurate estimates than iTIP-ER, and it is clear that the re-imputation estimate is at least as good as the initial estimate. In fact, $\wh{r}_{\text{re,1}}$ outperforms $\wh{r}_\text{full}$ which is based on full observation.

\subsubsection{Asymptotic normality}\label{subsec:empirical_asymp_normality}
We present the asymptotic normality results for $K=1,2,3$ respectively. When the data is a vector time series ($K=1$), our approach is similar to \cite{Xiong_Pelger}, but their proposed covariance estimator for the asymptotic normality includes information at lag 0 only (i.e., the estimator of the asymptotic variance of the loading estimator), while we use the HAC-type estimator facilitating more serial information. For all $K$ considered, we present the result on $(\wh\Q)_{11}$, with the parameter $\beta$ of our HAC-type estimator set as $\lfloor \frac{1}{5} (Td_1)^{1/4} \rfloor$. We use (M-i) as the missing pattern for all settings.

The data generating process is similar to the ones for assessing the factor loading matrix estimators and imputations, but the parameters are slightly adjusted. All elements in $\cF_t$, $\cF_{e,t}$, and $\bm\epsilon_t$ are now independent standardised AR(1) with AR coefficients 0.05, and we use i.i.d. $\cN(0,1)$ as the innovation process. We stress that we include contemporary and serial dependence among the noise variables through our construction following Assumption (E1) and (E2), while most existing literature demonstrating asymptotic normality display results only for i.i.d. Gaussian noise.

\begin{figure}[ht!]
\begin{center}
\centerline{\includegraphics[width=\columnwidth,scale=1]{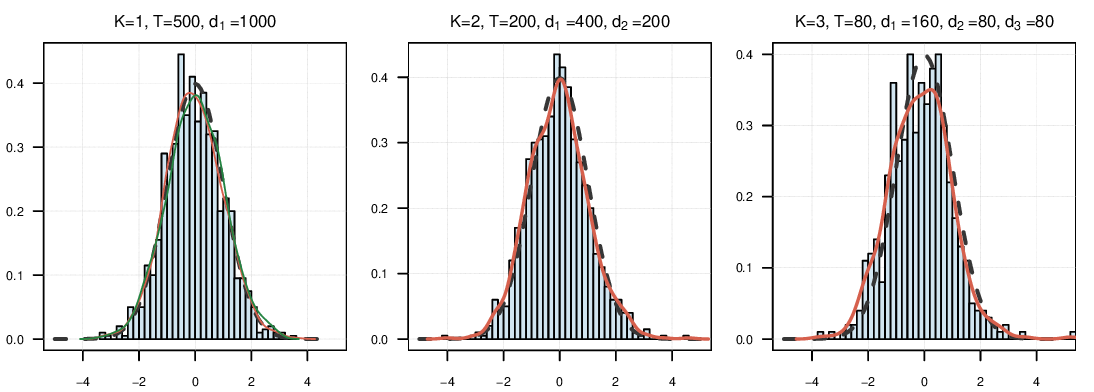}}
\caption{Histograms of the first entry of $(\wh\bSigma_{HAC} + \wh\bSigma_{HAC}^{\Delta})^{-1/2}\wh\D_1(\wh\Q_{1,1\cdot}-\H_1^a\Q_{1,1\cdot})$. In each panel, the curve (in red) is the empirical density, and the other curve (in green) in the left panel depicts the empirical density of the second entry of $(\wh\bSigma_{HAC} + \wh\bSigma_{HAC}^{\Delta})^{-1/2}\wh\D_1(\wh\Q_{1,1\cdot}-\H_1^a\Q_{1,1\cdot})$. The density curve for $\cN(0,1)$ (in black, dotted) is also superimposed on each histogram.}
\label{Fig: sim_asymp_hist}
\end{center}
\end{figure}

\begin{figure}[ht!]
\begin{center}
\centerline{\includegraphics[width=\columnwidth,scale=1]{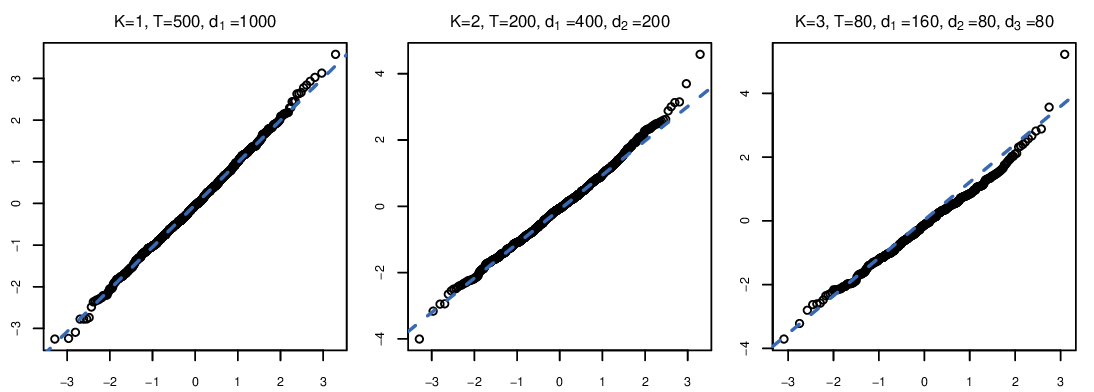}}
\caption{QQ plots of the first entry of $(\wh\bSigma_{HAC} + \wh\bSigma_{HAC}^{\Delta})^{-1/2}\wh\D_1(\wh\Q_{1,1\cdot}-\H_1^a\Q_{1,1\cdot})$. The horizontal and vertical axes are theoretical and empirical quantiles respectively.}
\label{Fig: sim_asymp_QQ}
\end{center}
\end{figure}

We assume all factors are pervasive in this section. For all $K=1,2,3$, given $d_1$, we set $T,d_i=d_1/2$, $i\neq 1$. We generate a two-factor model for $K=1$, and a one-factor model for $K=2,3$. For the settings $(K,d_1)=(1,1000),(2,400), (3,160)$, we consider
$(\wh\bSigma_{HAC} + \wh\bSigma_{HAC}^{\Delta})^{-1/2}\wh\D_1(\wh\Q_{1,1\cdot}-\H_1^a\Q_{1,1\cdot})$. In particular, we plot the histograms of the first and second entry
in Figure \ref{Fig: sim_asymp_hist}, whereas the corresponding QQ plots are presented in Figure \ref{Fig: sim_asymp_QQ}.

The plots in Figure \ref{Fig: sim_asymp_hist} provide empirical support to Theorem \ref{thm:asymp_normality_loadings} and result 3 of Theorem \ref{thm:covariance_estimator}. For $K=3$, there are some heavy-tail issues, as seen in the bump at the right tail in the histogram (confirmed by its corresponding QQ plot). The QQ plot for $K=2$ also hints on this, but the tail is thinned as the dimension increases.
Our simulation is similar to that in \cite{ChenFan2023} for $K=2$, but we allow partial data  unobserved and we generalize to any tensor order $K$. We remark that the convergence rate of the HAC-type estimator is not completely satisfactory, such that relatively large dimension is needed, and it becomes less feasible for some applications. We leave the improvements of the HAC-type estimator to future work.

\begin{figure}[ht!]
\begin{center}
\centerline{\includegraphics[width=0.8\columnwidth]{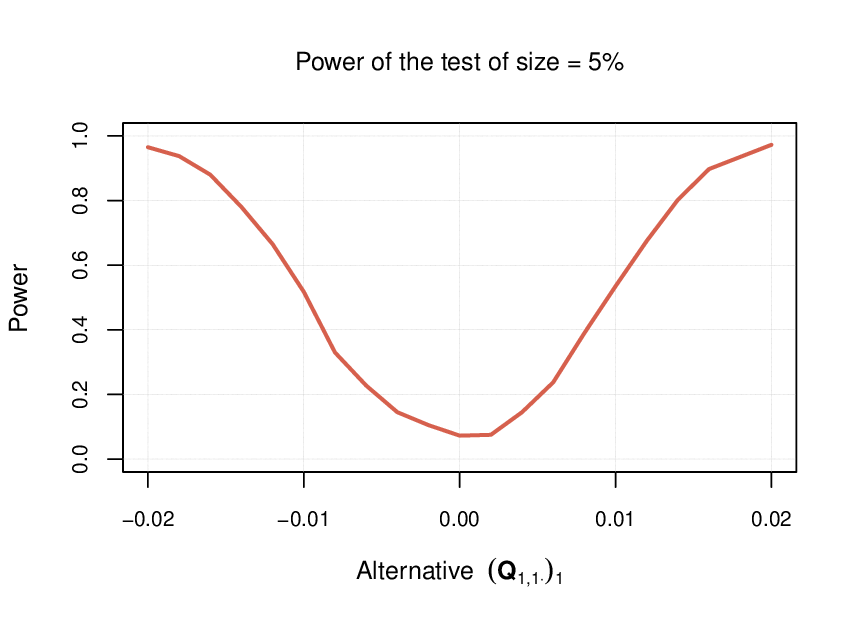}}
\caption{Statistical power of testing the null hypothesis $\c{H}_0:(\Q_{1,1\cdot})_1 = \Q_{1,11} =0$ against the general alternative. The null is rejected when $|(\wh\bSigma_{HAC} + \wh\bSigma_{HAC}^{\Delta})^{-1/2}\wh\D_1\wh\Q_{1,11}|> 1.96$.}
\label{Fig: sim_asymp_test_power}
\end{center}
\end{figure}

Lastly, we demonstrate an example of statistical testing for the above one-factor model for $K=2$. More precisely, we want to test the null hypothesis $\c{H}_0:\Q_{1,11}=0$ with a two-sided test. A $5\%$ significance level is used so that we reject the null if $(\wh\bSigma_{HAC} + \wh\bSigma_{HAC}^{\Delta})^{-1/2}\wh\D_1\wh\Q_{1,11}$ is not in $[-1.96,1.96]$. Each experiment is repeated 400 times and the power function for $\Q_{1,11}$ ranging from $-0.02$ to $0.02$ is presented in Figure \ref{Fig: sim_asymp_test_power}.
The power function is approximately symmetric, and suggests that our test can successfully reject the null if the true value for $\Q_{1,11}$ is away from 0. When $\Q_{1,11}=0$, the false positive probability is $7.25\%$ which is slightly higher than the designated size of test. This is due to the slow convergence of the HAC estimators, and an increase in dimensions would improve this.

{
\subsubsection{Comparison with an iterative vectorization-based approach}\label{subsec:compare_iter_vec}
}
We compare our proposed tensor factor-based imputation method with the following procedure.

\noindent \underline{\textit{Iterative vectorisation-based imputation}}
\begin{enumerate}
    \item Given an order-$K$ tensor with missing entries, $\c{Y}_t \in \b{R}^{d_1\times \dots\times d_K}$ for $t\in[T]$, obtain $\bf{y}_t = \vec{\c{Y}_t} \in\b{R}^d$ for all time stamps. Impute the vector time series $\{\bf{y}_t \}_{t\in[T]}$ by \cite{Xiong_Pelger} and denote by the tensorized imputation data $\{ \wh{ \c{Y}}_{\text{vec}, t}\}_{t \in[T]}$.
    \item Replace missing entries in $\c{Y}_t$ by the corresponding entries in $\wh{ \c{Y}}_{\text{vec}, t}$. For the resulting time series, estimate the loading matrices, core factors and hence the common components by \cite{Chen_Lam}. Denote the series of estimated common components by $\{ \wh{ \c{Y}}_{\text{preavg}, t}\}_{t \in[T]}$.
    \item Iterate from step 2, except that we replace missingness of $\c{Y}_t$ by $\wh{ \c{Y}}_{\text{preavg}, t}$ from the previous iteration.
\end{enumerate}

The above algorithm is a natural way of leveraging the vector imputation of \cite{Xiong_Pelger} to tensor time series, and the iteration step is akin to Appendix A of \cite{StockWatson2002b}. For demonstration, all innovation processes in constructing $\cF_t, \cF_{e,t}$ and $\bepsilon_t$ are i.i.d. standard normal, and all factors are pervasive. In particular, the following settings are considered:
\begin{itemize}
\item[(Va)] $K=2, T=20, d_1=d_2=40, r_1= r_2=2$, and missing pattern is (M-ii).
\item[(Vb)] Same as (Va), except that the missing pattern is (M-iii).
\item[(Vc)] $K=3, T=10, d_1=d_2=d_3=10, r_1= r_2= r3 =2$, and missing pattern is (M-ii).
\item[(Vd)] Same as (Vc), except that the missing pattern is (M-iii).
\end{itemize}

\begin{figure}[ht!]
\begin{center}
\centerline{\includegraphics[width=\columnwidth,scale=1]{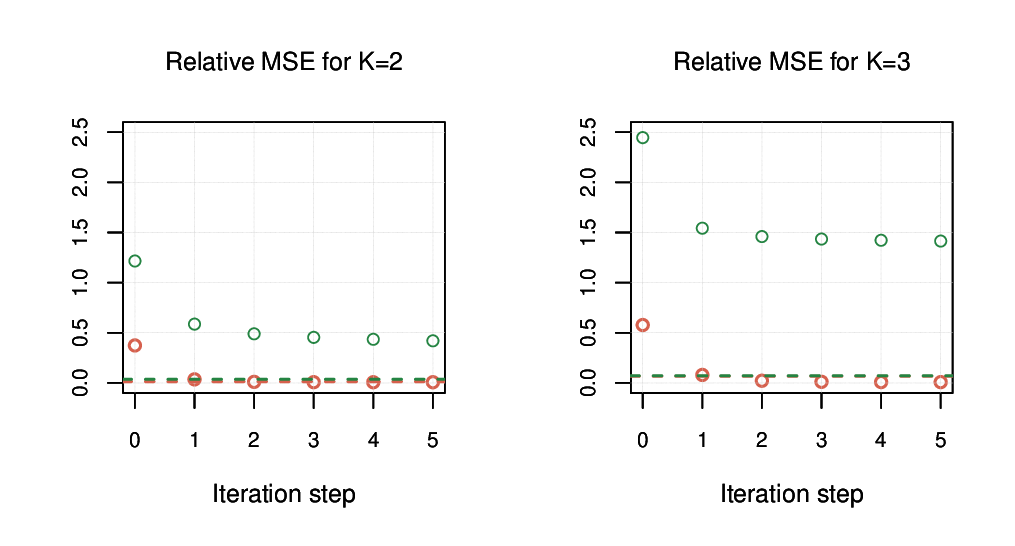}}
\caption{Plot of the relative MSE for settings (Va) to (Vd), averaged over 1000 replications. Dashed lines are the errors of our tensor factor-based approach (without iteration), points are errors of the iterative vectorization-based method at each step. Step 0 corresponds to the initial vectorized imputation. Results under missing patterns (M-ii) and (M-iii) are presented in red and green, respectively.}
\label{Fig: compare_iter_vec}
\end{center}
\end{figure}

The results for settings (Va) to (Vd) are shown in Figure \ref{Fig: compare_iter_vec}. From both panels, our proposed method (in dashed lines) performs better than the direct vectorized imputation. One intuition can be the following. Suppose we have a matrix-valued time series $\Y_t\in \b{R}^{d_1\times d_2}$ for $t\in[T]$, and assume $d_1 \asymp d_2$ and the data is asymptotically observed with the rate $\eta \asymp 1/\sqrt{T d_1}$. According to Corollary \ref{corollary:simplified_consistency}, the squared imputation error has rate $1/(Td_1) + 1/d_1^2$. In comparison, if we choose to vectorize the data and impute, the squared error rate is $1/T + 1/d_1^2$ which is inflated.

The performance of the vectorization-based imputation can be further improved by iterative imputation in the context of tensor data. However, Figure \ref{Fig: compare_iter_vec} demonstrates the low efficiency of such iterative method if the missing pattern is unbalanced to a certain extent. We also point out that the computation time of the initial vectorized imputations can be significantly larger than the our proposed method if the order of the data is large. In fact, the computational complexity (given the number of factors) of direct vectorized imputation is (ignoring the cost of vectorization and unfolding) $O(Td^2 + d^3)$, while our proposed method is $O(K\max_{k\in[K]}\{Tdd_k + d_k^3\})$, which can be of significantly smaller order than $d^3$.

{
\subsection{Real data analysis: Fama-French portfolio returns}\label{subsec:famafrench}
}
We analyze the set of Fama-French portfolio returns data described in Section \ref{subsec:fama-french_data}.
With sufficient observed samples of each category along its time series, Assumption (O1) in Section \ref{sec:Assumption_theories} can be satisfied and our imputation approach is applicable under such missing pattern.
Since the market factor is pervasive in financial returns, we remove the market effect by modelling the data with CAPM as
\[
\vec{\c{X}_t} = \vec{ \bar{\c{X}}} + \bm{\beta} (r_t -\bar{r}) + \vec{\c{Y}_t},
\]
where $\vec{\c{X}_t} \in\b{R}^{100}$ is the vectorized returns at time $t$, $\vec{ \bar{\c{X}}}$ is the sample mean of $\vec{\c{X}_t}$, $\bm{\beta}$ is the coefficient vector to be estimated, $r_t$ is the return of the NYSE composite index at time $t$, $\bar{r}$ is the sample mean of $r_t$, and $\vec{\c{Y}_t}$ is the CAPM residual. We compute the sample mean using only the observed data, and more sophisticated methods could be studied in the future. The least squares solution is
\[
\wh{\bm{\beta}} = \frac{\sum_{t=1}^T (r_t -\bar{r}) \{\vec{\c{X}_t} -\vec{ \bar{\c{X}}} \}}{\sum_{t=1}^T (r_t -\bar{r})^2}.
\]
Hence for the rest of this section, we focus on the matrix series $\{ \wh{\c{Y}}_t\}_{t \in[570]}$ with $\wh{\c{Y}}_t \in\b{R}^{10\times 10}$, constructed from the estimated CAPM residual $\{\vec{\c{X}_t} -\vec{ \bar{\c{X}}} - \wh{\bm{\beta}} (r_t -\bar{r}) \}_{t \in[570]}$.

To estimate the rank of the core factors, we first use our proposed rank estimator to obtain initial estimates $(\wh{r}_1, \wh{r}_2) =(1,1)$ for both series, followed by re-estimating the rank based on the imputed series using $(\wh{r}_1+ r_\ast, \wh{r}_2+ r_\ast)$ with some pre-specified integer $r_\ast$ to capture any omitted weak factors. We have seen in Table \ref{tab: sim_num_fac_IIa} and Table \ref{tab: sim_num_fac_II_to_IV} where such rank re-estimation with $r_\ast= 1$ is stable and accurate. However, factors can be empirically too weak to detect in the initial estimation under various missing patterns, see the NYC taxi tensor data analysis in the supplementary materials as an extra example. According to previous studies by e.g. \cite{Wangetal2019}, we choose $r_\ast= 3$ here to ensure sufficient information of factors is carried in the imputation, at the cost of including more noise. For re-estimation, in addition to our eigenvalue-ratio estimator, we also experiment BCorTh by \cite{Chen_Lam}, iTIP-ER by \cite{Hanetal2022} and RTFA-ER by \cite{Heetal2022b}. The results are presented in Table \ref{tab: famafrench_num_fac}. To ease demonstration, we use $(2,2)$ as the core factor rank for both series hereafter.

\begin{table}[t!]\centering
\small
\ra{1.3}
\begin{tabular}{@{}rrcrrcrrcrrcrrc@{}}\toprule
& \multicolumn{2}{c}{initial} & \phantom{ab}
& \multicolumn{2}{c}{ Miss-ER } & \phantom{ab} & \multicolumn{2}{c}{ BCorTh } & \phantom{ab} & \multicolumn{2}{c}{ iTIP-ER } & \phantom{ab} & \multicolumn{2}{c}{ RTFA-ER } \\
\cmidrule{2-3} \cmidrule{5-6} \cmidrule{8-9} \cmidrule{11-12} \cmidrule{14-15}
& $\wh{r}_1$ & $\wh{r}_2$  &&  $\wh{r}_1$ & $\wh{r}_2$ && $\wh{r}_1$ & $\wh{r}_2$ && $\wh{r}_1$ & $\wh{r}_2$ && $\wh{r}_1$ & $\wh{r}_2$ \\ \midrule
Value Weighted & $1$ & $1$  &&  $1$ & $1$ && $2$ & $1$ && $1$ & $1$ && $1$ & $2$ \\
Equal Weighted & $1$ & $1$  &&  $1$ & $1$ && $2$ & $1$ && $1$ & $1$ && $1$ & $2$ \\
\bottomrule
\end{tabular}
\caption {Rank estimators for Fama-French portfolios. Miss-ER represents the rank re-estimated by our proposed eigenvalue-ratio estimator for missing data.}
\label{tab: famafrench_num_fac}
\end{table}

\begin{table}[ht!]
\setlength{\tabcolsep}{4pt}
\begin{center}
\begin{tabular}{l||ccccccccccccccccccccccccccc}
\hline & ME1 & ME2 & ME3 & ME4 & ME5 & ME6 & ME7 & ME8 & ME9 & ME10 \\
\hline Factor 1 & \red{-15} & \red{-14} & -9 & -7 & -6 & -3 & -1 & 0 & 2 & {3} \\
Factor 2 & 5 & 3 & -3 & -6 & -7 & -9 & \red{-10} & \red{-11} & \red{-10} & \red{-10} \\
\hline
\hline  & BE1 & BE2 & BE3 & BE4 & BE5 & BE6 & BE7 & BE8 & BE9 & BE10 \\
\hline Factor 1 & 2 & -1 & -2 & -3 & -5 & -6 & -7 & -8 & \red{-10} & \red{-18} \\
Factor 2 & \red{16} & \red{12} & 9 & {7} & 5 & 3 & 3 & 2 & 1 & -7 \\
\hline
\end{tabular}
\end{center}
\caption{Estimated loading matrices $\wh\Q_1$ and $\wh\Q_2$ for the value weighted portfolio series, after varimax rotation and scaling (entries rounded to the nearest integer). Magnitudes larger than 9 are highlighted in red. All null hypotheses of a row of $\Q_1$ or $\Q_2$ being zero (see (\ref{eqn: famafrench_inference})) are rejected at $5\%$ significance level.}
\label{tab: famafrench_loading}
\end{table}

With the chosen rank, we perform imputation which is further refined by re-imputation. The results are similar on the two portfolio series, so we only present the one for the value weighted series. The estimated loading matrices are presented in Table \ref{tab: famafrench_loading}, after a varimax rotation and scaling. We can see from the table that on the size factor (i.e., ME loading), ME1 and ME2 form one group (``small size'') and ME7 to ME10 form another group (``large size''). On the book-to-equity factor (i.e., BE loading), BE1 and BE2 form a group and BE9 and BE10 form another, which can be interpreted as ``undervalued'' and ``overvalued'' respectively. This grouping effect turns out to be aligning with \cite{Wangetal2019}, but we have a clearer grouping effect seen.

Moreover, we apply our Theorem \ref{thm:asymp_normality_loadings} and Theorem \ref{thm:covariance_estimator} to test if any rows of the loading matrices are zero. For each $k\in[2], i\in[10]$, we test
\begin{equation}
\label{eqn: famafrench_inference}
\c{H}_0: \Q_{k,i\cdot} = \0, \,\,\,\,\,
\c{H}_1: \Q_{k,i\cdot} \neq \0.
\end{equation}
The above can be tested since $\H_k^a\Q_{k,i\cdot} = \0$ under the null, and no matter what varimax rotations we use, it retains its meaning. For instance, if $\Q_{1,i\cdot} = \0$, then it means that the $i$-th category of the row factor (here, the $i$-th Market Equity category) is useless in explaining any data variability.

It turns out that at 5\% significance level, we cannot reject any null hypotheses for $\Q_{1,i\cdot} = \0$ or $\Q_{2,i\cdot}=\0$, meaning that individual market equity and book-to-equity ratio categories are tested to be meaningful in explaining some variations of the data. See the NYC Taxi traffic data analysis in the supplementary materials for some similar null hypotheses not rejected. 
We remark that, since the dimensions of our data are not very large, the accuracy of the asymptotic normality and the HAC estimators are weakened, and there can be false positives as a result.

\begin{figure}[ht!]
\begin{center}
\centerline{\includegraphics[width=\columnwidth,scale=1]{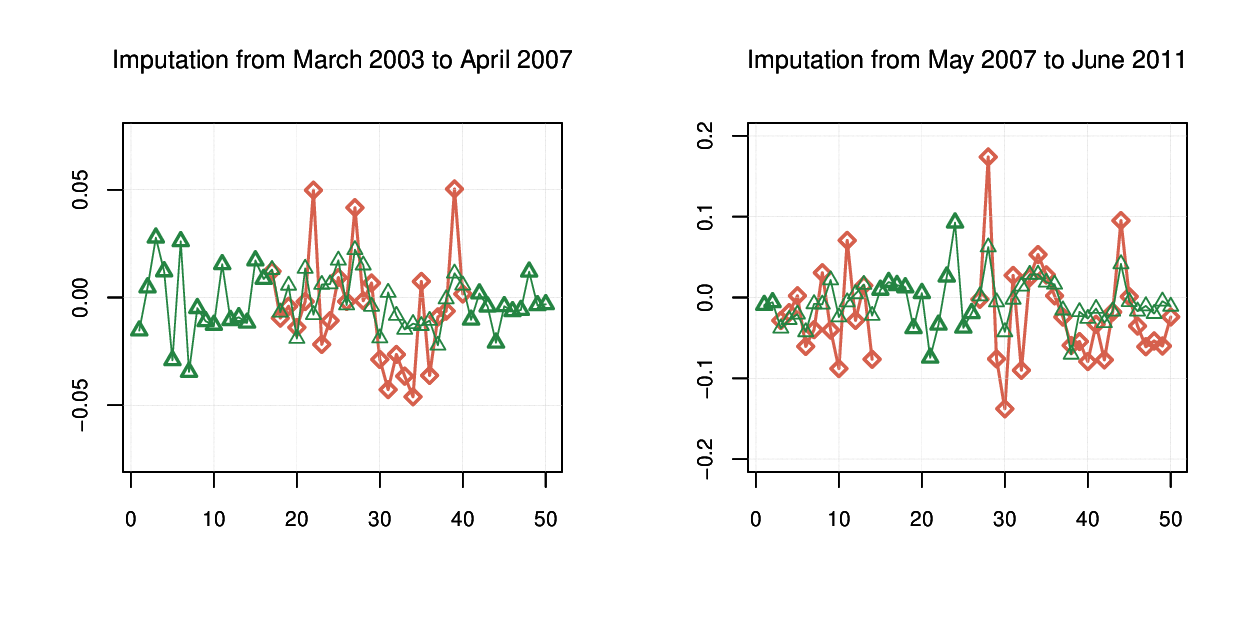}}
\caption{Imputations for the value weighted series in the category (ME10, BE10). The observed and estimated series are represented by red squares and green triangles, respectively. The imputed series (bold) is formed by replacing missing entries with the estimated data using the CAPM described in the section.}
\label{Fig: famafrench_imputation}
\end{center}
\end{figure}

Lastly, two imputation examples are displayed in Figure \ref{Fig: famafrench_imputation}, where we also show the estimated series on timestamps for which the portfolio series is observed. The estimated series does capture some of the pattern of fluctuations, and can be a good reference for the CAPM residual of portfolios consisted of large size, overvalued stocks. This is certainly more revealing than a naive imputation using zeros or local means. From the above discussions, the estimated factors can be potentially used to replace the Fama–French size factor (SMB) and book-to-equity factor (HML) in a Fama–French factor model for asset pricing, factor trading etc, with a more sophisticated further analysis of the data.

\vspace{-12pt}

{
\subsection{Real data analysis: OECD Economic indicators for countries}\label{subsec:OECD}
}
We analyze the OECD Economic data described in Section \ref{subsec:OECD_data}. After investigating the estimated number of factors (Table \ref{tab: OECD_num_fac}) in a similar re-imputation approach as in Section \ref{subsec:famafrench}, we decide to use $(\wh{r}_1, \wh{r}_2) =(3,3)$ for the rest of this section due to the potentially weak factors suggested by iTIP-ER and RTFA-ER. The estimated loading matrices for countries are presented in Table \ref{tab: OECD_loading_country} after a varimax rotation and scaling. The first factor is mainly formed by European countries except the Northern European ones which, together with Canada, form the third factor. Such regional grouping effects are also confirmed in the second factor which mainly consists of the United States, and the fact that Germany loads also heavily in this factor suggests their similar economic patterns as large economic entities. For the estimated loading for indicators reported in Table \ref{tab: OECD_loading_indicator}, CP, PRVM and TOVM form the first factor (``consumption factor''), PP and ULC form the second (``production factor''), and EX and IM form the third (``international trade factor'').

Moreover, we apply Theorem \ref{thm:asymp_normality_loadings} and Theorem \ref{thm:covariance_estimator} to test if a particular row in the two factor loading matrices is zero, meaning that if a country (if a row in $\Q_1$ is $\0$) or an economic indicator (if a row in $\Q_2$ is $\0$) cannot explain any variations in the data. The meaning here is independent of the varimax rotation performed. For each $k\in[2], i\in[d_k], j\in[3]$ with $(d_1,d_2) =(17,11)$, we form the hypothesis
\begin{equation}
\label{eqn: OECD_inference}
\c{H}_0: \Q_{k,i\cdot} = \0, \,\,\,\,\,
\c{H}_1: \Q_{k,i\cdot} \neq \0.
\end{equation}
Similar to the Fama-French data analysis, all null hypotheses of a row of $\Q_1$ or $\Q_2$ being zero are rejected at 5\% significance level. It means that all individual country and economic indicator are tested to be meaningful categories in explaining some variations of the data.
Similar to a reminder in Section \ref{subsec:famafrench}, there could be false positives due to the fact that the dimension of the data is not very large.

\begin{table}[t!]\centering
\small
\ra{1.3}
\begin{tabular}{@{}rrcrrcrrcrrcrrc@{}}\toprule
& \multicolumn{2}{c}{initial} & \phantom{ab}
& \multicolumn{2}{c}{ Miss-ER } & \phantom{ab} & \multicolumn{2}{c}{ BCorTh } & \phantom{ab} & \multicolumn{2}{c}{ iTIP-ER } & \phantom{ab} & \multicolumn{2}{c}{ RTFA-ER } \\
\cmidrule{2-3} \cmidrule{5-6} \cmidrule{8-9} \cmidrule{11-12} \cmidrule{14-15}
& $\wh{r}_1$ & $\wh{r}_2$  &&  $\wh{r}_1$ & $\wh{r}_2$ && $\wh{r}_1$ & $\wh{r}_2$ && $\wh{r}_1$ & $\wh{r}_2$ && $\wh{r}_1$ & $\wh{r}_2$ \\ \midrule
OECD & $1$ & $1$  &&  $1$ & $1$ && $1$ & $2$ && $4$ & $5$ && $3$ & $3$ \\
\bottomrule
\end{tabular}
\caption {Rank estimators for economic indicators. Refer to Table \ref{tab: famafrench_num_fac} for the definitions of different estimators.}
\label{tab: OECD_num_fac}
\end{table}

\begin{table}[ht!]
\footnotesize
\setlength{\tabcolsep}{4pt}
\begin{center}
\begin{tabular}{l||cccccccccccccccccccccccccccccccccc}
\hline & BEL & CAN & DNK & FIN & FRA & DEU & GRC & ITA & LUX & NLD & NOR & PRT & ESP & SWE & CHE & GBR & USA \\
\hline 1 & \red{-10} & 4 & -3 & -7 & -8 & -8 & -9 & -7 & -1 & -2 & 1 & -1 & \red{-10} & 0 & \red{-15} & \red{-13} & 1 \\
2 & -1 & -6 & 2 & 5 & -2 & \red{-12} & 7 & -1 & 2 & -7 & 0 & 2 & 1 & -2 & 0 & -1 & \red{-24} \\
3 & 1 & \red{-12} & -8 & -5 & {-2} & 3 & -6 & -4 & {\red{-11}} & {-6} & {\red{-12}} & \red{-11} & -2 & \red{-10} & {5} & 4 & -1 \\
\hline
\end{tabular}
\end{center}
\caption{Estimated loading matrix $\wh\Q_1$ on three country factors for the OECD data, after varimax rotation and scaling (entries rounded to the nearest integer). Magnitudes larger than 9 are highlighted in red. All null hypotheses of a row of $\Q_1$ being zero (see (\ref{eqn: OECD_inference})) are rejected at $5\%$ significance level.}
\label{tab: OECD_loading_country}
\end{table}

\begin{table}[ht!]
\footnotesize
\setlength{\tabcolsep}{4pt}
\begin{center}
\begin{tabular}{l||cccccccccccccccccccccccccccc}
\hline & CA-GDP & CP & EX & IM & IR3TIB & IRLT & IRSTCI & PP & PRVM & TOVM & ULC \\
\hline 1 & 0 & \red{-20} & 1 & 3 & 0 & 0 & 0 & 0 & \red{-20} & \red{-11} & -1 \\
2 & {0} & -6 & 2 & 3 & 1 & 1 & 1 & \red{20} & 1 & 9 & \red{18} \\
3 & {0} & 9 & \red{18} & \red{22} & -2 & -2 & -2 & {1} & -4 & -2 & {-1} \\
\hline
\end{tabular}
\end{center}
\caption{Estimated loading matrix $\wh\Q_2$ on three indicator factors for OECD data, after varimax rotation and scaling (entries rounded to the nearest integer). Magnitudes larger than 9 are highlighted in red. All null hypotheses of a row of $\Q_2$ being zero (see (\ref{eqn: OECD_inference})) are rejected at $5\%$ significance level.}
\label{tab: OECD_loading_indicator}
\end{table}

\begin{figure}[t!]
\begin{center}
\centerline{\includegraphics[width=\columnwidth,scale=1]{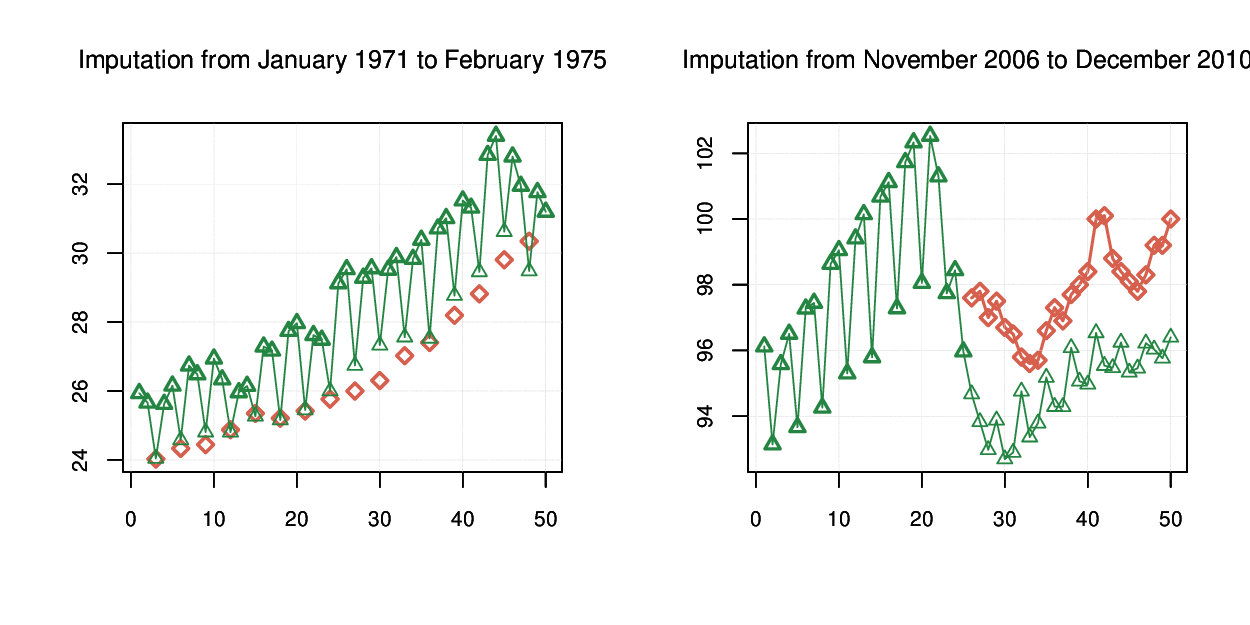}}
\caption{Demonstration of imputations for unit labour cost of the United States (left panel) and production price index of the United Kingdom (right panel). The observed and estimated series are represented by red squares and green triangles respectively. The imputed series (bold) is formed by replacing missing entries with the estimated data.}
\label{Fig: OECD_imputation}
\end{center}
\end{figure}

In Figure \ref{Fig: OECD_imputation}, we present two examples of the imputed series overlaid on the observed series. One panel plots ULC of the United States and the other plots PP of the United Kingdom. ULC is a quarterly observed index and the peak pattern in-between each reported timestamp suggests potentially high labour cost in the United States from 1971 to 1975. The PP data in our OECD data is unavailable for the United Kingdom until December 2008. Our imputation implies a gradual increase of the PP before the data is reported, which is reasonable by the impact of the financial crisis. Lastly, we compare between our tensor imputation (matrix imputation for this example) and the vectorized imputation using \cite{Xiong_Pelger}. We use different models to perform imputations whose results are summarized in Table \ref{tab: OECD_compare_vectorize} similar to \cite{Wangetal2019}, except that the reported residual sum of squares are computed on the observed entries. Although we require a larger number of factors in general for matrix models, the imputation by matrix models with less parameters can perform better than those by vector models with a much larger number of parameters. This is consistent with the conclusion of Table 11 in \cite{Wangetal2019}.

\begin{table}[ht!]
\setlength{\tabcolsep}{4pt}
\begin{center}
\begin{tabular}{llrrr}
\hline & Factor & RSS & $\#$ factors & $\#$ parameters  \\
\hline
Matrix model & (3,3) & 7,087,373 & 9 & 84  \\
Matrix model & (4,4) & 4,542,956 & 16 & 112  \\
Matrix model & (5,5) & 3,066,851 & 25 & 140  \\
Matrix model & (6,6) & 1,973,321 & 36 & 168  \\
\hline
Vector model & 2 & 8,240,976 & 2 & 374 \\
Vector model & 3 & 3,954,554 & 3 & 561 \\
Vector model & 4 & 2,093,001 & 4 & 748 \\
\hline
\end{tabular}
\end{center}
\caption{Comparison of different models for the OECD data. The total sum of squares of the observation is 324,402,709. }
\label{tab: OECD_compare_vectorize}
\end{table}

\section{Appendix}\label{sec:Appendix}
Proofs of all the theorems in this paper can be found in
the supplement of this paper at \linebreak \href{http://stats.lse.ac.uk/lam/Supp-tensorMiss.pdf}{http://stats.lse.ac.uk/lam/Supp-tensorMiss.pdf}.
Instruction in using our R package  \texttt{tensorMiss} can be found \href{http://stats.lse.ac.uk/lam/intro-to-tensorMiss.html}{here}.

\newpage
\appendix

\begin{center}
{\textbf{Supplementary materials for the paper ``Tensor Time Series Imputation through Tensor Factor Modelling'' - Extra simulations, NYC taxi data analysis, and proof of all theorems and propositions}}
\end{center}

\section{Appendix: Extra simulation results and data analysis}\label{sec: Appendix_sim_data}

\subsection{Sensitivity for the tuning parameter in Theorem \ref{thm:number_of_factors}}
\label{subsec:sensitivity_xi}

We provide some simulation results on the performance of our number of factor estimator $\wh{r}_k$ relative to the choice of $\xi$ in Theorem \ref{thm:number_of_factors}. For demonstration purpose, we adapt the general setup depicted in Section \ref{subsec:simulation} to generate the loading matrices, factor and noise series, except that only $\c{N}(0,1)$ is used to generate the innovation process. On the dimension of data, we consider order $K=1,2,3$ such that
\begin{itemize}
    \item $K=1$: $T=160$, $d_1 =80$;
    \item $K=2$: $T= d_1 =d_2 =40$;
    \item $K=3$: $T =d_1 =d_2 =d_3 =20$;
\end{itemize}
where we assume two pervasive factors on each mode, i.e. true $r_k =2$ and $\zeta_{k,j} =0$ for all $k,j$.

To show the robustness of our choice of $\xi$, each setting is repeated 1000 times under four missing patterns: fully observed, (M-i), (M-ii) and (M-iii). See Section \ref{subsec:simulation} for the description of these missing patterns.

\begin{table}[ht!]
\begin{center}
\begin{tabular}{|cccllc|}
\hline
\multicolumn{6}{|c|}{Correct Proportion of $\wh{r}_1 = r_1 =2$}              \\ \hline
\multicolumn{1}{|c|}{Missing Pattern} & \multicolumn{1}{c|}{$\xi=0.002$} & \multicolumn{1}{c|}{$\xi=0.02$} & \multicolumn{1}{c|}{$\xi=0.2$} & \multicolumn{1}{c|}{$\xi=2$} & \multicolumn{1}{c|}{$\xi=20$} \\ \hline
\multicolumn{6}{|c|}{$K=1$ ($T=160, \, d_1 =80$)}                             \\ \hline
\multicolumn{1}{|c|}{Fully observed}   & \multicolumn{1}{c|}{.999}  & \multicolumn{1}{c|}{999}      & \multicolumn{1}{c|}{1}      & \multicolumn{1}{c|}{.986}           & \multicolumn{1}{c|}{.909}       \\ \hline
\multicolumn{1}{|c|}{(M-i)}   & \multicolumn{1}{c|}{.999}  & \multicolumn{1}{c|}{.999}      & \multicolumn{1}{c|}{1}      & \multicolumn{1}{c|}{.985}           & \multicolumn{1}{c|}{.906}      \\ \hline
\multicolumn{1}{|c|}{(M-ii)}   & \multicolumn{1}{c|}{.997} & \multicolumn{1}{c|}{.997}     & \multicolumn{1}{c|}{.995}      & \multicolumn{1}{c|}{.985}           & \multicolumn{1}{c|}{.903}     \\ \hline
\multicolumn{1}{|c|}{(M-iii)}   & \multicolumn{1}{c|}{.99} & \multicolumn{1}{c|}{.989}     & \multicolumn{1}{c|}{.987}      & \multicolumn{1}{c|}{.934}           & \multicolumn{1}{c|}{.791}     \\ \hline
\multicolumn{6}{|c|}{$K=2$ ($T= d_1 =d_2 =40$)}                             \\ \hline
\multicolumn{1}{|c|}{Fully observed}  & \multicolumn{1}{c|}{1} & \multicolumn{1}{c|}{1}     & \multicolumn{1}{c|}{1}      & \multicolumn{1}{c|}{.992}           & \multicolumn{1}{c|}{.843}         \\ \hline
\multicolumn{1}{|c|}{(M-i)}  & \multicolumn{1}{c|}{1} & \multicolumn{1}{c|}{1}     & \multicolumn{1}{c|}{1}      & \multicolumn{1}{c|}{.993}          & \multicolumn{1}{c|}{.842}       \\ \hline
\multicolumn{1}{|c|}{(M-ii)}  & \multicolumn{1}{c|}{1} & \multicolumn{1}{c|}{1}     & \multicolumn{1}{c|}{1}     & \multicolumn{1}{c|}{.986}          & \multicolumn{1}{c|}{.842}         \\ \hline
\multicolumn{1}{|c|}{(M-iii)}   & \multicolumn{1}{c|}{.994} & \multicolumn{1}{c|}{.995}     & \multicolumn{1}{c|}{.996}      & \multicolumn{1}{c|}{.972}           & \multicolumn{1}{c|}{.815}     \\ \hline
\multicolumn{6}{|c|}{$K=3$ ($T= d_1 =d_2 =d_3 =20$)}                              \\ \hline
\multicolumn{1}{|c|}{Fully observed}   & \multicolumn{1}{c|}{.999} & \multicolumn{1}{c|}{.999}     & \multicolumn{1}{c|}{.997}     & \multicolumn{1}{c|}{.967}          & \multicolumn{1}{c|}{.732}           \\ \hline
\multicolumn{1}{|c|}{(M-i)}   & \multicolumn{1}{c|}{.999} & \multicolumn{1}{c|}{.999}     & \multicolumn{1}{c|}{.997}     & \multicolumn{1}{c|}{.967}          & \multicolumn{1}{c|}{.73}     \\ \hline
\multicolumn{1}{|c|}{(M-ii)}   & \multicolumn{1}{c|}{.998} & \multicolumn{1}{c|}{.998}     & \multicolumn{1}{c|}{.996}     & \multicolumn{1}{c|}{.959}          & \multicolumn{1}{c|}{.718}         \\ \hline
\multicolumn{1}{|c|}{(M-iii)}   & \multicolumn{1}{c|}{.995} & \multicolumn{1}{c|}{.997}     & \multicolumn{1}{c|}{.995}      & \multicolumn{1}{c|}{.95}           & \multicolumn{1}{c|}{.695}     \\ \hline
\end{tabular}
\end{center}
\caption {Results of correct proportion for the number of factor estimator on mode-1 in 1000 replications.}
\label{tab: sensitivity_xi}
\end{table}

As the dimension and the number of factors along each tensor mode is the same (within any setting), it suffices to study the correct proportion of $\wh{r}_1 = r_1 = 2$. The result for different values of $\xi \in \{0.002, 0.02, 0.2, 2, 20\}$ is shown in Table \ref{tab: sensitivity_xi}. It is clear from the results that relatively small values of $\xi$ should help to estimate the number of factors consistently. In particular, $\xi$ ranges from $0.02$ to $0.2$ should work sufficiently well.

\subsection{Introduction of new measure}
\label{subsec:q-RSE}
We propose a new measure to gauge the imputation performance in this section. Let the observed data be $\bf{y} = (y_1, y_2, \dots, y_N)'$ and a vector of imputed data be $\wh{\bf{y}} = (\wh{y}_1, \wh{y}_2, \dots, \wh{y}_N)'$.
WLOG, we let the entries in $\bf{y}$ be ordered such that $y_1\leq y_2\leq \dots\leq y_N$, and each $\wh{y}_i$ is the corresponding imputed value for $y_i$. Given $q\in[N]$ and any integer $j=0,1,\ldots, q$, denote by $q_j\in[N]$ the index of entry with value closest to the $(j/q)$-th quantile of $\bf{y}$. Hereafter, quantiles are always with respect to $\bf{y}$ when referred to. For $n\in[q]$, define
\begin{align*}
    \mu_n := \frac{1}{q}\sum_{i=q_{n-1}}^{q_{n}} y_i , \;\;\; \wh\mu_n &:= \frac{1}{q}\sum_{i=q_{n-1}}^{q_{n}} \wh{y}_i .
\end{align*}
Hence, $\mu_n$ is the average of the entries in $\bf{y}$ between the $((n-1)/q)$-th quantile and the $(n/q)$-th quantile.

With these, we define the \textit{$q$-quantile relative squared error (q-{\em RSE})} of $\wh{\bf{y}}$ on $\bf{y}$ as
\begin{equation*}
    q\text{-RSE} := \sum_{n=1}^q (\mu_n - \wh\mu_n)^2 \bigg/ \sum_{n=1}^q \mu_n^2 = \sum_{n=1}^q \bigg\{\sum_{i=q_n}^{q_{n+1}} (y_i -\wh{y}_i) \bigg\}^2 \bigg/ \sum_{n=1}^q \bigg\{\sum_{i=q_n}^{q_{n+1}} y_i\bigg\}^2.
\end{equation*}
It is obvious that $q$-RSE is always non-negative, and $q\text{-RSE}=0$ if $\wh{y}_i=y_i$ for all $i\in[N]$. Note that the relative MSE defined in (\ref{eqn:relative_MSE}) is equivalent to $q$-RSE with $q=N$, and essentially the average within each quantile interval smooths out the idiosyncratic noise in each data point. Intuitively, $q$-RSE can be considered as a measure of how good the imputation of $y_i$ is, compared with the overall level of $y_i$ in each quantile interval.

As an example to demonstrate the pitfall of relative MSE, consider a one-factor model for vector time series and $y_i$ represents some data entry for $i\in[N]$. Suppose $y_i=c_i +\epsilon_i$ and $\wh{y}_i=c_i$, where $c_i$ is the common component to be imputed and $\epsilon_i$ is the noise. Hence we have perfect imputation here, and a good measure should have value close to 0. However,
\[
\text{relative MSE} = \frac{\sum_{i=1}^N \epsilon_i^2}{\sum_{i=1}^N (c_i + \epsilon_i)^2},
\]
which might even be larger than 1 if the noise term is comparable to the common component. For simplicity if we set $q=1$, then
\[
q\text{-RSE} = \frac{(\sum_{i=1}^N \epsilon_i)^2}{(\sum_{i=1}^N c_i + \sum_{i=1}^N \epsilon_i)^2},
\]
which goes to 0 if we assumes sensibly that as $N$ gets large, $\sum_{i=1}^N \epsilon_i$ goes to 0 at a faster rate than $\sum_{i=1}^N c_i$. Thus, given that we cannot distinguish between the common component and the noise term for a given set of data, the $q$-RSE is a better measure to gauge imputation performances.

\subsection{NYC taxi traffic}
\label{subsec:taxi_traffic}
We analyze a set of taxi traffic data in New York city in this example. The data includes all individual taxi rides operated by Yellow Taxi in New York City, published at

\url{https://www1.nyc.gov/site/tlc/about/tlc-trip-record-data.page}.

For simplicity, we only consider the rides within Manhattan Island, which comprises most of the data. The dataset contains 842 million trip records within the period of January 1, 2013 to December 31, 2022. Each trip record includes features such as pick-up and drop-off dates/times, pick-up and drop-off locations, trip distances, itemized fares, rate types, payment types, and driver-reported passenger counts. Our example here focuses on the pick-up and drop-off dates/times, and pick-up and drop-off locations.

To classify the pick-up and drop-off locations in Manhattan, they are coded according to 69 predefined zones in the dataset. Moreover, each day is divided into 24 hourly periods to represent the pick-up and drop-off times each day, with the first hourly period from 0 a.m. to 1 a.m. The total number of rides moving among the zones within each hour are recorded, yielding data $\cY_t\in\b{R}^{69\times 69\times 24}$ each day, where $y_{i_1,i_2,i_3,t}$ is the number of trips from zone $i_1$ to zone $i_2$ and the pick-up time is within the $i_3$-th hourly period on day $t$. We consider business days and non-business days separately, so that we will analyze two tensor time series. The business-day series is 2,519 days long, while the non-business-day series is 1,133 days long, within the period of January 1, 2013 to December 31, 2022.

\begin{figure}[ht!]
\begin{center}
\centerline{\includegraphics[width=0.95\columnwidth]{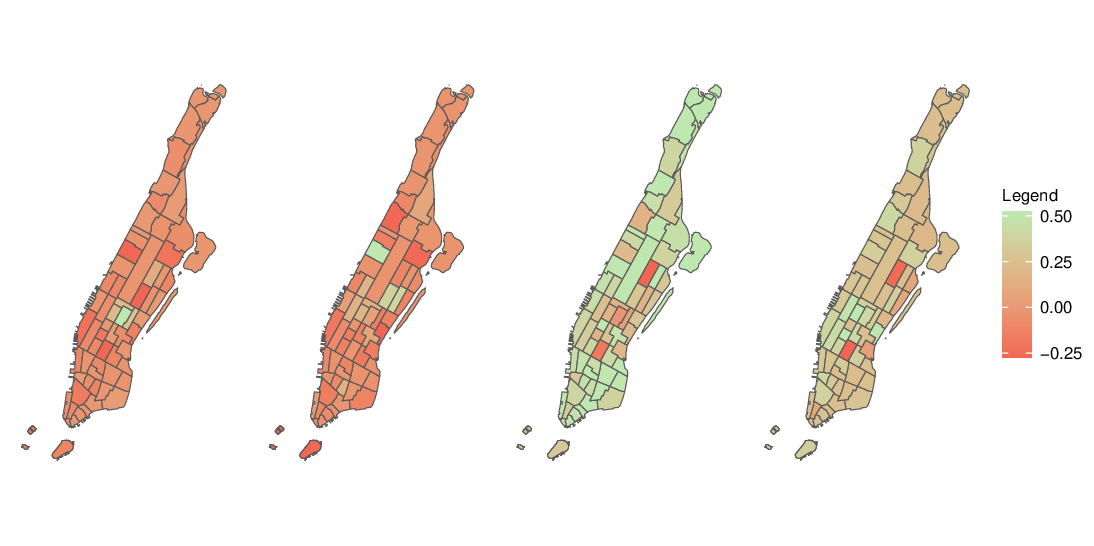}}
\caption{Estimated loading matrix $\wh\Q_1$ on four pick-up factors for business-day series. The data is block missing according to (M-iii).}
\label{Fig: taxi_A1_business}
\end{center}
\end{figure}

\begin{figure}[ht!]
\begin{center}
\centerline{\includegraphics[width=0.8\columnwidth]{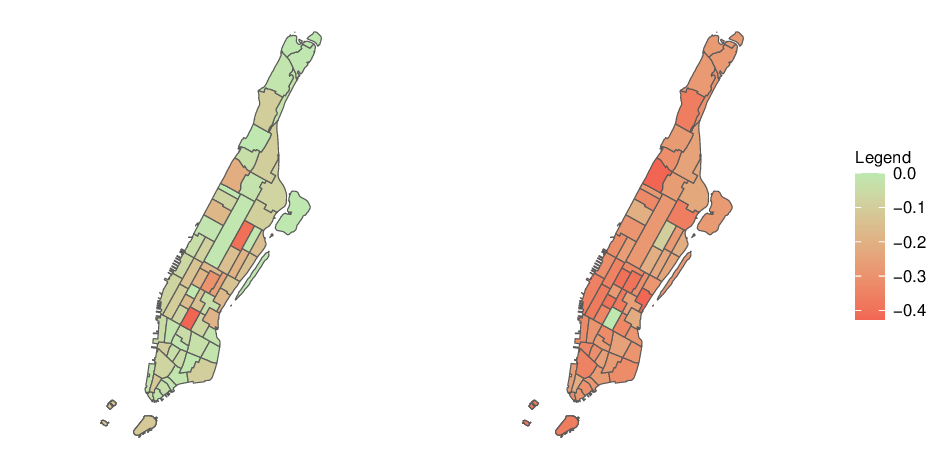}}
\caption{Estimated loading matrix $\wh\Q_2$ on two drop-off factors for business-day series. The data is block missing according to (M-iii).}
\label{Fig: taxi_A2_business}
\end{center}
\end{figure}

\begin{figure}[ht!]
\begin{center}
\centerline{\includegraphics[width=0.8\columnwidth]{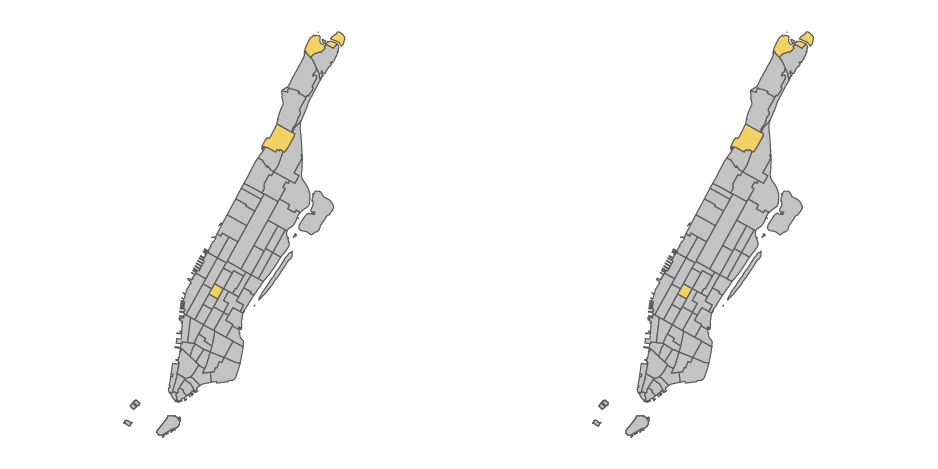}}
\caption{Testing zeros of $\Q_2$ on two drop-off factors for business-day series. Rows rejecting the null with $5\%$ significance level are in grey, otherwise in yellow.}
\label{Fig: taxi_test_A2_nonbusiness}
\end{center}
\end{figure}

As we are more interested in estimating the tensor factor structure with missingness and hence imputation, we demean the data before imposing different missing patterns, and investigate their effects. See descriptions of missing pattern (M-i) to (M-iii) in Section \ref{subsec:simulation}. We first estimate the rank of the core tensor using our proposed eigenvalue ratio estimator. Inspired by the simulation results in Section \ref{subsec:performance_of_number_of_factors}, one approach is to obtain the initial rank estimates $\wh r_k$ for $k\in[K]$, and re-estimate the number of factors based on the imputation using $\wh{r}_k+1$ as the number of factors along mode-$k$.

{
It appears that the re-estimated rank is $(1,1,1)$ for both the business-day and non-business-day series, and both are consistent under missing patterns (M-i), (M-ii), (M-iii) and no missing. This also coincides with iTIP-ER of \cite{Hanetal2022}. However, BCorTh of \cite{Chen_Lam} estimates $(4,2,2)$ for business-day series and $(3,2,2)$ for non-business-day series on fully observed data. Such discrepancy is due to the existence of weak factors in the NYC traffic data, and the rank $(1,1,1)$ is not sensible according to some exploratory studies in \cite{Chenetal2022}. Although our factor number estimator has taken weak factors into account, the requirements for the factor strength are more stringent than the requirements when using BCorTh. To see this, consider the case of $K=3$ and all factor strength are the same, then $\alpha_{k,j} =2/3$ satisfies the rate assumption in Theorem 5.2 of \cite{Chen_Lam} but fails our assumption in Theorem \ref{thm:number_of_factors}.
}

{
To cope with such weak factors, we use the rank $\wh{r}_k +r_\ast$ in the initial imputation for re-estimating the rank, where $r_\ast$ is a pre-specified integer characterising the maximum number of omitted factors. As discussed above, $r_\ast =1$ might not be sufficient if weak factors are likely to exist. We use $r_\ast =4$ here. As mentioned in Remark \ref{remark:2}, we may re-estimate the rank using BCorTh of \cite{Chenetal2022} based on the imputed data. The resulting estimated ranks for (M-i), (M-ii) and (M-iii) align with those using full observations.
}

For convenience, we use the estimated rank $(\wh{r}_1, \wh{r}_2, \wh{r}_3) =(4, 2, 2)$ in the rest of the analysis to estimate the factor loading matrices for different modes. We note the resulting heatmaps of the estimated loadings are almost the same under different missing patterns, so we only display one particular example here. Figures \ref{Fig: taxi_A1_business} and \ref{Fig: taxi_A2_business} show the heatmaps of $\wh\Q_1$ and $\wh\Q_2$ for the business-day series under missing pattern (M-iii).

\begin{figure}[ht!]
\begin{center}
\centerline{\includegraphics[width=0.95\columnwidth]{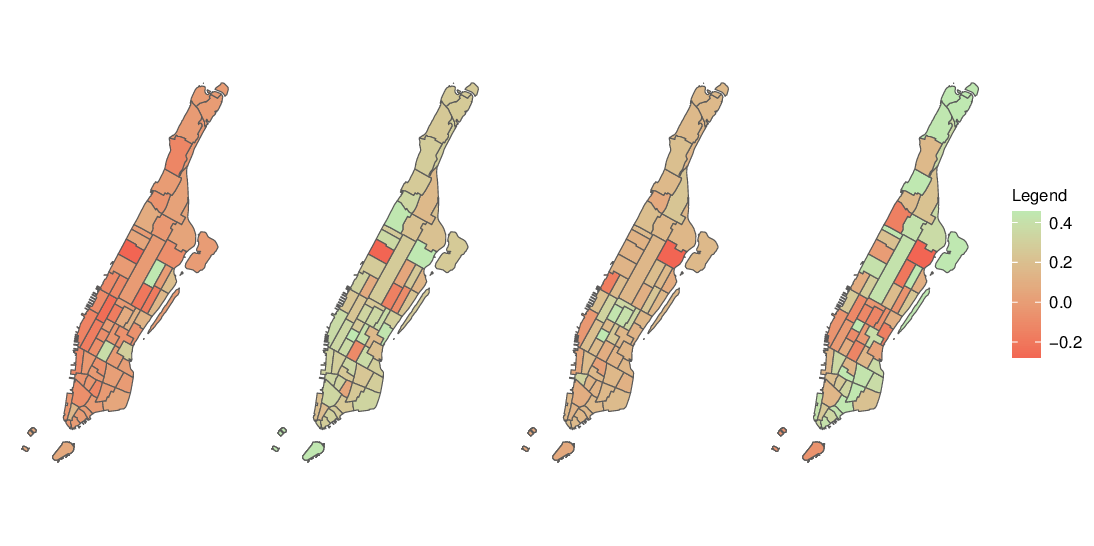}}
\caption{Estimated loading matrix $\wh\Q_1$ on four pick-up factors for non-business-day series. The data is block missing according to (M-iii).}
\label{Fig: taxi_A1_nonbusiness}
\end{center}
\end{figure}

\begin{figure}[ht!]
\begin{center}
\centerline{\includegraphics[width=0.8\columnwidth]{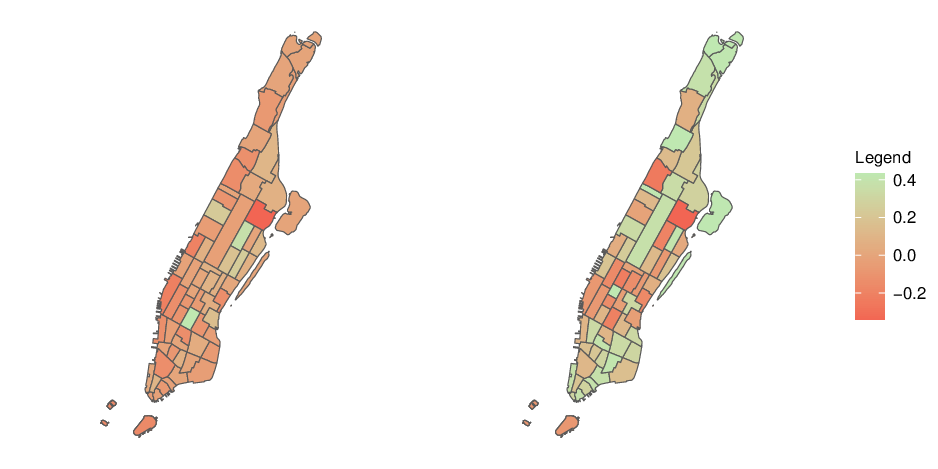}}
\caption{Estimated loading matrix $\wh\Q_2$ on two drop-off factors for non-business-day series. The data is block missing according to (M-iii).}
\label{Fig: taxi_A2_nonbusiness}
\end{center}
\end{figure}

The two heatmaps reveal the same pattern in Manhattan, that traffic is mainly within the downtown areas (the upper east and around the Empire State Building), whereas recreation park or amusement areas (the Central Park and upper west) load slightly on the factor structure. In contrast, dining and sports areas (especially East Harlem) stand out in the factor loadings for non-business-day series. See Figures \ref{Fig: taxi_A1_nonbusiness} and \ref{Fig: taxi_A2_nonbusiness}.
{
In particular, the first diagram in Figure \ref{Fig: taxi_A2_business} is revealing that taxi drop-offs in business days are concentrated in the midtown center. To further find out active drop-off areas for business days, we statistically test if the loading entry is zero by leveraging Theorem \ref{thm:asymp_normality_loadings} and Theorem \ref{thm:covariance_estimator}. We test for each $i\in[69]$,
\[
\c{H}_0: \Q_{2, i\cdot} =(0,\, 0)', \,\,\,\,\,
\c{H}_1: \text{otherwise}.
\]
Test statistics obtained by $\big\|( \wh\bSigma_{HAC} + \wh\bSigma_{HAC}^{\Delta} )^{-1/2} \wh\D_2 \wh\Q_{2, i\cdot} \big\|^2$ follow a $\chi_3^2$ distribution under the null. Using $5\%$ significance level, the inactive zones (i.e., rows failing to reject the null) are highlighted in yellow, as presented in Figure \ref{Fig: taxi_test_A2_nonbusiness}. As we would expect, the common zones for the two factors represent parks and fashion production areas in general.
}

The mode-3 loading matrix conveys the information of pick-up time, and the unveiled pattern is the same over different missing patterns. For better interpretation, we show $\wh\Q_3$ after a varimax rotation and scaling under missing pattern (M-iii) in Table \ref{tab: taxi_loading_hour_business} and Table \ref{tab: taxi_loading_hour_nonbusiness} for the business-day and non-business-day series respectively. The highlighted parts indicate the rush-hours in business and non-business days, which depict drastically different patterns. From Table \ref{tab: taxi_loading_hour_business}, we can see our first and second factors respectively capture the morning rush-hours (roughly from 9 a.m. to 10 a.m.) and the evening rush-hours (from 8 p.m. to 12 p.m.) for business days. Interestingly, we can see from the first factor that the taxi volume spikes at 4 p.m. which might be interpreted as an afternoon period. Table \ref{tab: taxi_loading_hour_nonbusiness} gives similar patterns, except that all active hours are delayed. We may notice the night life in business days almost end after 1 a.m., while the traffic remains busy even after 2 a.m. in non-business days.

\begin{table}[ht!]
\footnotesize
\setlength{\tabcolsep}{4pt}
\begin{center}
\begin{tabular}{l||ccccccccccccccccccccccccccc}
\hline $0_\text{am}$ & & 2 & & 4 & & 6 & & 8 & & 10 & & $12_\text{pm}$ & & 2 & & 4 & & 6 & & 8 & & 10 & & $12_\text{am}$ \\
\hline 1 & 2 & 1 & 0 & 0 & 0 & -1 & -3 & -7 & \red{-10} & -8 & -7 & -7 & -9 & -9 & -9 & \red{-10} & -8 & -9 & -8 & -5 & -2 & -2 & 0 & 1 \\
2 & -9 & -4 & -2 & -1 & -1 & 0 & 2 & 5 & 4 & 0 & -1 & -1 & 2 & 0 & 1 & 1 & -2 & 0 & -3 & \red{-10} & \red{-13} & \red{-12} & \red{-13} & \red{-13} \\
\hline
\end{tabular}
\end{center}
\caption{Estimated loading matrix $\wh\Q_3$ for business-day series, after varimax rotation and scaling. The data is block missing according to (M-iii). Magnitudes larger than 9 are highlighted in red.}
\label{tab: taxi_loading_hour_business}
\end{table}

\begin{table}[ht!]
\footnotesize
\setlength{\tabcolsep}{4pt}
\begin{center}
\begin{tabular}{l||ccccccccccccccccccccccccccc}
\hline $0_\text{am}$ & & 2 & & 4 & & 6 & & 8 & & 10 & & $12_\text{pm}$ & & 2 & & 4 & & 6 & & 8 & & 10 & & $12_\text{am}$ \\
\hline 1 & 0 & 2 & 3 & 2 & 1 & 0 & -1 & -1 & -4 & -6 & -9 & \red{-10} & \red{-11} & \red{-10} & \red{-10} & -9 & -8 & -8 & -8 & -6 & -3 & -2 & -2 & 0 \\
2 & \red{-10} & \red{-12} & \red{-11} & -9 & -6 & -2 & 0 & 0 & 1 & 2 & 2 & 2 & 1 & 0 & 0 & 0 & -1 & -1 & -3 & -6 & -9 & -9 & \red{-10} & \red{-11} \\
\hline
\end{tabular}
\end{center}
\caption{Estimated loading matrix $\wh\Q_3$ for non-business-day series, after varimax rotation and scaling. The data is block missing according to (M-iii). Magnitudes larger than 9 are highlighted in red.}
\label{tab: taxi_loading_hour_nonbusiness}
\end{table}

From the above discussion, it is easy to see that with our proposed imputation procedure, the tensor factor structure could be well studied even if we cannot observe the full data set. Nevertheless, Table \ref{tab: taxi_qRSE} reveals that different missing patterns have a certain degree of influences on the imputation accuracy, which is also suggested by previous simulation results. The $q$-RSE, {introduced in Section \ref{subsec:q-RSE}}, are reported for $q=5,10,20,100$, so that the overall performance is gauged accounting for large magnitudes of noise. Lastly, the analysis on the Manhattan Taxi data is for demonstration of our imputation method, and more sophisticated investigations are required to further understand the traffic pattern.

\begin{table}[ht!]
\begin{center}
\begin{tabular}{|c||ccc||ccc|}
\hline & \multicolumn{3}{c||}{business-day} & \multicolumn{3}{c|}{non-business-day}   \\ \hline
$q$   & \multicolumn{1}{c|}{(M-i)} & \multicolumn{1}{c|}{(M-ii)} & (M-iii) & \multicolumn{1}{c|}{(M-i)} & \multicolumn{1}{c|}{(M-ii)} & (M-iii) \\ \hline
5   & \multicolumn{1}{c|}{0.231}    & \multicolumn{1}{c|}{0.231}     & 0.217        & \multicolumn{1}{c|}{0.206}    & \multicolumn{1}{c|}{0.206}     & 0.226      \\ \hline
10   & \multicolumn{1}{c|}{0.225}   & \multicolumn{1}{c|}{0.225}    & 0.204      & \multicolumn{1}{c|}{0.213}   & \multicolumn{1}{c|}{0.213}    & 0.221     \\ \hline
20  & \multicolumn{1}{c|}{0.225}   & \multicolumn{1}{c|}{0.225}    & 0.195     & \multicolumn{1}{c|}{0.225}   & \multicolumn{1}{c|}{0.225}    & 0.228     \\ \hline
100 & \multicolumn{1}{c|}{0.234}   & \multicolumn{1}{c|}{0.234}    & 0.186 & \multicolumn{1}{c|}{0.247}   & \multicolumn{1}{c|}{0.247}    & 0.247     \\ \hline
\end{tabular}
\end{center}
\caption {$q$-RSE for business-day and non-business-day series under different missing patterns, with $q=5,10,20,100$.}
\label{tab: taxi_qRSE}
\end{table}

\newpage

\section{Appendix: Proof of all theorems and propositions}\label{sec: Appendix_proof}

From Section \ref{subsec:est_Q}, $\wh\Q_k$ contains the eigenvectors corresponding to the first $r_k$ largest eigenvalues of $\wh\S_k$. Hence with $\wh\D_k$ an $r_k\times r_k$ diagonal matrix containing all the eigenvalues of $\wh\S_k$ (WLOG from the largest on the top-left element to the smallest on the bottom right element), we have $\wh\S_k\wh\Q_k = \wh\Q_k\wh\D_k$, so that
\begin{equation}\label{eqn: Q_pre}
    \wh\Q_k = \wh\S_k\wh\Q_k\wh\D_k^{-1}.
\end{equation}
To simplify notations, hereafter we fix $k$ and only focus on the mode-$k$ unfolded data. Define
\begin{equation}
\label{eqn: notation_short}
    \begin{split}
    \wh\D &:= \wh\D_k, \; \Y_t := \mat{k}{\cY_t}, \; \wh\S := \wh\S_k, \;
    \psi_{ij,h} := \psi_{k,ij,h}, \; \Q := \Q_k, \\
    \bLambda &:= \bLambda_k, \;
    \F_{Z,t} := \mat{k}{\cF_{Z,t} }, \; \E_t := \mat{k}{\cE_t}, \; \H_{j} := \H_{k,j}, \; \H^a := \H_{k}^a,
    \end{split}
\end{equation}
where $\H_{k,j}$ and $\H_{k}^a$ are defined in (\ref{eqn:H_j}) and (\ref{eqn:H_k(all)}) respectively,   and similarly to all respective hat versions of the above.

Before proving any theorems, we present and prove the following Proposition first.
\begin{proposition}\label{Prop:assumption_implications}
Let Assumption (E1), (E2) and (F1) hold. Then
\begin{itemize}
\item[1.] there exists a constant $c>0$ so that
for any $k\in[K], t\in[T], i_k\in[d_k]$ and $h\in[\dmk]$, we have $\b{E}\c{E}_{t,i_1,\dots,i_K}=0$, $\b{E}\c{E}_{t,i_1,\dots,i_K}^4 \leq c$, and
\begin{align*}
  &\sum_{j=1}^{d_k}\sum_{l=1}^{\dmk}
    \Big|\b{E}[\textnormal{mat}_k(\c{E}_t)_{ih}
    \textnormal{mat}_k(\c{E}_t)_{jl}]\Big|
    \leq c,\\
    &\sum_{l=1}^{d_{\text{-}k}}
\sum_{s\in\psi_{k,ij,l}}
    \bigg| \textnormal{cov}\Big(
\textnormal{mat}_k(\c{E}_t)_{ih}
\textnormal{mat}_k(\c{E}_t)_{jh},
\textnormal{mat}_k(\c{E}_s)_{il}
\textnormal{mat}_k(\c{E}_s)_{jl}
\Big) \bigg| \leq c;
\end{align*}

\item[2.] there exists a constant $c>0$ so that for any $k\in[K],i,j\in[d_k]$, and any deterministic vectors $\bf{u}\in\b{R}^{r_k}$ and $\bf{v}\in\b{R}^{\rmk}$ with constant magnitudes,
\[
\b{E}\Bigg(
        \frac{1}{\dmk^{1/2}}\sum_{h=1}^{\dmk}
        \frac{1}{|\psi_{k,ij,h}|^{1/2}}\sum_{t\in\psi_{k,ij,h}}
        \textnormal{mat}_k(\c{E}_{t})_{jh}
        \bf{u}' \textnormal{mat}_k(\c{F}_{t})
        \bf{v}
    \Bigg)^2 \leq c;
\]

\item[3.] for any $k\in[K],i,j\in[d_k],h\in[d_{\text{-}k}]$,
\[
\frac{1}{|\psi_{k,ij,h}|}\sum_{t\in\psi_{k,ij,h}}
\textnormal{mat}_k(\c{F}_t)
\textnormal{mat}_k(\c{F}_t)', \;\;\; \frac{1}{T}\sum_{t=1}^T
\textnormal{mat}_k(\c{F}_t)
\textnormal{mat}_k(\c{F}_t)'
\xrightarrow{p} \bf{\Sigma}_k := \rmk\I_{r_k},
\]
with the number of factors $r_k$ fixed as $\min\{T,d_1,\dots,d_K\}\to\infty$. For each $t\in[T]$, all elements in $\c{F}_t$ are independent of each other, with mean 0 and unit variance.
\end{itemize}
\end{proposition}
Our consistency results in Theorem \ref{thm:factor_loading_consistency} can be proved assuming the three implied results from Proposition \ref{Prop:assumption_implications}, on top of Assumption (O1), (M1), (L1) and (R1). Result 1 from Proposition \ref{Prop:assumption_implications} can be a stand alone assumption on the weak correlation of the noise $\c{E}_t$ across different dimensions and times, while result 2 can be on the weak dependence between the factor $\c{F}_t$ and the noise $\c{E}_t$. Finally, result 3 can be a stand alone assumption on the factors $\cF_t$.

\noindent\textbf{\textit{Proof of Proposition \ref{Prop:assumption_implications}.}}
We have $\b{E}(\c{E}_{t}) = 0$ from Assumption (E1). Next we want to show that for any $k\in[K], t\in[T]$ and $i_k\in[d_k]$, $\b{E}\c{E}_{t,i_1,\dots,i_K}^4$ is bounded uniformly. From (\ref{eqn: NE1}), each entry in $\c{E}_t$ is a sum of two parts: a linear combination of the elements in $\c{F}_{e,t}$, and the corresponding entry in $\bepsilon_t$. By Assumption (E2), we have
\begin{equation*}
\begin{split}
    \b{E}[(\bepsilon_t)_{i_1,\dots,i_K}^4]
    &= \b{E}\bigg[\bigg(\sum_{q\geq 0}a_{\epsilon,q}(\c{X}_{\epsilon,t-q})
    _{i_1,\dots,i_K}\bigg)^4\bigg]
    \leq
    \bigg(\sum_{q\geq 0}|a_{\epsilon,q}|\bigg)^4
    \sup_{t}\b{E}\Big[(\c{X}_{\epsilon,t})_{i_1,\dots,i_K}^4\Big]\leq C,
\end{split}
\end{equation*}
where $C > 0$ is a generic constant. It holds similarly that $\b{E}[(\c{F}_{e,t})_{i_1,\dots,i_K}^4]\leq C$ uniformly on all indices. With this, defining $\A^{\circ m} := \A\circ\cdots\circ\A$ (elementwise $m$-th power),
\begin{align*}
  \|\b{E}\cE_t^{\circ 4}\|_{\max} &= \|\b{E}[\mat{k}{\cE_t}]^{\circ 4}\|_{\max}
  \leq 8(\|\b{E}[\A_{e,k}\mat{k}{\c{F}_{e,t}}\Aemk']^{\circ 4}\|_{\max} + \|\bSigma_{\epsilon}^{\circ 4}\|_{\max}\cdot \|\b{E}\bepsilon_t^{\circ 4}\|_{\max})\\
  &\leq 8(\norm{\A_{e,k}}_{\infty}^4\cdot\norm{\Aemk}_{\infty}^4\cdot\norm{\b{E}\c{F}_{e,t}^{\circ 4}}_{\max} + \|\bSigma_{\epsilon}^{\circ 4}\|_{\max}\cdot \|\b{E}\bepsilon_t^{\circ 4}\|_{\max})\\
  &=8(\prod_{j=1}^K\norm{\A_{e,j}}_{\infty}^4\cdot\norm{\b{E}\c{F}_{e,t}^{\circ 4}}_{\max}+\|\bSigma_{\epsilon}^{\circ 4}\|_{\max}\cdot \|\b{E}\bepsilon_t^{\circ 4}\|_{\max})
  \leq C,
\end{align*}
where $C>0$ is again a generic constant, and we used Assumption (E1) in the last line and the fact that $r_j$ is a constant for $j\in[K]$. This is equivalent to $\b{E}\c{E}_{t,i_1,\dots,i_K}^4\leq c$ for some constant $c$.

With (\ref{eqn: NE1}) in Assumption (E1), we have \[\mat{k}{\c{E}_t}=\A_{e,k}\mat{k}{\c{F}_{e,t}}\A_{e,\text{-}k}' + \mat{k}{\bSigma_{\epsilon}}\circ \mat{k}{\bepsilon_t},\]
 where $\A_{e,\text{-}k}:=\A_{e,K}\otimes\cdots\otimes\A_{e,k+1}\otimes\A_{e,k-1}\otimes\cdots\otimes\A_{e,1}$. Each mode-$k$ noise fibre $\bf{e}_{t,k,l}$ for $l\in[d_{\text{-}k}]$ can then be decomposed as
\begin{equation}
    \label{eqn: noise_fibre}
    \bf{e}_{t,k,l}:=
    \A_{e,k}\mat{k}{\c{F}_{e,t}}\bf{A}_{e,\text{-}k,l\cdot}+\bf{\Sigma}_{\epsilon,k,l}^{1/2}\bepsilon_{t,k,l},
\end{equation}
where $\bf{\Sigma}_{\epsilon,k,l} = \diag((\mat{k}{\bSigma_{\epsilon}})_{\cdot l}(\mat{k}{\bSigma_{\epsilon}})_{\cdot l}')$, and $\bepsilon_{t,k,l}$ contains independent elements each with mean $0$ and variance $1$.

Given $h\in[d_{\text{-}k}], i\in[d_k]$, from (\ref{eqn: noise_fibre}) and Assumption (E1) and (E2), we have
\begin{equation*}
\begin{split}
    \sum_{l\neq h}^{\dmk}\sum_{j=1}^{d_k}\Big|\b{E}[\mat{k}{\c{E}_t}_{ih}\mat{k}{\c{E}_t}_{jl}]\Big| &\leq
    \|\bf{A}_{e,\text{-}k,h\cdot}\|
    \|\bf{A}_{e,\text{-}k,l\cdot}\|
    \|\A_{e,k}\|_1\|\A_{e,k}\|_\infty
    = O(1).
\end{split}
\end{equation*}
Moreover,
\begin{equation*}
\begin{split}
    \sum_{j=1}^{d_k}\Big|\b{E}[\mat{k}{\c{E}_t}_{ih}\mat{k}{\c{E}_t}_{jh}]\Big| &\leq
    \|\text{cov}(\bf{e}_{t,k, h},\bf{e}_{t,k,h})\|_1 \\
    &\leq
    \|\bf{A}_{e,\text{-}k,h\cdot}\|^2
    \|\A_{e,k}\|_1\|\A_{e,k}\|_\infty
    +\|\bf{\Sigma}_{\epsilon,k,h}\|_1
    = O(1),
\end{split}
\end{equation*}
where the last equality is from Assumption (E1).

To finish the proof of the first result in the Proposition, fix indices $k,t,i,j,h$. From Assumption (E2) and (\ref{eqn: noise_fibre}), we have
\begin{equation}
\begin{split}
\label{eqn: noise_entry_decomp}
    [\mat{k}{\c{E}_t}]_{il}
    &=\sum_{q\geq 0}a_{e,q}
    \A_{e,k,i\cdot}'\mat{k}{\c{X}_{e,t-q}}\bf{A}_{e,\text{-}k,l\cdot}+
    [\mat{k}{\bSigma_\epsilon}]_{il}\sum_{q\geq 0}a_{\epsilon,q}
    \mat{k}{\c{X}_{\epsilon,t-q}}_{il}.
\end{split}
\end{equation}
Hence when $l\neq h$, from the independence between $\{\c{X}_{e,t}\}$ and $\{\c{X}_{\epsilon,t}\}$ and the independence of the elements within $\{\c{X}_{\epsilon,t}\}$ in Assumption (E2), we have for any $s\in\psi_{k,ij,l}$,
\begin{equation}
\begin{split}
\label{eqn: cov_eeee}
    &\hspace{5mm}
    \text{cov}\Big(
    \mat{k}{\c{E}_t}_{ih}
    \mat{k}{\c{E}_t}_{jh},
    \mat{k}{\c{E}_s}_{il}
    \mat{k}{\c{E}_s}_{jl}
    \Big) \\
    &=
    \text{cov}\Bigg(
    \Big(
    \sum_{q\geq 0}a_{e,q}
    \A_{e,k,i\cdot}'\mat{k}{\c{X}_{e,t-q}}\bf{A}_{e,\text{-}k,h\cdot}
    \Big)\Big(
    \sum_{q\geq 0}a_{e,q}
    \A_{e,k,j\cdot}'\mat{k}{\c{X}_{e,t-q}}\bf{A}_{e,\text{-}k,h\cdot}
    \Big),\\
    &
    \Big(
    \sum_{q\geq 0}a_{e,q}
    \A_{e,k,i\cdot}'\mat{k}{\c{X}_{e,s-q}}\bf{A}_{e,\text{-}k,l\cdot}
    \Big)\Big(
    \sum_{q\geq 0}a_{e,q}
    \A_{e,k,j\cdot}'\mat{k}{\c{X}_{e,s-q}}\bf{A}_{e,\text{-}k,l\cdot}
    \Big) \Bigg).
\end{split}
\end{equation}
By (E2) again, all mixed covariance terms are zero except for $\text{cov}\big(\mat{k}{\c{X}_{e,t-q}}_{nm}^2,\mat{k}{\c{X}_{e,t-q}}_{nm}^2\big)$ for all $q\geq 0, n\in[r_{e,k}],m\in[r_{e,\text{-}k}]$, with coefficient
$
    a_{e,q}^2a_{e,q+|t-s|}^2
    A_{e,k,in}^2
    A_{e,k,jn}^2
    A_{e,\text{-}k,hm}^2
    A_{e,\text{-}k,lm}^2.
$
Thus we have
\begin{equation*}
\begin{split}
    &\hspace{5mm}
    \sum_{l=1,l\neq h}^{d_{\text{-}k}}
    \bigg|\text{cov}\Big(
    \mat{k}{\c{E}_t}_{ih}
    \mat{k}{\c{E}_t}_{jh},
    \mat{k}{\c{E}_s}_{il}
    \mat{k}{\c{E}_s}_{jl}
    \Big)\bigg| \\
    &=
    \sum_{l=1,l\neq h}^{d_{\text{-}k}}
    \bigg|\sum_{n=1}^{r_{e,k}}
    \sum_{m=1}^{r_{e,\text{-}k}}
    \sum_{q\geq 0}
    a_{e,q}^2a_{e,q+|t-s|}^2
    A_{e,k,in}^2
    A_{e,k,jn}^2
    A_{e,\text{-}k,hm}^2
    A_{e,\text{-}k,lm}^2
    \cdot \text{cov}\big(\mat{k}{\c{X}_{e,t-q}}_{nm}^2,\mat{k}{\c{X}_{e,t-q}}_{nm}^2\big)\bigg| \\
    &=
    O(1)\cdot
    \sum_{q\geq 0}
    a_{e,q}^2a_{e,q+|t-s|}^2
    \bigg(
    \sum_{n=1}^{r_{e,k}}
    A_{e,k,in}^2
    A_{e,k,jn}^2 \bigg)\bigg(
    \sum_{m=1}^{r_{e,\text{-}k}}
    A_{e,\text{-}k,hm}^2
    \sum_{l=1,l\neq h}^{d_{\text{-}k}}
    A_{e,\text{-}k,lm}^2 \bigg)
    =
    \sum_{q\geq 0}
    O(a_{e,q}^2a_{e,q+|t-s|}^2) ,
\end{split}
\end{equation*}
where we use Assumption (E2) in the second last equality, and (E1) in the last. Consequently,
\begin{equation*}
\begin{split}
    &\hspace{5mm}
    \sum_{l=1,l\neq h}^{d_{\text{-}k}}
    \sum_{s\in\psi_{k,ij,l}}
    \bigg|\text{cov}\Big(
    \mat{k}{\c{E}_t}_{ih}
    \mat{k}{\c{E}_t}_{jh},
    \mat{k}{\c{E}_s}_{il}
    \mat{k}{\c{E}_s}_{jl}
    \Big)\bigg|
    =
    \sum_{q\geq 0}\sum_{s=1}^T
    O(a_{e,q}^2a_{e,q+|t-s|}^2) =
    O(1),
\end{split}
\end{equation*}
where the last equality uses Assumption (E2). Now consider lastly $l=h$.
All arguments starting from (\ref{eqn: noise_entry_decomp}) follow exactly, except the following term is added in (\ref{eqn: cov_eeee}):
\begin{equation*}
\begin{split}
    &\hspace{5mm}
    \sum_{q\geq 0}a_{\epsilon,q}^2
    a_{\epsilon,q+|t-s|}^2
    \Sigma_{\epsilon,k,h,ii}
    \Sigma_{\epsilon,k,h,jj}
    \cdot \text{cov}\Big(
    \mat{k}{\c{X}_{\epsilon,t-q}}_{ih}
    \mat{k}{\c{X}_{\epsilon,t-q}}_{jh},
    \mat{k}{\c{X}_{\epsilon,s-q}}_{ih}
    \mat{k}{\c{X}_{\epsilon,s-q}}_{jh}
    \Big)\\
    &=
    O(1)\cdot
    \sum_{q\geq 0}a_{\epsilon,q}^2
    a_{\epsilon,q+|t-s|}^2
    = O(1).
\end{split}
\end{equation*}
where we used again Assumptions (E2) in the last line. Finally,
\begin{equation*}
\begin{split}
    &\hspace{5mm}
    \sum_{s\in\psi_{k,ij,h}}
    \bigg| \text{cov}\Big(
    \mat{k}{\c{E}_t}_{ih}
    \mat{k}{\c{E}_t}_{jh},
    \mat{k}{\c{E}_s}_{ih}
    \mat{k}{\c{E}_s}_{jh}
    \Big) \bigg|
    = O(1).
\end{split}
\end{equation*}
This completes the proof of result 1 in the Proposition.

To prove the second result, fix $k\in[K], i,j\in[d_k]$ and  deterministic vectors $\bf{u}\in\b{R}^{r_k}$ and $\bf{v}\in\b{R}^{\rmk}$ with $\|\bf{u}\|,\|\bf{v}\|=O(1)$. Note that $\b{E}[\mat{k}{\c{F}_{t}}\bf{v}\bf{v}'\mat{k}{\c{F}_{s}}']=\bf{v}'\bf{v}(r_{\text{-}k}\sum_{q\geq 0}a_{f,q}a_{f,q+|t-s|})\bf{I}_{r_k}$, as the series $\{\c{X}_{f,t}\}$ has i.i.d. elements from Assumption (F1). Similarly, from (\ref{eqn: noise_fibre}) and Assumption (E1) and (E2),
\begin{equation*}
\begin{split}
    &\hspace{5mm}
    \text{cov}(\mat{k}{\c{E}_{t}}_{jh},
    \mat{k}{\c{E}_{s}}_{jl}) \\
    &=
    \b{E}[\A_{e,k,j\cdot}'\mat{k}{\c{F}_{e,t}}\A_{e,\text{-}k,h\cdot}\A_{e,\text{-}k,l\cdot}'\mat{k}{\c{F}_{e,s}}'\A_{e,k,j\cdot}]
    +
    \b{E}[\bepsilon_{t,k,h}'
    (\Sigma_{\epsilon,k,h,j\cdot}
    \Sigma_{\epsilon,k,l,j\cdot}')^{1/2}
    \bepsilon_{s,k,l}]\\
    &=
    \A_{e,\text{-}k,l\cdot}'\A_{e,\text{-}k,h\cdot}
    \cdot\|\A_{e,k,j\cdot}\|^2
    \cdot \sum_{q\geq 0}a_{e,q}a_{e,q+|t-s|} +
    \b{1}_{\{h=l\}}\cdot \Sigma_{\epsilon,k,h,jj}\sum_{q\geq 0}a_{\epsilon,q}a_{\epsilon,q+|t-s|}.
\end{split}
\end{equation*}
Hence if we fix $h\in[\dmk], t\in\psi_{k,ij,h}$, then together with Assumption (E2), we have
\begin{equation}
\label{proof: prop3_main}
\begin{split}
    &\hspace{5mm}
    \sum_{l=1}^{\dmk}
    \sum_{s\in\psi_{k,ij,l}}
    \frac{1}{|\psi_{k,ij,l}|}
    \cdot \b{E}
    \Big[\mat{k}{\c{E}_{t}}_{jh}
    \bf{u}'\mat{k}{\c{F}_{t}}\bf{v}
    \cdot \mat{k}{\c{E}_{s}}_{jl}
    \bf{v}'\mat{k}{\c{F}_{s}}'
    \bf{u}\Big] \\
    &=
    \sum_{l=1}^{\dmk}
    \sum_{s\in\psi_{k,ij,l}}
    \frac{1}{|\psi_{k,ij,l}|}
    \cdot \text{cov}(\mat{k}{\c{E}_{t}}_{jh},
    \mat{k}{\c{E}_{s}}_{jl})
    \cdot \b{E}\Big[
    \bf{u}'\mat{k}{\c{F}_{t}}\bf{v}
    \bf{v}'\mat{k}{\c{F}_{s}}'\bf{u}
    \Big] \\
    &=
    \sum_{l=1}^{\dmk}
    \frac{1}{|\psi_{k,ij,l}|}\bigg[
    O(\A_{e,\text{-}k,l\cdot}'\A_{e,\text{-}k,h\cdot} \cdot\|\A_{e,k,j\cdot}\|^2) \cdot
    \sum_{q\geq 0}\sum_{p\geq 0}\sum_{s\in\psi_{k,ij,l}}
    a_{e,q}a_{e,q+|t-s|} a_{f,p}a_{f,p+|t-s|} \\
    & + O(\b{1}_{\{h=l\}}\cdot \Sigma_{\epsilon,k,h,jj})\cdot
    \sum_{q\geq 0}\sum_{p\geq 0}\sum_{s\in\psi_{k,ij,l}}a_{\epsilon,q}a_{\epsilon,q+|t-s|}a_{f,p}a_{f,p+|t-s|} \bigg] \\
    &=
    \sum_{l=1}^{\dmk}
    \frac{1}{|\psi_{k,ij,l}|} \cdot O(\A_{e,\text{-}k,l\cdot}'\A_{e,\text{-}k,h\cdot} \cdot\|\A_{e,k,j\cdot}\|^2 +
    \b{1}_{\{h=l\}}\cdot \Sigma_{\epsilon,k,h,jj})
    = O\Big(\frac{1}{T}\Big),
\end{split}
\end{equation}
where for the second last equality, we argue for the first term in the second last line only, as the second term could be shown similarly:
\begin{equation*}
\begin{split}
    &\sum_{q\geq 0}
    \sum_{p\geq 0}\sum_{s\in\psi_{k,ij,h}}
    a_{e,q}a_{e,q+|t-s|}
    a_{f,p}a_{f,p+|t-s|} =
    \sum_{q\geq 0}
    \sum_{p\geq 0}
    a_{e,q}a_{f,p}
    \sum_{s\in\psi_{k,ij,h}}
    a_{e,q+|t-s|}a_{f,p+|t-s|}\\
    &\leq
    \sum_{q\geq 0}
    \sum_{p\geq 0}
    |a_{e,q}|\cdot|a_{f,p}| \cdot
    \Bigg(\sum_{s\in\psi_{k,ij,h}}
    a_{e,q+|t-s|}^2\Bigg)^{1/2}
    \Bigg(\sum_{s\in\psi_{k,ij,h}}
    a_{f,p+|t-s|}^2\Bigg)^{1/2} \leq \sum_{q\geq 0}
    \sum_{p\geq 0}
    |a_{e,q}|\cdot|a_{f,p}| \leq c^2,
\end{split}
\end{equation*}
where the constant $c$ is from Assumptions (F1) and (E2). Finally,
\begin{equation*}
\begin{split}
    &\hspace{5mm}
    \b{E}\Bigg(
        \sum_{h=1}^{\dmk}
        \sum_{t\in\psi_{k,ij,h}}
        \frac{1}{\dmk \cdot |\psi_{k,ij,h}|}
        \mat{k}{\c{E}_{t}}_{jh}
        \bf{u}' \mat{k}{\c{F}_{t}}
        \bf{v}
    \Bigg)^2 \\
    &=
    \frac{1}{\dmk^2}
    \sum_{h=1}^{\dmk}\sum_{l=1}^{\dmk}
    \sum_{t\in\psi_{k,ij,h}}
    \sum_{s\in\psi_{k,ij,l}}
    \frac{1}{|\psi_{k,ij,h}| \cdot |\psi_{k,ij,l}|}
    \b{E} \Big[\mat{k}{\c{E}_{t}}_{jh}
    \bf{u}'\mat{k}{\c{F}_{t}}\bf{v}
    \cdot \mat{k}{\c{E}_{s}}_{jl}
    \bf{v}'\mat{k}{\c{F}_{s}}'
    \bf{u}\Big] \\
    &=
    \frac{1}{\dmk^2 T}
    \sum_{h=1}^{\dmk}\sum_{t\in\psi_{k,ij,h}}
    O\Big(\frac{1}{T}\Big) = O\Big(\frac{1}{\dmk T}\Big),
\end{split}
\end{equation*}
which then implies result 2 of the Proposition.

Finally, we prove result 3 of the Proposition.
From Assumption (F1), we have $\b{E}[\c{F}_t]={0}$. Next, for any $t\in[T]$, it is direct from Assumption (F1) that all elements in $\c{F}_t$ are independent. Moreover,
\begin{equation*}
\begin{split}
    \b{E}[\mat{k}{\c{F}_t}\mat{k}{\c{F}_t}']
    &= \b{E}\bigg[
    \bigg(\sum_{q\geq 0}a_{f,q}\mat{k}{\c{X}_{f,t-q}}\bigg)
    \bigg(\sum_{q\geq 0}a_{f,q}\mat{k}{\c{X}_{f,t-q}}'\bigg)\bigg] \\
    &=
    \sum_{q\geq 0}a_{f,q}^2
    \b{E}[\mat{k}{\c{X}_{f,t-q}}
    \mat{k}{\c{X}_{f,t-q}}']
    =
    \bigg(\sum_{q\geq 0}a_{f,q}^2\bigg)
    \cdot r_{\text{-}k}\bf{I}_{r_k}
    = r_{\text{-}k}\bf{I}_{r_k},
\end{split}
\end{equation*}
where we use Assumption (F1) in the last line. To complete the proof, WLOG consider
the variance of the $j$-th diagonal element of $\mat{k}{\cF_t}\mat{k}{\cF_t}'$. From Assumption (F1), we have
\begin{equation*}
\begin{split}
    &\text{Var}\bigg(\frac{1}{T}
    \sum_{t=1}^T [\mat{k}{\cF_t}]_{j\cdot}'[\mat{k}{\cF_t}]_{j\cdot}\bigg) =
    \frac{1}{T^2}\text{Var}\bigg(
    \sum_{t=1}^T
    \bigg[\sum_{q\geq 0}a_{f,q}[\mat{k}{\c{X}_{f,t-q}}]_{j\cdot}'\bigg]
    \bigg[\sum_{q\geq 0}a_{f,q}[\mat{k}{\c{X}_{f,t-q}}]_{j\cdot}\bigg]\bigg)
    \\ &=
    \frac{1}{T^2}\sum_{t=1}^T
    \sum_{s=1}^T\text{cov}\bigg(
    \bigg[\sum_{q\geq 0}a_{f,q}[\mat{k}{\c{X}_{f,t-q}}]_{j\cdot}'\bigg]
    \bigg[\sum_{q\geq 0}a_{f,q}[\mat{k}{\c{X}_{f,t-q}}]_{j\cdot}\bigg],\\
    &\hspace{1.2in} \bigg[\sum_{q\geq 0}a_{f,q}[\mat{k}{\c{X}_{f,s-q}}]_{j\cdot}'\bigg]
    \bigg[\sum_{q\geq 0}a_{f,q}[\mat{k}{\c{X}_{f,s-q}}]_{j\cdot}\bigg]\bigg)
    \\ &=
    \frac{1}{T^2}\sum_{t=1}^T
    \sum_{q\geq 0}
    a_{f,q}^4\text{Var}(
    [\mat{k}{\c{X}_{f,t-q}}]_{j\cdot}'
    [\mat{k}{\c{X}_{f,t-q}}]_{j\cdot}) \\
    &+
    \frac{1}{T^2}\sum_{t=1}^T
    \sum_{q\geq 0}\sum_{p\neq q}
    a_{f,q}^2 a_{f,p}^2
    \text{Var}([\mat{k}{\c{X}_{f,t-q}}]_{j\cdot}' [\mat{k}{\c{X}_{f,t-p}}]_{j\cdot}) \\
    &=
    \frac{\rmk}{T^2}\sum_{t=1}^T
    \sum_{q\geq 0}
    a_{f,q}^4\text{Var}(
    [\mat{k}{\c{X}_{f,t-q}}]_{j1}^2) +
    \frac{\rmk}{T^2}\sum_{t=1}^T
    \sum_{q\geq 0}\sum_{p\neq q}
    a_{f,q}^2 a_{f,p}^2 \\
    &=
    O(1) \cdot \frac{1}{T^2}\sum_{t=1}^T
    \sum_{q\geq 0}a_{f,q}^4 +
    O(1) \cdot \frac{1}{T^2}\sum_{t=1}^T
    \sum_{q\geq 0}\sum_{p\neq q}a_{f,q}^2a_{f,p}^2
    =
    O(1) \cdot \frac{1}{T^2}\sum_{t=1}^T
    (\sum_{q\geq 0}a_{f,q}^2)^2 =
    O(\frac{1}{T})=o(1),
\end{split}
\end{equation*}
where the third equality uses the independence in Assumption (E2). This completes the proof of result 3, and hence the Proposition. $\square$


To prove Theorem \ref{thm:factor_loading_consistency}, we are going to present some lemmas and prove them first. From (\ref{eqn: Q_pre}),
\begin{equation}
\label{eqn: Q_simplify}
\begin{split}
    \wh{\bf{Q}}_{j\cdot} =& \wh{\bf{D}}^{-1}\sum_{i=1}^{d_k} \wh{\bf{Q}}_{i\cdot}\wh{S}_{ij}
    = \wh{\bf{D}}^{-1}\sum_{i=1}^{d_k} \wh{\bf{Q}}_{i\cdot}
    \sum_{h=1}^{d_{\text{-}k}}
        \frac{1}{|\psi_{ij,h}|}
        \sum_{t\in\psi_{ij,h}}
        Y_{t,ih}Y_{t,jh}.
\end{split}
\end{equation}
With the notations in (\ref{eqn: notation_short}), (\ref{eqn: tenfac-all}) can be written as $\Y_t = \Q\F_{Z,t}\bLambda' + \E_t$, and hence for $i,j\in[d_k],h\in[d_{\text{-}k}]$,
\begin{equation*}
\begin{split}
    Y_{t,ih} =&
    \Big(
    \sum_{n=1}^{r_k} \sum_{m=1}^{r_{\text{-}k}} Q_{in}\Lambda_{hm} F_{Z,t,nm}
    \Big)+  E_{t,ih}
    =
    \Q_{i\cdot}'
    \Big(
    \sum_{m=1}^{r_{\text{-}k}}\Lambda_{hm} \F_{Z,t,\cdot m}
    \Big)+ E_{t,ih}\\
    =&
    \Big(
    \sum_{m=1}^{r_{\text{-}k}}\Lambda_{hm} \F_{Z,t,\cdot m}
    \Big)'\Q_{i\cdot}+ E_{t,ih}.
\end{split}
\end{equation*}
Hence the product $Y_{t,ih}Y_{t,jh}$ in (\ref{eqn: Q_simplify}) can be written as
\begin{equation}
\label{eqn: Q_simplify2}
\begin{split}
    Y_{t,ih}Y_{t,jh} &=
    \Big(
    \sum_{m=1}^{r_{\text{-}k}}\Lambda_{hm} \bf{F}_{Z,t,\cdot m}
    \Big)'\bf{Q}_{i\cdot}\Big(
    \sum_{m=1}^{r_{\text{-}k}}\Lambda_{hm} \bf{F}_{Z,t,\cdot m}
    \Big)'\bf{Q}_{j\cdot} +E_{t,ih} E_{t,jh}
    \\
    &+
    E_{t,jh} \Big(
    \sum_{m=1}^{r_{\text{-}k}}\Lambda_{hm} \bf{F}_{Z,t,\cdot m}
    \Big)'\bf{Q}_{i\cdot}
    + E_{t,ih}\Big(
    \sum_{m=1}^{r_{\text{-}k}}\Lambda_{hm} \bf{F}_{Z,t,\cdot m}
    \Big)'\bf{Q}_{j\cdot}.
\end{split}
\end{equation}
We then have, from (\ref{eqn: Q_simplify}) and (\ref{eqn: Q_simplify2}) that
\begin{equation}
\label{eqn: Q_j_consistency}
\begin{split}
    \wh{\bf{Q}}_{j\cdot} - \bf{H}_j \bf{Q}_{j\cdot} &=
    \wh{\bf{D}}^{-1}
    \Bigg(
    \sum_{i=1}^{d_k} \wh{\bf{Q}}_{i\cdot}
    \sum_{h=1}^{d_{\text{-}k}}
        \frac{1}{|\psi_{ij,h}|}
        \sum_{t\in\psi_{ij,h}}
        E_{t,ih} E_{t,jh} \\
        &+
        \sum_{i=1}^{d_k} \wh{\bf{Q}}_{i\cdot}
    \sum_{h=1}^{d_{\text{-}k}}
        \frac{1}{|\psi_{ij,h}|}
        \sum_{t\in\psi_{ij,h}}
        E_{t,jh} \Big(
    \sum_{m=1}^{r_{\text{-}k}}\Lambda_{hm} \bf{F}_{Z,t,\cdot m}
    \Big)'\bf{Q}_{i\cdot} \\
        &+
        \sum_{i=1}^{d_k} \wh{\bf{Q}}_{i\cdot}
    \sum_{h=1}^{d_{\text{-}k}}
        \frac{1}{|\psi_{ij,h}|}
        \sum_{t\in\psi_{ij,h}}
        E_{t,ih}\Big(
    \sum_{m=1}^{r_{\text{-}k}}\Lambda_{hm} \bf{F}_{Z,t,\cdot m}
    \Big)'\bf{Q}_{j\cdot}
        \Bigg)\\
        &=:
        \wh{\bf{D}}^{-1}
        \Big(\c{I}_j+\c{II}_j +\c{III}_j
        \Big),
\end{split}
\end{equation}
where
\begin{equation*}
\begin{split}
    \c{I}_j :&=
    \sum_{i=1}^{d_k} \wh{\bf{Q}}_{i\cdot}
    \sum_{h=1}^{d_{\text{-}k}}
        \frac{1}{|\psi_{ij,h}|}
        \sum_{t\in\psi_{ij,h}}
        E_{t,ih} E_{t,jh},
        \\
    \c{II}_j :&=
    \sum_{i=1}^{d_k} \wh{\bf{Q}}_{i\cdot}
    \sum_{h=1}^{d_{\text{-}k}}
        \frac{1}{|\psi_{ij,h}|}
        \sum_{t\in\psi_{ij,h}}
        E_{t,jh} \Big(
    \sum_{m=1}^{r_{\text{-}k}}\Lambda_{hm} \bf{F}_{Z,t,\cdot m}
    \Big)'\bf{Q}_{i\cdot},
        \\
    \c{III}_j :&=
    \sum_{i=1}^{d_k} \wh{\bf{Q}}_{i\cdot}
    \sum_{h=1}^{d_{\text{-}k}}
        \frac{1}{|\psi_{ij,h}|}
        \sum_{t\in\psi_{ij,h}}
        E_{t,ih}\Big(
    \sum_{m=1}^{r_{\text{-}k}}\Lambda_{hm} \bf{F}_{Z,t,\cdot m}
    \Big)'\bf{Q}_{j\cdot}.
\end{split}
\end{equation*}
The following lemma bounds the terms $\c{I}_j,\c{II}_j$, and $\c{III}_j$.
\begin{lemma}\label{lemma:1}
Under Assumptions (O1), (F1), (L1), (E1) and (E2), we have
\begin{align}
    \frac{1}{d_k}\sum_{j=1}^{d_k}
    \|\c{I}_j\|_F^2 &= O_P\Big(
    \frac{d}{T} + d_{\text{-}k}^2
    \Big), \label{eqn:I_j_bound}\\
    \frac{1}{d_k}\sum_{j=1}^{d_k}
    \|\c{II}_j\|_F^2 &= O_P\Big(
    \frac{\dmk d_k^{\alpha_{k,1}}}
    {T} \Big) = \frac{1}{d_k}\sum_{j=1}^{d_k}
    \|\c{III}_j\|_F^2. \label{eqn:II_j_bound}
\end{align}
\end{lemma}

\noindent\textbf{\textit{Proof of Lemma \ref{lemma:1}.}}
To prove (\ref{eqn:I_j_bound}), we decompose
\begin{equation}
\label{eqn: proof_I_decomp}
\begin{split}
    \c{I}_j = &
    \sum_{i=1}^{d_k} \wh{\bf{Q}}_{i\cdot}
    \sum_{h=1}^{d_{\text{-}k}}
        \frac{\sum_{t\in\psi_{ij,h}}
        E_{t,ih} E_{t,jh}}{|\psi_{ij,h}|} \\
    =&
    \sum_{i=1}^{d_k} \wh{\bf{Q}}_{i\cdot}
    \sum_{h=1}^{d_{\text{-}k}}
    \Big(
        \frac{\sum_{t\in\psi_{ij,h}}
        \big(E_{t,ih}E_{t,jh} -\b{E}[E_{t,ih}E_{t,jh}]\big)}{|\psi_{ij,h}|}
        +
        \frac{\sum_{t\in\psi_{ij,h}}
        \b{E}[E_{t,ih}E_{t,jh}]}{|\psi_{ij,h}|}
    \Big)\\
    =& :
    \sum_{i=1}^{d_k} \wh{\bf{Q}}_{i\cdot}\xi_{ij}
    +\sum_{i=1}^{d_k} \wh{\bf{Q}}_{i\cdot}
    \eta_{ij},
\end{split}
\end{equation}
where
\begin{equation*}
\begin{split}
    \xi_{ij} :=&
    \sum_{h=1}^{d_{\text{-}k}}
        \frac{\sum_{t\in\psi_{ij,h}}
        \big(E_{t,ih}E_{t,jh} -\b{E}[E_{t,ih}E_{t,jh}]\big)}{|\psi_{ij,h}|}, \;\;\;
    \eta_{ij} :=
    \sum_{h=1}^{d_{\text{-}k}}
        \frac{\sum_{t\in\psi_{ij,h}}
        \b{E}[E_{t,ih}E_{t,jh}]}{|\psi_{ij,h}|}.
\end{split}
\end{equation*}
We want to show the following:
\begin{align}
    \sum_{j=1}^{d_k}\Big\|
    \sum_{i=1}^{d_k} \wh{\bf{Q}}_{i\cdot}\xi_{ij}\Big\|_F^2 &=
    O_P\Big(\frac{dd_k}{T}\Big), \label{eqn:R1}\\
    \sum_{j=1}^{d_k}\Big\|
    \sum_{i=1}^{d_k} \wh{\bf{Q}}_{i\cdot}\eta_{ij}\Big\|_F^2 &=
    O_P\big( d d_{\text{-}k} \big). \label{eqn:R2}
\end{align}
To show (\ref{eqn:R1}),
first note that $\b{E}\xi_{ij}=0$, and also by Assumption (O1),
\begin{align}
    \b{E}|\xi_{ij}|^2
    &= \text{var}(\xi_{ij}) \leq
    \frac{1}{\psi_0^2 T^2}
    \text{var}\Bigg(\sum_{h=1}^{d_{\text{-}k}}
    \sum_{t\in\psi_{ij,h}}
\big(E_{t,ih}E_{t,jh} -\b{E}\big[E_{t,ih}E_{t,jh}\big]\big)\Bigg) \notag\\
    & \leq
    \frac{1}{\psi_0^2 T^2}
    \sum_{h=1}^{d_{\text{-}k}}
    \sum_{l=1}^{d_{\text{-}k}}
    \sum_{t\in\psi_{ij,h}}
    \sum_{s\in\psi_{ij,l}}
    \bigg| \text{cov}\Big(
\big(E_{t,ih}E_{t,jh} -\b{E}\big[E_{t,ih}E_{t,jh}\big]\big),
\big(E_{s,il}E_{s,jl} -\b{E}\big[E_{s,il}E_{s,jl}\big]\big)\Big) \bigg| \notag\\
    &=
    \frac{1}{\psi_0^2 T^2}
    \sum_{h=1}^{d_{\text{-}k}}
    \sum_{t\in\psi_{ij,h}}
    \sum_{l=1}^{d_{\text{-}k}}
    \sum_{s\in\psi_{ij,l}}
    \bigg| \text{cov}\Big(
E_{t,ih}E_{t,jh},
E_{s,il}E_{s,jl}\Big)
\bigg|
     \leq
    \frac{cd_{\text{-}k}}{\psi_0^2T}, \label{eqn:xi_ij_sq}
\end{align}
where the last inequality and the constant $c$ are from result 1 of Proposition \ref{Prop:assumption_implications} (hereafter Proposition \ref{Prop:assumption_implications}.1, etc).
Then by the Cauchy--Schwarz inequality,
\begin{equation*}
\begin{split}
    \sum_{j=1}^{d_k}\Big\|
    \sum_{i=1}^{d_k} \wh{\bf{Q}}_{i\cdot}\xi_{ij}\Big\|_F^2
    \leq &
    \sum_{j=1}^{d_k}
    \bigg(
    \sum_{i=1}^{d_k}
    \big\|\wh{\bf{Q}}_{i\cdot}\big\|_F^2
    \bigg)\bigg(
    \sum_{i=1}^{d_k}
    \xi_{ij}^2 \bigg) =
    O_P\Big( \frac{dd_k}{T}\Big),
\end{split}
\end{equation*}
which is (\ref{eqn:R1}). To show (\ref{eqn:R2}),
note that if we define
\begin{equation*}
\begin{split}
    \rho_{ij,h}:&=
        \frac{
        \frac{1}{|\psi_{ij,h}|}\sum_{t\in\psi_{ij,h}}
        \b{E}[E_{t,ih}E_{t,jh}]
        }{
        \big(
        \frac{1}{|\psi_{ij,h}|}\sum_{t\in\psi_{ij,h}}
        \b{E}[E_{t,ih}^2]
        \big)^{1/2}\big(
        \frac{1}{|\psi_{ij,h}|}\sum_{t\in\psi_{ij,h}}
        \b{E}[E_{t,jh}^2]\big)^{1/2}
        },
\end{split}
\end{equation*}
then $|\rho_{ij,h}|<1$ and hence $\rho_{ij,h}^2\leq |\rho_{ij,h}|$. It is then easy to prove also that
\begin{equation*}
\begin{split}
    \rho_{ij}:&=
        \frac{
        \sum_{h=1}^{d_{\text{-}k}}
        \frac{1}{|\psi_{ij,h}|}\sum_{t\in\psi_{ij,h}}
        \b{E}[E_{t,ih}E_{t,jh}]
        }{
        \big(
        \sum_{h=1}^{d_{\text{-}k}}
        \frac{1}{|\psi_{ij,h}|}\sum_{t\in\psi_{ij,h}}
        \b{E}[E_{t,ih}^2]
        \big)^{1/2}\big(
        \sum_{h=1}^{d_{\text{-}k}}
        \frac{1}{|\psi_{ij,h}|}\sum_{t\in\psi_{ij,h}}
        \b{E}[E_{t,jh}^2]\big)^{1/2}
        },
\end{split}
\end{equation*}
also satisfy $|\rho_{ij}|\leq 1$ and $\rho_{ij}^2\leq |\rho_{ij}|$.  By Proposition \ref{Prop:assumption_implications}.1,
\begin{equation*}
\bigg|\sum_{h=1}^{d_{\text{-}k}}
    \frac{\sum_{t\in\psi_{ij,h}}
        \b{E}[E_{t,ih}^2]}
        {|\psi_{ij,h}|} \bigg| = O_P(d_{\text{-}k}),
\end{equation*}
and hence
\begin{equation*}
\begin{split}
    \eta_{ij}^2 &= \bigg(
    \sum_{h=1}^{d_{\text{-}k}}
        \frac{\sum_{t\in\psi_{ij,h}}
        \b{E}[E_{t,ih}E_{t,jh}]}{|\psi_{ij,h}|}
        \bigg)^2
    = \rho_{ij}^2
    \bigg(
    \sum_{h=1}^{d_{\text{-}k}}
    \frac{\sum_{t\in\psi_{ij,h}}
        \b{E}[E_{t,ih}^2]}
        {|\psi_{ij,h}|}\bigg)
    \bigg(
    \sum_{h=1}^{d_{\text{-}k}}
    \frac{\sum_{t\in\psi_{ij,h}}
        \b{E}[E_{t,jh}^2]}
        {|\psi_{ij,h}|}\bigg) \\
    &=
    |\rho_{ij}|
    \bigg(\sum_{h=1}^{d_{\text{-}k}}
    \frac{\sum_{t\in\psi_{ij,h}}
        \b{E}[E_{t,ih}^2]}
        {|\psi_{ij,h}|}\bigg)^{1/2}
    \bigg(\sum_{h=1}^{d_{\text{-}k}}
    \frac{\sum_{t\in\psi_{ij,h}}
        \b{E}[E_{t,jh}^2]}
        {|\psi_{ij,h}|}\bigg) ^{1/2}
        \cdot O_P(d_{\text{-}k}) \\
    &=
    \bigg|\sum_{h=1}^{d_{\text{-}k}}
        \frac{\sum_{t\in\psi_{ij,h}}
        \b{E}[E_{t,ih}E_{t,jh}]}{|\psi_{ij,h}|}
    \bigg| \cdot O_P(d_{\text{-}k})
    =
    \big|\eta_{ij}
    \big| O_P(d_{\text{-}k}).
\end{split}
\end{equation*}
Using the above, we then have
\begin{align}
    \sum_{j=1}^{d_k}\Big\|
    \sum_{i=1}^{d_k} \wh{\bf{Q}}_{i\cdot}\eta_{ij}\Big\|_F^2
    &\leq
    \sum_{j=1}^{d_k}
    \bigg(\sum_{i=1}^{d_k}
    \big\|\wh{\bf{Q}}_{i\cdot}\big\|_F^2
    \bigg)\bigg(
    \sum_{i=1}^{d_k}
    \eta_{ij}^2 \bigg)
    \leq
    r_k\sum_{j=1}^{d_k}
    \sum_{i=1}^{d_k}
    \eta_{ij}^2
    =
    O_P(d_{\text{-}k}) \cdot
    \sum_{j=1}^{d_k}
    \sum_{i=1}^{d_k}\big|\eta_{ij}\big| \notag\\
    &=
    O_P(d_{\text{-}k}) \cdot
    \frac{1}{|\psi_{0}T|}
    \sum_{t\in\psi_{ij,h}}
    \sum_{h=1}^{d_{\text{-}k}}
    \sum_{j=1}^{d_k}
    \sum_{i=1}^{d_k}
    \Big|\b{E}[E_{t,ih}E_{t,jh}]\Big|
     =
    O_P(d d_{\text{-}k}), \label{eqn:eta_ij_sq}
\end{align}
where the second last equality used Assumption (O1), and the last equality used Proposition \ref{Prop:assumption_implications}.1.
This proves (\ref{eqn:R2}).
Using (\ref{eqn:R1}) and (\ref{eqn:R2}), from (\ref{eqn: proof_I_decomp}) we have
\begin{equation*}
\begin{split}
    \frac{1}{d_k}\sum_{j=1}^{d_k}\big\|\c{I}_j\big\|_F^2 &=
    \frac{1}{d_k}\sum_{j=1}^{d_k}
    \Big\|
    \sum_{i=1}^{d_k} \wh{\bf{Q}}_{i\cdot}\sum_{h=1}^{d_{\text{-}k}}\xi_{ij,h}
    +\sum_{i=1}^{d_k} \wh{\bf{Q}}_{i\cdot}
    \sum_{h=1}^{d_{\text{-}k}}\eta_{ij,h}
    \Big\|_F^2 \\
    &\leq
    \frac{2}{d_k}\sum_{j=1}^{d_k}
    \Big\|
    \sum_{i=1}^{d_k} \wh{\bf{Q}}_{i\cdot}\xi_{ij}\Big\|_F^2
    +
    \frac{2}{d_k}\sum_{j=1}^{d_k}
    \Big\|
    \sum_{i=1}^{d_k} \wh{\bf{Q}}_{i\cdot}
    \eta_{ij}
    \Big\|_F^2
    =
    O_P\Big(
    \frac{d}{T} + d_{\text{-}k}^2
    \Big).
\end{split}
\end{equation*}
This completes the proof of (\ref{eqn:I_j_bound}). To prove (\ref{eqn:II_j_bound}), consider
\begin{equation}
\label{eqn: proof_lemma1.2}
\begin{split}
    \frac{1}{d_k}\sum_{j=1}^{d_k}
    \|\c{II}_j\|_F^2 &=
    \frac{1}{d_k}\sum_{j=1}^{d_k}
    \bigg\|
    \sum_{i=1}^{d_k} \wh{\bf{Q}}_{i\cdot}
    \sum_{h=1}^{d_{\text{-}k}}
        \frac{1}{|\psi_{ij,h}|}
        \sum_{t\in\psi_{ij,h}}
        E_{t,jh} \Big(
    \sum_{m=1}^{r_{\text{-}k}}\Lambda_{hm} \bf{F}_{Z,t,\cdot m}
    \Big)'\bf{Q}_{i\cdot}
    \bigg\|_F^2 \\
    &\leq
    \frac{1}{d_k}\sum_{j=1}^{d_k}
    \bigg(
    \sum_{i=1}^{d_k}
    \|\wh{\bf{Q}}_{i\cdot}\|_F^2 \bigg)
    \cdot \sum_{i=1}^{d_k} \bigg(
    \sum_{h=1}^{d_{\text{-}k}}
        \frac{1}{|\psi_{ij,h}|}
        \sum_{t\in\psi_{ij,h}}
        E_{t,jh} \Big(
    \sum_{m=1}^{r_{\text{-}k}}\Lambda_{hm} \bf{F}_{Z,t,\cdot m}
    \Big)'\bf{Q}_{i\cdot}
    \bigg)^2  \\
    &=
    \frac{r_k}{d_k}\sum_{j=1}^{d_k}
    \sum_{i=1}^{d_k}\Bigg(
    \sum_{h=1}^{d_{\text{-}k}}
        \frac{1}{|\psi_{ij,h}|}
        \sum_{t\in\psi_{ij,h}}
        E_{t,jh} \Big(
    \sum_{m=1}^{r_{\text{-}k}}\Lambda_{hm} \bf{F}_{Z,t,\cdot m}
    \Big)'\bf{Q}_{i\cdot}
    \Bigg)^2 \\
    &=
    \frac{r_k\dmk^2}{d_k}\sum_{j=1}^{d_k}
    \sum_{i=1}^{d_k}\Bigg(
        \sum_{h=1}^{\dmk}
        \sum_{t\in\psi_{ij,h}}
        \frac{1}{\dmk \cdot |\psi_{ij,h}|}
        E_{t,jh} [\otimes_{l\in[K]\setminus \{k\}}\bf{A}_l]_{h\cdot}'
        \bf{F}_{t}' \bf{A}_{i\cdot}
    \Bigg)^2\\
    &=
    \frac{r_k\dmk^2}{d_k}\sum_{j=1}^{d_k}
    \sum_{i=1}^{d_k}\|\bf{u}_i\|_F^2\Bigg(
        \sum_{h=1}^{\dmk}
        \sum_{t\in\psi_{ij,h}}
        \frac{1}{\dmk \cdot |\psi_{ij,h}|}
        E_{t,jh}\bf{v}_h' \bf{F}_{t}' \frac{1}
        {\|\bf{u}_i\|_F}\bf{u}_i
    \Bigg)^2 ,
\end{split}
\end{equation}
where $\bf{A} = \bf{A}_k$ and $\bf{F}_t\ = \mat{k}{\c{F}_t}$ above, and we define
$
    \bf{v}_h :=
    [\otimes_{l\in[K]\setminus \{k\}}\bf{A}_l]_{h\cdot}, \;
    \bf{u}_i :=
    \bf{A}_{i\cdot} .
$
By Proposition \ref{Prop:assumption_implications}.2, the last bracket in the last line of (\ref{eqn: proof_lemma1.2}) is $O_P(\dmk^{-1}T^{-1})$, and hence
\begin{equation*}
\begin{split}
    \frac{1}{d_k}\sum_{j=1}^{d_k}
    \|\c{II}_j\|_F^2 &=
    O_P(\frac{\dmk}{d_k T})
    \sum_{j=1}^{d_k}
    \sum_{i=1}^{d_k}\|\bf{u}_i\|_F^2
    =
    O_P\Big(
    \frac{d_{\text{-}k}}
    {d_kT} \Big)\sum_{j=1}^{d_k}
    \big\|\bf{A}\big\|_F^2
    =
    O_P\Big(
    \frac{\dmk d_k^{\alpha_{k,1}}}
    {T} \Big),
\end{split}
\end{equation*}
where the last equality follows since for any $l\in[K], \|\bf{A}_l\|_F^2=O_P\big(\tr(\bf{Z}_l)\big)=O_P\big(d_l^{\alpha_{l,1}}\big)$ by Assumption (L1). The bound corresponding to $\c{III}_j$ can be proved similarly (omitted), and hence (\ref{eqn:II_j_bound}) is established. This concludes the proof of Lemma \ref{lemma:1}.
$\square$

\begin{lemma}\label{lemma:2}
Under Assumptions (O1), (M1), (F1), (L1), (E1), (E2) and (R1), with $\H_j$ and $\wh\D$ from (\ref{eqn: notation_short}), we have
\begin{align}
    \big\|\wh\D^{-1}\big\|_F &=
     O_P\Bigg( d_k^{\alpha_{k,1}-\alpha_{k,r_k}}\prod_{j=1}^Kd_j^{-\alpha_{j,1}} \Bigg),
     \label{eqn: lemma2.1}\\
    \frac{1}{d_k}
    \sum_{j=1}^{d_k}\Big\|
    \wh\Q_{j\cdot}-\H_j
    \Q_{j\cdot}
    \Big\|_F^2 &=
    O_P\Bigg(d_k^{2(\alpha_{k,1} - \alpha_{k,r_k})-1}
    \bigg(\frac{1}{T\dmk}+\frac{1}{d_k}\bigg)
    \prod_{j=1}^Kd_j^{2(1-\alpha_{j,1})}
    \Bigg)
    . \label{eqn:_lemm2.2}
\end{align}
\end{lemma}

\noindent\textbf{\textit{Proof of Lemma \ref{lemma:2}.}}
First, we bound the term $\big\|\wh\D^{-1}\big\|_F^2$ by finding the lower bound of $\lambda_{r_k}(\wh\D)$. To do this, define $\omega_k := d_k^{\alpha_{k,r_k}-\alpha_{k,1}}\prod_{j=1}^Kd_j^{\alpha_{j,1}}$, and consider the decomposition
\begin{align}
  \wh\S &= \R^* + (\wt\R - \R^*) + \R_1 + \R_2 + \R_3, \label{eqn:S_decomposition}
\end{align}
where for a unit vector $\bgamma$,
\begin{align}
    R(\bm{\gamma}) &:=
    \frac{1}{\omega_k}\bm{\gamma}' \wh\S \bm{\gamma} =
    \frac{1}{\omega_k}
    \sum_{i=1}^{d_k}\sum_{j=1}^{d_k}
    \bm{\gamma}_i\bm{\gamma}_j
    \wh{S}_{ij} =
    \frac{1}{\omega_k}
    \sum_{i=1}^{d_k}\sum_{j=1}^{d_k}
    \bm{\gamma}_i\bm{\gamma}_j
    \sum_{h=1}^{d_{\text{-}k}}
        \frac{1}{|\psi_{ij,h}|}
        \sum_{t\in\psi_{ij,h}}
        Y_{t,ih}Y_{t,jh} \notag\\
        &=: R^*(\bgamma) + (\wt{R}(\bgamma)-R^*(\bgamma)) + R_1 + R_2 + R_3, \; \text{ with} \notag\\
    \wt{R}(\bgamma) &:= \frac{1}{\omega_k}\bgamma'\wt\R\bgamma :=
    \frac{1}{\omega_k}
    \sum_{i=1}^{d_k}\sum_{j=1}^{d_k}
    \bm{\gamma}_i\bm{\gamma}_j
    \sum_{h=1}^{d_{\text{-}k}}
        \frac{1}{|\psi_{ij,h}|}
        \sum_{t\in\psi_{ij,h}}
        \Big(
    \sum_{m=1}^{r_{\text{-}k}}\Lambda_{hm} \F_{Z,t,\cdot m}
    \Big)'\Q_{i\cdot}\Big(
    \sum_{m=1}^{r_{\text{-}k}}\Lambda_{hm} \F_{Z,t,\cdot m}
    \Big)'\Q_{j\cdot}, \notag\\
    R^*(\bgamma) &:= \frac{1}{\omega_k}\bgamma'\R^*\bgamma
    :=    \frac{1}{\omega_k}
    \sum_{i=1}^{d_k}\sum_{j=1}^{d_k}
    \bm{\gamma}_i\bm{\gamma}_j
    \sum_{h=1}^{d_{\text{-}k}}
        \frac{1}{T}
        \sum_{t=1}^T
        \Big(
    \sum_{m=1}^{r_{\text{-}k}}\Lambda_{hm} \F_{Z,t,\cdot m}
    \Big)'\Q_{i\cdot}
    \Big(
    \sum_{m=1}^{r_{\text{-}k}}\Lambda_{hm} \F_{Z,t,\cdot m}
    \Big)'\Q_{j\cdot},\notag\\
    R_1 &:= \frac{1}{\omega_k}\bgamma'\R_1\bgamma :=
    \frac{1}{\omega_k}
    \sum_{i=1}^{d_k}\sum_{j=1}^{d_k}
    \bm{\gamma}_i\bm{\gamma}_j
    \sum_{h=1}^{d_{\text{-}k}}
        \frac{1}{|\psi_{ij,h}|}
        \sum_{t\in\psi_{ij,h}}
    E_{t,ih}E_{t,jh},
    \notag\\
    R_2 &:= \frac{1}{\omega_k}\bgamma'\R_2\bgamma :=
    \frac{1}{\omega_k}
    \sum_{i=1}^{d_k}\sum_{j=1}^{d_k}
    \bm{\gamma}_i\bm{\gamma}_j
    \sum_{h=1}^{d_{\text{-}k}}
        \frac{1}{|\psi_{ij,h}|}
        \sum_{t\in\psi_{ij,h}}
    E_{t,jh}\Big(
    \sum_{m=1}^{r_{\text{-}k}}\Lambda_{hm} \F_{Z,t,\cdot m}
    \Big)'\Q_{i\cdot}, \notag\\
    R_3 &:= \frac{1}{\omega_k}\bgamma'\R_3\bgamma :=
    \frac{1}{\omega_k}
    \sum_{i=1}^{d_k}\sum_{j=1}^{d_k}
    \bm{\gamma}_i\bm{\gamma}_j
    \sum_{h=1}^{d_{\text{-}k}}
        \frac{1}{|\psi_{ij,h}|}
        \sum_{t\in\psi_{ij,h}}
    E_{t,ih}\Big(
    \sum_{m=1}^{r_{\text{-}k}}\Lambda_{hm} \F_{Z,t,\cdot m}
    \Big)'\Q_{j\cdot}, \label{eqn:different_R}
\end{align}
and we used (\ref{eqn: Q_simplify2}) for the expansion above.
Then we have the decomposition
\begin{equation}
\label{eqn: proof_RR_decom}
\begin{split}
    R(\bm{\gamma})-R^{\ast}(\bm{\gamma})
    =R(\bm{\gamma}) -\tilde{R}(\bm{\gamma})
    +\tilde{R}(\bm{\gamma}) - R^{\ast}(\bm{\gamma}) .
\end{split}
\end{equation}
Similar to the treatment of the term $\c{I}_j$ in the proof of Lemma \ref{lemma:1}, since $\|\bgamma\|=1$,
\begin{align}
    |R_1|
    &\leq
    \frac{1}{\omega_k}
    \bigg|\sum_{i=1}^{d_k}\sum_{j=1}^{d_k}
    \bm{\gamma}_i\bm{\gamma}_j
    \xi_{ij} \bigg| +
    \frac{1}{\omega_k}\bigg|
    \sum_{i=1}^{d_k}\sum_{j=1}^{d_k}
    \bm{\gamma}_i\bm{\gamma}_j
    \eta_{ij} \bigg|
    \leq
    \frac{1}{\omega_k}
    \bigg(\sum_{i=1}^{d_k}\sum_{j=1}^{d_k}
    \xi_{ij}^2 \bigg)^{1/2} +
    \frac{1}{\omega_k}\bigg(
    \sum_{i=1}^{d_k}\sum_{j=1}^{d_k}
    \eta_{ij}^2 \bigg)^{1/2} \notag\\
    &=
    O_P\bigg(d_k^{\alpha_{k,1} - \alpha_{k,r_k}}\Big[\frac{1}{T^{1/2}d_{\text{-}k}^{1/2}}
    + \frac{1}{d_k^{1/2}}\Big]\prod_{j=1}^Kd_j^{1-\alpha_{j,1}}\bigg) = O_P(d[(T\dmk)^{-1/2} + d_k^{-1/2}]/\omega_k), \label{eqn:lem2_withoutE1}
\end{align}
where the second last equality is from (\ref{eqn:xi_ij_sq}) and part of (\ref{eqn:eta_ij_sq}).
Together with Assumption (R1), (\ref{eqn:lem2_withoutE1}) implies that as $T,d_k,d_{\text{-}k}\to\infty$, $R_1 \xrightarrow{p} 0$.

From (\ref{eqn: proof_lemma1.2}) and the arguments for $\c{II}_j$ immediately afterwards, we see that
\begin{align}
  |R_2|
    &\leq \frac{1}{\omega_k}\Bigg\{\sum_{j=1}^{d_k}
    \sum_{i=1}^{d_k}\Bigg(
    \sum_{h=1}^{d_{\text{-}k}}
        \frac{1}{|\psi_{ij,h}|}
        \sum_{t\in\psi_{ij,h}}
        E_{t,jh} \Big(
    \sum_{m=1}^{r_{\text{-}k}}\Lambda_{hm} \bf{F}_{Z,t,\cdot m}
    \Big)'\bf{Q}_{i\cdot}
    \Bigg)^2\Bigg\}^{1/2}\notag\\
    &= O_P(\omega_k^{-1})\cdot O_P\bigg( T^{-1/2}d^{1/2}\prod_{j=1}^Kd_j^{\alpha_{j,1}/2} \bigg)
    = O_P((dg_s)^{1/2}T^{-1/2}/\omega_k) = o_P(1), \label{eqn:R2rate}
\end{align}
where the last equality is from Assumption (R1). The term $R_3$ can be proved to have the same rate with same lines of proof as for $R_2$. Hence we have
\begin{equation}\label{eqn:RminusRtildeto0}
\sup_{\|\bgamma\|=1} |R(\bgamma) - \wt{R}(\bgamma)| = \sup_{\|\bgamma\|=1}|R_1+R_2+R_3| \xrightarrow{p} 0.
\end{equation}

Similar to the proof of (R6) in Lemma 4 in \cite{Xiong_Pelger}, using the definition $\v_h := [\otimes_{l\in[K]\setminus \{k\}}\bf{A}_l]_{h\cdot}$ as before,
\begin{align}
    &\hspace{12pt}
    \tilde{R}(\bm{\gamma})
    -R^{\ast}(\bm{\gamma}) \notag \\
    &=
    \frac{1}{\omega_k}
    \sum_{i=1}^{d_k}\sum_{j=1}^{d_k}
    \bm{\gamma}_i\bm{\gamma}_j
    \sum_{h=1}^{d_{\text{-}k}}
    \Big(
        \frac{1}{|\psi_{ij,h}|}
        \sum_{t\in\psi_{ij,h}}
        \Q_{i\cdot}'
        \F_{Z,t}
        \bf{\Lambda}_{h\cdot}
    \bf{\Lambda}_{h\cdot}'
    \F_{Z,t}'\Q_{j\cdot}
    - \frac{1}{T}
        \sum_{t=1}^T
        \Q_{i\cdot}'
        \F_{Z,t}
        \bf{\Lambda}_{h\cdot}
    \bf{\Lambda}_{h\cdot}'
    \F_{Z,t}'\Q_{j\cdot}
     \Big) \notag\\
    &=
    \frac{1}{\omega_k}
    \sum_{i=1}^{d_k}\sum_{j=1}^{d_k}
    \bm{\gamma}_i\bm{\gamma}_j
    \sum_{h=1}^{d_{\text{-}k}}
    \A_{i\cdot}'
    \Big(
        \frac{1}{|\psi_{ij,h}|}
        \sum_{t\in \psi_{ij,h}}
        \F_{t}\v_h\v_h'\F_{t}'
     -
    \frac{1}{T}\sum_{t=1}^T\F_{t}\v_h\v_h'\F_{t}' \Big)\A_{j\cdot} \notag\\
     &=:
     \frac{1}{\omega_k}
    \sum_{i=1}^{d_k}\sum_{j=1}^{d_k}
    \bm{\gamma}_i\bm{\gamma}_j
    \sum_{h=1}^{d_{\text{-}k}}
    \A_{i\cdot}'\Delta_{F,k,ij,h}
    \A_{j\cdot}, \;\text{ where } \label{eqn: proof_delta_def}\\
&\Delta_{F,k,ij,h} :=
\frac{1}{|\psi_{ij,h}|}
        \sum_{t\in\psi_{ij,h}}
        \F_{t}\v_h\v_h'\F_{t}'
- \frac{1}{T}
        \sum_{t=1}^T
        \F_{t}\v_h\v_h'\F_{t}'. \notag
\end{align}

By the Cauchy--Schwarz inequality, we then have
\begin{equation}
\label{eqn: proof_R_tilde_R_ast}
\begin{split}
    \big|\tilde{R}(\bm{\gamma})
    -R^{\ast}(\bm{\gamma})\big|
    &\leq
    \frac{1}{\omega_k} \Bigg(
    \sum_{i=1}^{d_k}\sum_{j=1}^{d_k}
    \Big[\A_{i\cdot}'
    \Big(\sum_{h=1}^{d_{\text{-}k}}
    \Delta_{F,k,ij,h}\Big)
    \A_{j\cdot}\Big]^2
    \Bigg)^{1/2}.
\end{split}
\end{equation}
With Assumption (M1), using the standard rate of convergence in the weak law of large number for $\alpha$-mixing sequence and the fact that the elements in $\c{F}_t$ are independent from Assumption (F1), since $\Delta_{F,k,ij,h}$ has fixed dimension, we have for each $k\in[K], i,j\in[d_k]$ and $h\in[\dmk]$,
\begin{align}
  \|\Delta_{F,k,ij,h}\|_F \leq  \bigg\|\frac{1}{|\psi_{ij,h}|}
        \sum_{t\in\psi_{ij,h}}
        \F_{t}\v_h\v_h'\F_{t}' - \v_h'\v_h\bSigma_k\bigg\|_F
 + \bigg\|\frac{1}{T}
        \sum_{t=1}^T
        \F_{t}\v_h\v_h'\F_{t}' - \v_h'\v_h\bSigma_k\bigg\|_F = O_P\Big(\frac{\|\v_h\|^2}{T^{1/2}}\Big).
        \label{eqn:Delta_{F,k,ij,h}_rate1}
\end{align}
From Assumption (L1), we then have
\begin{align}
    \sum_{i=1}^{d_k}\sum_{j=1}^{d_k}
    \Big[\A_{i\cdot}'
    \Big(\sum_{h=1}^{d_{\text{-}k}}
    \Delta_{F,k,ij,h}\Big)
    \A_{j\cdot}\Big]^2
    &\leq
    \sum_{i=1}^{d_k}\sum_{j=1}^{d_k}
    \|\A_{i\cdot}\|_F^2
    \|\A_{j\cdot}\|_F^2
    \Big(\sum_{h=1}^{\dmk}
    \|\Delta_{F,k,ij,h}\|_F\Big)^2 \notag\\
    &\leq
    \|\A_k\|_F^4\cdot O_P\bigg(\frac{1}{T}\Big(\sum_{h=1}^{\dmk}\|\v_h\|^2\Big)^2\bigg)
    = \|\A_k\|_F^4\cdot O_P\bigg( \frac{1}{T} \|\otimes_{j\in[K]\setminus \{k\}}\A_j\|_F^4 \bigg) \notag\\
    &= O_P\bigg(\frac{1}{T}\prod_{j=1}^K \norm{\A_j}_F^4 \bigg)
     = O_P\bigg(\frac{1}{T}\prod_{j=1}^K d_j^{2\alpha_{j,1}}\bigg), \label{eqn:DeltaFrobenius}
\end{align}
and hence from (\ref{eqn: proof_R_tilde_R_ast}), we have by Assumption (R1) that
\begin{align}
    \big|\tilde{R}(\bm{\gamma})
    -R^{\ast}(\bm{\gamma}) \big|
    =
    O_P\Big(\frac{1}{\sqrt{T}}d_k^{\alpha_{k,1}-\alpha_{k,r_k}}\Big) = o_P(d^{-1/2}g_s^{1/2}) = o_P(1), \label{eqn:RtildeminusR*to0}
\end{align}
where the second last equality is from Assumption (R1).
Next, with Proposition \ref{Prop:assumption_implications}.3, consider
\begin{align*}
  \lambda_{r_k}(\R^*) &= \lambda_{r_k}\Bigg(
    \frac{1}{T}\sum_{t=1}^T
    \Q \F_{Z,t} \bf{\Lambda}'
    \bf{\Lambda} \F_{Z,t}' \Q'
    \Bigg) = \lambda_{r_k}\Bigg(
    \frac{1}{T}\sum_{t=1}^T
    \A_k\F_{t}
    \big[\otimes_{j\in[K]\setminus\{k\}}
    \A_j\big]'
    \big[\otimes_{j\in[K]\setminus\{k\}}
    \A_j\big]
    \F_{t}'\A_k'
    \Bigg)\\
    &\geq \lambda_{r_k}(\A_k'\A_k)\cdot \lambda_{r_k}\Bigg(\frac{1}{T}\sum_{t=1}^T
    \F_{t}
    \big[\otimes_{j\in[K]\setminus\{k\}}
    \A_j\big]'
    \big[\otimes_{j\in[K]\setminus\{k\}}
    \A_j\big]
    \F_{t}'\Bigg)\\
    &\asymp_P d_k^{\alpha_{k,r_k}} \cdot \lambda_{r_k}(\tr(\otimes_{j\in[K]\setminus\{k\}}\A_j'\A_j)\bSigma_k) \asymp_P d_k^{\alpha_{k,r_k}}\prod_{j\in[K]\setminus\{k\}}d_j^{\alpha_{j,1}} = \omega_k.
\end{align*}

With this, going back to the decomposition (\ref{eqn:S_decomposition}),
\begin{align*}
  \omega_k^{-1}\lambda_{r_k}(\wh\D) &= \omega_k^{-1}\lambda_{r_k}(\wh\S) \geq \omega_k^{-1}\lambda_{r_k}(\R^*) - \sup_{\|\bgamma\|=1}|\wt{R}(\bgamma) - R^*(\bgamma)| - \sup_{\|\bgamma\|=1}|R_1+R_2+R_3|\asymp_P 1,
\end{align*}
where we used (\ref{eqn:RminusRtildeto0}) and (\ref{eqn:RtildeminusR*to0}). Hence finally,
\begin{align*}
  \norm{\wh\D^{-1}}_F = O_P\big( \lambda_{r_k}^{-1}(\wh\D) \big) = O_P(\omega_k^{-1})
  = O_P\bigg(d_k^{\alpha_{k,1}-\alpha_{k,r_k}}\prod_{j=1}^Kd_j^{-\alpha_{j,1}}\bigg),
\end{align*}
which completes the proof of (\ref{eqn: lemma2.1}).

To prove (\ref{eqn:_lemm2.2}),
from (\ref{eqn: Q_j_consistency}) we obtain
\begin{equation*}
\label{eqn: lemma2-decom}
\begin{split}
    \frac{1}{d_k}
    \sum_{j=1}^{d_k}\Big\|
    \wh\Q_{j\cdot}-\H_j
    \Q_{j\cdot}
    \Big\|_F^2  &=
    \frac{1}{d_k}
    \sum_{j=1}^{d_k}\Big\|
    \wh\D^{-1}
    \Big(\c{I}_j+\c{II}_j +\c{III}_j
    \Big) \Big\|_F^2
    \leq
    \big\|\wh\D^{-1}
    \big\|_F^2
    \Bigg[
    \frac{1}{d_k}
    \sum_{j=1}^{d_k}\Big\|
        \c{I}_j+\c{II}_j +\c{III}_j
    \Big\|_F^2 \Bigg] \\
    &\leq
    \big\|\wh\D^{-1}
    \big\|_F^2
    \Bigg[
    \bigg(
    \frac{2}{d_k}
    \sum_{j=1}^{d_k}\big\|
        \c{I}_j
    \big\|_F^2 \bigg) +
    \bigg(
    \frac{4}{d_k}
    \sum_{j=1}^{d_k}\big\|
        \c{II}_j
    \big\|_F^2 \bigg) +
    \bigg(
    \frac{4}{d_k}
    \sum_{j=1}^{d_k}\big\|
        \c{III}_j
    \big\|_F^2 \bigg)
    \Bigg]\\
    &=
    O_P\Bigg(d_k^{2(\alpha_{k,1} - \alpha_{k,r_k})-1}
    \bigg(\frac{1}{T\dmk}+\frac{1}{d_k}
    \bigg)\prod_{j=1}^Kd_j^{2(1-\alpha_{j,1})}
    \Bigg) ,
\end{split}
\end{equation*}
where the last line follows from (\ref{eqn: lemma2.1}) and Lemma \ref{lemma:1}. This concludes the proof of Lemma \ref{lemma:2}.
$\square$

\begin{lemma}\label{lemma:3}
Under the Assumptions in Lemma \ref{lemma:2}, for each $i,j\in[d_k]$ and $h\in[\dmk]$, if the set $\psi_{ij,h}$ has cardinality satisfying
\[\frac{|\psi_{ij,h}|}{T} = 1-\eta_{i,j,h} \geq 1-\eta,\]
then for all $j\in[d_k]$ with $\H_j$ and $\H^a$ from (\ref{eqn: notation_short}),
\begin{equation}
    \big\|\bf{H}_j - \H^a\big\|_F^2 =
    O_P\bigg(
    \min\bigg(\frac{1}{T}, \frac{\eta^2}{(1-\eta)^2}\bigg)
    d_k^{2(\alpha_{k,1}-\alpha_{k,r_k})}
    \bigg) = o_P(1).
\end{equation}
\end{lemma}
\textbf{\textit{Proof of Lemma \ref{lemma:3}.}}
Firstly, consider $\Delta_{F,k,ij,h}$ from (\ref{eqn: proof_delta_def}), where
\begin{align}
\|\Delta_{F,k,ij,h}\|_F &= \bigg\|\frac{1}{|\psi_{ij,h}|}\sum_{t\in\psi_{ij,h}}\F_t\v_h\v_h'\F_t' - \frac{1}{T}\sum_{t=1}^T\F_t\v_h\v_h'\F_t'\bigg\|_F\notag\\
&= \bigg\| \frac{1}{|\psi_{ij,h}|}\sum_{t\in\psi_{ij,h}^c}\F_t\v_h\v_h'\F_t' \bigg\|_F +
\bigg\|\bigg(\frac{1}{|\psi_{ij,h}|} - \frac{1}{T}\bigg)\sum_{t=1}^T\F_t\v_h\v_h'\F_t'\bigg\|_F \notag\\
&\leq O_P\bigg(\frac{T\eta}{T-T\eta}\norm{\v_h}^2\bigg) + \frac{T\eta}{T(T-T\eta)}\cdot O_P(T\norm{\v_h}^2) \notag\\
&= O_P\bigg(\frac{\eta}{1-\eta}\norm{\v_h}^2\bigg). \notag
\end{align}
Combining this with (\ref{eqn:Delta_{F,k,ij,h}_rate1}), we have
\begin{equation}\label{eqn:Delta_final_rate}
\|\Delta_{F,k,ij,h}\|_F = O_P\bigg(\min\bigg(\frac{1}{\sqrt{T}}, \frac{\eta}{1-\eta}\bigg)\norm{\v_h}^2\bigg).
\end{equation}

Note also that
\begin{equation*}
\begin{split}
    \bf{H}_j=&
    \wh{\bf{D}}^{-1}\sum_{i=1}^{d_k} \wh{\bf{Q}}_{i\cdot}
    \sum_{h=1}^{d_{\text{-}k}}
        \frac{1}{|\psi_{ij,h}|}
        \sum_{t\in\psi_{ij,h}}
        \bf{\Lambda}_{h\cdot}'
        \bf{F}_{Z,t}'\bf{Q}_{i\cdot}
        \bf{\Lambda}_{h\cdot}'
        \bf{F}_{Z,t}' \\
    =&
    \wh{\bf{D}}^{-1}\sum_{i=1}^{d_k} \wh{\bf{Q}}_{i\cdot}
    \sum_{h=1}^{d_{\text{-}k}}
        \frac{1}{|\psi_{ij,h}|}
        \sum_{t\in\psi_{ij,h}}
        \bf{Q}_{i\cdot}'
        \bf{F}_{Z,t}
        \bf{\Lambda}_{h\cdot}
        \bf{\Lambda}_{h\cdot}'
        \bf{F}_{Z,t}' ,
\end{split}
\end{equation*}
and also
\begin{equation*}
\begin{split}
    \H^a=&
    \wh{\bf{D}}^{-1}\sum_{i=1}^{d_k} \wh{\bf{Q}}_{i\cdot}
    \sum_{h=1}^{d_{\text{-}k}}
        \frac{1}{T}
        \sum_{t=1}^T
        \bf{Q}_{i\cdot}'
        \bf{F}_{Z,t}
        \bf{\Lambda}_{h\cdot}
        \bf{\Lambda}_{h\cdot}'
        \bf{F}_{Z,t}' .
\end{split}
\end{equation*}
We then have, for $\wh\A_{i\cdot} = \Z_k^{1/2}\wh\Q_{i\cdot}$ and using (\ref{eqn:Delta_final_rate}),
\begin{equation*}
\begin{split}
    \big\|\bf{H}_j-\H^a\big\|_F^2
    &= \bigg\|
    \wh{\bf{D}}^{-1}\sum_{i=1}^{d_k} \wh{\bf{Q}}_{i\cdot}
    \sum_{h=1}^{d_{\text{-}k}}
    \Big(
    \frac{1}{|\psi_{ij,h}|}
        \sum_{t\in\psi_{ij,h}}
        \bf{Q}_{i\cdot}'
        \bf{F}_{Z,t}
        \bf{\Lambda}_{h\cdot}
        \bf{\Lambda}_{h\cdot}'
        \bf{F}_{Z,t}'
        -
        \frac{1}{T}
        \sum_{t=1}^T
        \bf{Q}_{i\cdot}'
        \bf{F}_{Z,t}
        \bf{\Lambda}_{h\cdot}
        \bf{\Lambda}_{h\cdot}'
        \bf{F}_{Z,t}'
    \Big)
    \bigg\|_F^2 \\
     &=
    \bigg\|
    \wh{\bf{D}}^{-1}\sum_{i=1}^{d_k} \wh{\Q}_{i\cdot}\A_{i\cdot}'
    \sum_{h=1}^{\dmk}
    \Delta_{F,k,ij,h}\Z_k^{1/2}
    \bigg\|_F^2  \leq
    \big\|\wh{\bf{D}}^{-1}\big\|_F^2
    \bigg\|
    \sum_{i=1}^{d_k} \wh{\A}_{i\cdot}\A_{i\cdot}'
    \sum_{h=1}^{d_{\text{-}k}}
    \Delta_{F,k,ij,h}
    \bigg\|_F^2 \\
    &= O_P(\omega_k^{-2})\cdot \max_{i\in[d_k]}\bigg\|\sum_{h=1}^{\dmk}\Delta_{F,k,ij,h}\bigg\|_F^2\bigg(\sum_{i=1}^{d_k}
    \|\wh\A_{i\cdot}\|\cdot\|\A_{i\cdot}\|\bigg)^2\\
    &= O_P(\omega_k^{-2})\cdot O_P\bigg( \min\bigg(\frac{1}{T}, \frac{\eta^2}{(1-\eta)^2}\bigg)\|\otimes_{j\in[K]\setminus\{k\}}\A_j\|_F^4 \bigg)
       \cdot \norm{\wh\A_k}_F^2\cdot\|\A_k\|_F^2\\
    &= O_P\bigg( \min\bigg(\frac{1}{T},\frac{\eta^2}{(1-\eta)^2}\bigg)d_k^{2(\alpha_{k,1} - \alpha_{k,r_k})} \bigg) = o_P(1),
\end{split}
\end{equation*}
where we used (\ref{eqn:Delta_final_rate}) in the second last line, Assumption (L1) in the second last equality, and Assumption (R1) in the last equality.
This completes the proof of Lemma \ref{lemma:3}. $\square$


\textbf{\textit{Proof of Theorem \ref{thm:factor_loading_consistency}.}}

The first result is shown in Lemma \ref{lemma:2}. Together with Lemma \ref{lemma:3}, we have
\begin{equation*}
\begin{split}
    \frac{1}{d_k}
    \sum_{j=1}^{d_k}\big\|
    \wh\Q_{j\cdot}-\H^a
    \Q_{j\cdot}
    \big\|^2
    &\leq
    \frac{1}{d_k}
    \sum_{j=1}^{d_k}\big\|
    \wh\Q_{j\cdot}-\H_j \Q_{j\cdot}
    \big\|_F^2 +
    \frac{1}{d_k}
    \sum_{j=1}^{d_k}\Big\|
    \Big(\H_j-\H^a\Big)
    \Q_{j\cdot} \Big\|_F^2 \\
    &\leq
    \frac{1}{d_k}
    \sum_{j=1}^{d_k}\big\|
    \wh\Q_{j\cdot}-\H_j
    \Q_{j\cdot}\big\|_F^2 +
    \frac{1}{d_k}
    \sum_{j=1}^{d_k}\big\|
    \H_j-\H^a\big\|_F^2
    \big\|\Q_{j\cdot}\big\|_F^2 \\
    &=
    O_P\Bigg(d_k^{2(\alpha_{k,1} - \alpha_{k,r_k})-1}
    \bigg(\frac{1}{T\dmk}+\frac{1}{d_k}
    \bigg)\prod_{j=1}^Kd_j^{2(1-\alpha_{j,1})}
    \Bigg)\\
    &\;\;\; +
    O_P\bigg(\min\bigg(\frac{1}{T},\frac{\eta^2}{(1-\eta)^2}\bigg)d_k^{2(\alpha_{k,1} - \alpha_{k,r_k})-1}\bigg) \\
    &=
    O_P\Bigg(d_k^{2(\alpha_{k,1} - \alpha_{k,r_k})-1}
    \Bigg[
    \bigg(\frac{1}{T\dmk}+\frac{1}{d_k}
    \bigg)\frac{d^2}{g_s^2}
    + \min\bigg(\frac{1}{T},\frac{\eta^2}{(1-\eta)^2}\bigg)\Bigg]\Bigg),
\end{split}
\end{equation*}
where we used $d_k^{-1}\sum_{j=1}^{d_k}\|\bf{Q}_{j\cdot}\|_F^2=d_k^{-1}\|\bf{Q}_k\|_F^2=O(d_k^{-1})$ by Assumption (L1). This completes the proof of the theorem. $\square$

Before we prove the consistency results for our imputations, we want to prove asymptotic normality for our factor loading estimators first. Consistency for the imputations will then use the rate obtained from  asymptotic normality of the estimated factor loading matrices. We present a lemma first before proving Theorem \ref{thm:asymp_normality_loadings}.

\begin{lemma}\label{lemma:limit}
 Let Assumption (O1), (M1), (F1), (L1), (L2), (E1), (E2) and (R1) hold. For a given $k\in[K]$, let $\R^\ast$ be from (\ref{eqn:S_decomposition}) and $\omega_k := d_k^{\alpha_{k,r_k} - \alpha_{k,1}}g_s$. Then
\begin{align*}
    \omega_k^{-1}\R^\ast &\xrightarrow{p}
    \textnormal{tr}(\Amk'\Amk) \cdot \omega_k^{-1} \A_k\A_k' ,
    \\
    \omega_k^{-1}\wh\D_k &\xrightarrow{p}
    \omega_k^{-1}\D_k
    := \omega_k^{-1}
    \textnormal{tr}(\Amk'\Amk) \cdot
    \textnormal{diag}
    \{\lambda_j(\A_k'\A_k)
    \mid j\in[r_k]\} ,
    \\
    \H_k^a &\xrightarrow{p} \H_{k}^{a,\ast} :=
    (\textnormal{tr}(\Amk'\Amk))^{1/2} \cdot \D_k^{-1/2}\Upsilon_k'\Z_k^{1/2} ,
\end{align*}
where $\Upsilon_k$ is the eigenvector matrix of $\textnormal{tr}(\Amk'\Amk) \cdot\omega_k^{-1} \Z_k^{1/2} \bSigma_{A,k}\Z_k^{1/2}$.
\end{lemma}
\textbf{\textit{Proof of Lemma \ref{lemma:limit}.}}
First, let $\wh\S,\wh\Q,\wh\D,\H^a$ be from (\ref{eqn: notation_short}), and $\R^\ast,\wt\R,\R_1,\R_2,\R_3$ from (\ref{eqn:S_decomposition}). Define also $\H^{a,\ast}:=\H_{k}^{a,\ast}$, then we have
\begin{equation}
\begin{split}
\label{eqn: D_hat_QRQ}
    \frac{1}{\omega_k}(\wh\D-
    \wh\Q'\R^\ast\wh\Q) &=\frac{1}{\omega_k}\wh\Q'(\wh\S-\R^\ast)\wh\Q \\ &=
    \frac{1}{\omega_k}\wh\Q'(\wt\R-\R^\ast)\wh\Q
    + \frac{1}{\omega_k}\wh\Q'\R_1\wh\Q
    + \frac{1}{\omega_k}\wh\Q'\R_2\wh\Q
    + \frac{1}{\omega_k}\wh\Q'\R_3\wh\Q
    = o_P(1),
\end{split}
\end{equation}
where the last equality follows from the proof of Lemma \ref{lemma:2}.

Using the structure in Assumption (F1), we have
\begin{equation*}
\begin{split}
    \b{E}\Big[\frac{1}{\omega_k}
    \R^\ast\Big] &=
    \frac{1}{\omega_k T}\sum_{t=1}^T
    \b{E}\bigg[\A_k\Big(\sum_{q\geq 0}a_{f,q}\X_{f,t-q}\Big)\A_{\text{-}k}'
    \A_{\text{-}k}\Big(\sum_{q\geq 0}a_{f,q}\X_{f,t-q}'\Big)\A_k'
    \bigg] \\
    &=
    \tr(\A_{\text{-}k}'
    \A_{\text{-}k}) \cdot
    \frac{1}{\omega_k T}\sum_{t=1}^T
    \sum_{q\geq 0}a_{f,q}^2 \A_k\A_k'
    =
    \tr(\A_{\text{-}k}'
    \A_{\text{-}k}) \cdot \frac{1}{\omega_k} \A_k\A_k'.
\end{split}
\end{equation*}
Meanwhile, we have
\begin{equation*}
\begin{split}
    \text{Var}\Big(\Big(\frac{1}{\omega_k}
    \R^\ast\Big)_{ij}\Big) &=
    \frac{1}{\omega_k^2 T^2}\sum_{t=1}^T
    \sum_{s=1}^T\text{cov}\bigg(
    \A_{k,i\cdot}'\Big(\sum_{q\geq 0}a_{f,q}\X_{f,t-q}\Big)\A_{\text{-}k}'
    \A_{\text{-}k}\Big(\sum_{q\geq 0}a_{f,q}\X_{f,t-q}'\Big)
    \A_{k,j\cdot},\\
    &
    \A_{k,i\cdot}'\Big(\sum_{q\geq 0}a_{f,q}\X_{f,s-q}\Big)\A_{\text{-}k}'
    \A_{\text{-}k}\Big(\sum_{q\geq 0}a_{f,q}\X_{f,s-q}'\Big)
    \A_{k,j\cdot}\bigg) \\
    &=
    \frac{1}{\omega_k^2 T^2}
    \sum_{t=1}^T\sum_{q\geq 0}
    \sum_{p\geq 0}a_{f,q}^2
    a_{f,p}^2\cdot \text{Var}\bigg(
    \A_{k,i\cdot}'\X_{f,t-q}\A_{\text{-}k}'
    \A_{\text{-}k}\X_{f,t-p}'
    \A_{k,j\cdot}\bigg)\\
    &=
    \frac{1}{\omega_k^2 T^2}
    \sum_{t=1}^T\sum_{q\geq 0}
    \sum_{p\geq 0}a_{f,q}^2
    a_{f,p}^2\cdot O_P(\|\Amk\|_F^2) \\
    &=
    O_P\Bigg(
    d_k^{-2\alpha_{k,r_k}}
    \frac{1}{T}
    \prod_{j\in[K]\setminus\{k\}}
    d_j^{-\alpha_{j,1}}\Bigg) = o_P(1),
\end{split}
\end{equation*}
where we used Assumption (E2) in the third last equality, both (L1) and (F1) in the second last, and (R1) in the last. We can then conclude $\omega_k^{-1}\R^\ast \xrightarrow{p} \tr(\Amk'\Amk) \cdot \omega_k^{-1} \A_k\A_k'$. Together with (\ref{eqn: D_hat_QRQ}) and Assumption (L2), we obtain the limit of $\omega_k^{-1}\wh\D$ as
\begin{equation}
\label{eqn: D_hat_limit}
    \omega_k^{-1}\wh\D \xrightarrow{p}
    \omega_k^{-1}\D
    = \omega_k^{-1}
    \tr(\Amk'\Amk) \cdot \diag
    \{\lambda_j(\A_k'\A_k)
    \mid j\in[r_k]\} .
\end{equation}

Define further
\begin{equation*}
    \R_\text{res} := \omega_k^{-1}\Z_k^{1/2} \Q'((\wt\R-\R^\ast)+\R_1+\R_2+\R_3) \wh\Q.
\end{equation*}
With similar arguments in the proof of Lemma \ref{lemma:2}, we have $\|\R_\text{res}\|_F = o_P(\|\Z_k^{1/2}\|_F)$. Left-multiply both sides of $\wh\S\wh\Q=\wh\Q\wh\D $ by $\omega_k^{-1}\Z_k^{1/2}\Q'$, we can write
\begin{equation*}
\begin{split}
    (\Z_k^{1/2}\Q'\wh\Q)(\omega_k^{-1}\wh\D)
    &= \omega_k^{-1}\Z_k^{1/2}
    \Q'\wh\S\wh\Q
    =
    \omega_k^{-1} \Z_k^{1/2}\Q'\R^\ast \wh\Q + \R_\text{res}\\
    &=
    \Big[\omega_k^{-1} \Z_k^{1/2}\Q'\R^\ast \wh\Q
    (\Z_k^{1/2}\Q'\wh\Q)^{-1} + \R_\text{res}
    (\Z_k^{1/2}\Q'\wh\Q)^{-1}\Big](\Z_k^{1/2}\Q'\wh\Q).
\end{split}
\end{equation*}
Hence, each column of $\Z_k^{1/2}\Q'\wh\Q$ is an eigenvector of the matrix
\[\omega_k^{-1} \Z_k^{1/2} \Q'\R^\ast \wh\Q(\Z_k^{1/2}\Q'\wh\Q)^{-1} + \R_\text{res}(\Z_k^{1/2}\Q'\wh\Q)^{-1}.\]
We have $(\Z_k^{1/2}\Q'\wh\Q)'(\Z_k^{1/2}\Q'\wh\Q) \xrightarrow{p} (\tr(\Amk'\Amk))^{-1}\cdot\D$, since
\begin{equation*}
\begin{split}
    &\hspace{5mm}
    \omega_k^{-1}\Big[(\Z_k^{1/2}\Q'\wh\Q)'(\Z_k^{1/2}\Q'\wh\Q) - (\tr(\Amk'\Amk))^{-1}\cdot\D \Big]\\
    &= \Big\{ \frac{1}{\omega_k} \wh\Q'\A_k'\A_k \wh\Q - \frac{1}{\omega_k}(\tr(\Amk'\Amk))^{-1} \wh\Q' \R^{\ast} \wh\Q \Big\}
    + \Big\{ \frac{1}{\omega_k}(\tr(\Amk'\Amk))^{-1} \wh\Q'\R^{\ast} \wh\Q  - \frac{1}{\omega_k}(\tr(\Amk'\Amk))^{-1} \, \D \Big\} ,
\end{split}
\end{equation*}
which is $o_P(1)$ from the limit of $\omega_k^{-1}\R^{\ast}$ (for the first square bracket) and from (\ref{eqn: D_hat_QRQ}) and (\ref{eqn: D_hat_limit}) (for the second square bracket).
Hence the eigenvalues of $(\Q'\wh\Q)'(\Q'\wh\Q)$ are asymptotically bounded away from zero and infinity by Assumption (L1), and also $\|(\Z_k^{1/2}\Q'\wh\Q)^{-1}\|_F = O_P(\|\Z_k^{-1/2}\|_F)$.
Let
\begin{equation*}
    \Upsilon_k^\ast :=
    (\tr(\Amk'\Amk))^{1/2}
    \cdot(\Z_k^{1/2}\Q'\wh\Q)\D^{-1/2}.
\end{equation*}
Using the limit of $\omega_k^{-1}\R^\ast$, we have
\begin{equation*}
\begin{split}
    \omega_k^{-1} \Z_k^{1/2}\Q'\R^\ast \wh\Q(\Z_k^{1/2}\Q'\wh\Q)^{-1}
    &\xrightarrow{p}
    \tr(\Amk'\Amk)\cdot \omega_k^{-1} \Z_k^{1/2}\Q'
    \Q \Z_k \Q'
    \wh\Q(\Q'\wh\Q)^{-1}\Z_k^{-1/2}
    \\
    &=
    \tr(\Amk'\Amk) \cdot
    \omega_k^{-1} \Z_k^{1/2}\Q'\Q\Z_k^{1/2}
    \\
    &= \tr(\Amk'\Amk) \cdot
    \omega_k^{-1} \Z_k^{1/2}
    \bSigma_{A,k}\Z_k^{1/2},
\end{split}
\end{equation*}
and $\|\R_\text{res}(\Z_k^{1/2}\Q'\wh\Q)^{-1}\|_F=o_P(1)$ from the above. By Assumption (L2) and eigenvector perturbation theories, there exists a unique eigenvector matrix $\Upsilon_k$ of $\tr(\Amk'\Amk) \cdot\omega_k^{-1} \Z_k^{1/2} \bSigma_{A,k}\Z_k^{1/2}$ such that \linebreak $\|\Upsilon_k - \Upsilon_k^\ast\|=o_P(1)$. Therefore, we have
\begin{equation*}
    \Q'\wh\Q =
    (\tr(\Amk'\Amk))^{-1/2}
    \cdot \Z_k^{-1/2}\Upsilon_k^\ast\D^{1/2}
    \xrightarrow{p}
    (\tr(\Amk'\Amk))^{-1/2}
    \cdot \Z_k^{-1/2}\Upsilon_k\D^{1/2}.
\end{equation*}
Thus, we have
\begin{equation}
\begin{split}
\label{eqn: H_result}
    \H^a &= \wh\D^{-1}\frac{1}{T}
    \sum_{t=1}^T \wh\Q'\A_k\F_t\A_{\text{-}k}'
    \A_{\text{-}k}\F_t'\bf{Z}_k^{1/2}
    =
    \wh\D^{-1}\wh\Q'
    \R^\ast\A_k(\A_k'\A_k)^{-1}
    \Z_k^{1/2}
    = \wh\D^{-1}\wh\Q'
    \R^\ast\Q\bSigma_{A,k}^{-1} \\
    &=
    \tr(\Amk'\Amk) \cdot \D^{-1}\wh\Q'
    \Q\Z_k\Q'\Q\bSigma_{A,k}^{-1} + o_P(1) \\
    &\xrightarrow{p}
    (\tr(\Amk'\Amk))^{1/2} \cdot \D^{-1/2}\Upsilon_k'\Z_k^{1/2}.
\end{split}
\end{equation}
This completes the proof of Lemma \ref{lemma:limit}.
$\square$

\textbf{\textit{Proof of Theorem \ref{thm:asymp_normality_loadings}.}}
Suppose we focus on the $k$-th mode, and hence we adapt all notations by omitting the subscript $k$ for the ease of notational simplicity; see (\ref{eqn: notation_short}) for example. Moreover, we set  $\X_{e,t}:=\mat{k}{\c{X}_{e,t}}, \X_{\epsilon,t}:=\mat{k}{\c{X}_{\epsilon,t}}$ and $\X_{f,t}:=\mat{k}{\c{X}_{f,t}}$.

To proceed, we first decompose
\begin{equation}
\label{eqn: asymp_loading_expansion}
\begin{split}
    \wh\Q_{j\cdot}-\H^a\Q_{j\cdot}=
    (\wh\Q_{j\cdot}-\H_j\Q_{j\cdot}) + (\H_j-\H^a)\Q_{j\cdot} .
\end{split}
\end{equation}
Consider the first term $(\wh\Q_{j\cdot}-\H_j\Q_{j\cdot})$. Using the decomposition in (\ref{eqn: Q_j_consistency}),
\begin{equation}
\begin{split}
\label{eqn: Q_expansion}
    \wh\Q_{j\cdot}-\H_j\Q_{j\cdot} &=
    \wh{\bf{D}}^{-1}\sum_{i=1}^{d_k} (\wh{\bf{Q}}_{i\cdot}-\H_j\Q_{i\cdot})
    \sum_{h=1}^{d_{\text{-}k}}
    \frac{1}{|\psi_{ij,h}|}
    \sum_{t\in\psi_{ij,h}}
        E_{t,jh} (\A_{\text{-}k})_{h\cdot}'
        \F_t' \A_{i\cdot} \\
    &+
    \wh{\bf{D}}^{-1}\sum_{i=1}^{d_k} (\wh{\bf{Q}}_{i\cdot}-\H_j\Q_{i\cdot})
    \sum_{h=1}^{d_{\text{-}k}}
    \frac{1}{|\psi_{ij,h}|}
    \sum_{t\in\psi_{ij,h}}
        E_{t,ih} (\A_{\text{-}k})_{h\cdot}'
        \F_t' \A_{j\cdot} \\
    &+
    \wh{\bf{D}}^{-1}\sum_{i=1}^{d_k} (\wh{\bf{Q}}_{i\cdot}-\H_j\Q_{i\cdot})
    \sum_{h=1}^{d_{\text{-}k}}
    \frac{1}{|\psi_{ij,h}|}
    \sum_{t\in\psi_{ij,h}}
        E_{t,ih} E_{t,jh} \\
    &+
    \c{I}_{H,j}+\c{I}_{H,j}^\ast
    +\c{II}_{H,j}
    +\c{III}_{H,j} ,\;\text{ where }
\end{split}
\end{equation}
\vspace{-24pt}
\begin{align*}
    \c{I}_{H,j} &:=
    \wh{\bf{D}}^{-1}\sum_{i=1}^{d_k} \H^a\Q_{i\cdot}
    \sum_{h=1}^{d_{\text{-}k}}
    \frac{1}{|\psi_{ij,h}|}
    \sum_{t\in\psi_{ij,h}}
        E_{t,jh} (\A_{\text{-}k})_{h\cdot}'
        \F_t' \A_{i\cdot}, \\
    \c{I}_{H,j}^\ast &:=
    \wh{\bf{D}}^{-1}\sum_{i=1}^{d_k}
    (\H_j-\H^a)\Q_{i\cdot}
    \sum_{h=1}^{d_{\text{-}k}}
    \frac{1}{|\psi_{ij,h}|}
    \sum_{t\in\psi_{ij,h}}
        E_{t,jh} (\A_{\text{-}k})_{h\cdot}'
        \F_t' \A_{i\cdot}, \\
    \c{II}_{H,j} &:=
    \wh{\bf{D}}^{-1}\sum_{i=1}^{d_k} \H_j\Q_{i\cdot}
    \sum_{h=1}^{d_{\text{-}k}}
    \frac{1}{|\psi_{ij,h}|}
    \sum_{t\in\psi_{ij,h}}
        E_{t,ih} (\A_{\text{-}k})_{h\cdot}'
        \F_t' \A_{j\cdot} , \\
    \c{III}_{H,j} &:=
    \wh{\bf{D}}^{-1}\sum_{i=1}^{d_k}
    \H_j\Q_{i\cdot}
    \sum_{h=1}^{d_{\text{-}k}}
    \frac{1}{|\psi_{ij,h}|}
    \sum_{t\in\psi_{ij,h}}
        E_{t,ih} E_{t,jh}, \; \text{ with } \A_{\text{-}k} :=
    \otimes_{l\in[K]\setminus\{k\}}\A_l.
\end{align*}
We want to show that $\c{I}_{H,j}$ is the leading term among those in (\ref{eqn: Q_expansion}). To this end, we will show $\sqrt{T\omega_B} \cdot\c{I}_{H,j}$ converges to a normal distribution with mean zero and variance of constant order (see (\ref{eqn: lim_I_H_j}) later), so that $\c{I}_{H,j}$ is of order $(T\omega_B)^{-1/2}$ exactly. Then it suffices to show that the rate $(T\omega_B)^{-1}$ is dominating the following rates multiplied by the rate of $\|\wh\D^{-1} \|_F^2 =O_P\Big(d_k^{2(\alpha_{k,1} -\alpha_{k,r_k})} g_s^{-2} \Big)$ from Lemma \ref{lemma:2}:
\begin{align}
    &\bigg\|
    \sum_{i=1}^{d_k}(\H_j-\H^a)\Q_{i\cdot}
    \sum_{h=1}^{d_{\text{-}k}}
    \frac{1}{|\psi_{ij,h}|}
    \sum_{t\in\psi_{ij,h}}
    E_{t,jh} (\A_{\text{-}k})_{h\cdot}'
    \F_t' \A_{i\cdot}\bigg\|^2 ,
    \label{eqn: I_H_diff_to_prove}\\
    & \bigg\|\sum_{i=1}^{d_k}\H_j\Q_{i\cdot}
    \sum_{h=1}^{d_{\text{-}k}}
    \frac{1}{|\psi_{ij,h}|}
    \sum_{t\in\psi_{ij,h}}
    E_{t,ih} (\A_{\text{-}k})_{h\cdot}'
    \F_t' \A_{j\cdot}\bigg\|^2 ,
    \label{eqn: II_H_to_prove}\\
    &
    \bigg\|\sum_{i=1}^{d_k}\H_j\Q_{i\cdot}
    \sum_{h=1}^{\dmk}
    \frac{1}{|\psi_{ij,h}|}
    \sum_{t\in\psi_{ij,h}}
    E_{t,ih} E_{t,jh}\bigg\|^2 ,
    \label{eqn: III_H_to_prove}\\
    &\bigg\|
    \sum_{i=1}^{d_k}(\wh{\bf{Q}}_{i\cdot}-\H_j\Q_{i\cdot})
    \sum_{h=1}^{d_{\text{-}k}}
    \frac{1}{|\psi_{ij,h}|}
    \sum_{t\in\psi_{ij,h}}
    E_{t,jh} (\A_{\text{-}k})_{h\cdot}'
    \F_t' \A_{i\cdot}\bigg\|^2 ,
    \label{eqn: I_HQ_diff_to_prove}\\
    &
    \bigg\|\sum_{i=1}^{d_k}(\wh{\bf{Q}}_{i\cdot}-\H_j\Q_{i\cdot})
    \sum_{h=1}^{d_{\text{-}k}}
    \frac{1}{|\psi_{ij,h}|}
    \sum_{t\in\psi_{ij,h}}
    E_{t,ih} (\A_{\text{-}k})_{h\cdot}'
    \F_t' \A_{j\cdot}\bigg\|^2 ,
    \label{eqn: II_HQ_diff_to_prove}\\
    &
    \bigg\|\sum_{i=1}^{d_k}(\wh{\bf{Q}}_{i\cdot}-\H_j\Q_{i\cdot})
    \sum_{h=1}^{\dmk}
    \frac{1}{|\psi_{ij,h}|}
    \sum_{t\in\psi_{ij,h}}
    E_{t,ih} E_{t,jh}\bigg\|^2 .
    \label{eqn: III_HQ_diff_to_prove}
\end{align}

Note that we can easily see the rates of (\ref{eqn: II_H_to_prove}) and (\ref{eqn: III_H_to_prove}) are greater than those of (\ref{eqn: II_HQ_diff_to_prove}) and (\ref{eqn: III_HQ_diff_to_prove}) respectively, using Lemma \ref{lemma:2} and the Cauchy--Schwarz inequality.

Consider (\ref{eqn: I_H_diff_to_prove}) first. We have
\begin{equation*}
\begin{split}
    &\hspace{5mm} \bigg\|
    \sum_{i=1}^{d_k}(\H_j-\H^a)\Q_{i\cdot}
    \sum_{h=1}^{d_{\text{-}k}}
    \frac{1}{|\psi_{ij,h}|}
    \sum_{t\in\psi_{ij,h}}
    E_{t,jh} (\A_{\text{-}k})_{h\cdot}'
    \F_t' \A_{i\cdot}\bigg\|^2 \\
    &\leq
    \bigg( \sum_{i=1}^{d_k}
    \|\H_j-\H^a\|_F^2 \cdot
    \|\Q_{i\cdot}\|^2 \bigg)
    \cdot \sum_{i=1}^{d_k} \bigg(
    \sum_{h=1}^{d_{\text{-}k}}
        \frac{1}{|\psi_{ij,h}|}
        \sum_{t\in\psi_{ij,h}}
        E_{t,jh} (\A_{\text{-}k})_{h\cdot}'
    \F_t' \A_{i\cdot} \bigg)^2 \\
    &=
    O(\dmk^2) \cdot \|\H_j-\H^a\|_F^2
    \cdot \sum_{i=1}^{d_k}
    \bigg(\sum_{h=1}^{\dmk}\sum_{t\in\psi_{ij,h}}
    \frac{1}{\dmk \cdot |\psi_{ij,h}|}
    E_{t,jh} (\A_{\text{-}k})_{h\cdot}'
    \F_t' \A_{i\cdot} \bigg)^2 \\
    &=
    O(\dmk^2) \cdot \|\H_j-\H^a\|_F^2
    \cdot \sum_{i=1}^{d_k} \|\bf{u}_i\|^2
    \bigg(\sum_{h=1}^{\dmk}\sum_{t\in\psi_{ij,h}}
    \frac{1}{\dmk \cdot |\psi_{ij,h}|}
    E_{t,jh} \bf{v}_h'
    \F_t' \frac{1}{\|\bf{u}_i\|}\bf{u}_i\bigg)^2 \\
    &=
    O_P\Bigg(d_k^{2(\alpha_{k,1}-\alpha_{k,r_k})}
    \Big( \frac{\dmk}{T^2} \Big)
    d_k^{\alpha_{k,1}} \Bigg),
\end{split}
\end{equation*}
where $\bf{v}_h := (\Amk)_{h\cdot}, \; \bf{u}_i := \bf{A}_{i\cdot}$, and we used Lemma \ref{lemma:3}, Proposition \ref{Prop:assumption_implications}.2 and (L1) in the last equality.

To bound (\ref{eqn: II_H_to_prove}), note from Assumption (E1) and (E2) that we can write
\begin{equation*}
    E_{t,ih} = \sum_{q\geq 0}a_{e,q}
    \A_{e,k,i\cdot}'\X_{e,t-q}\bf{A}_{e,\text{-}k,h\cdot}+
    [\mat{k}{\bSigma_\epsilon}]_{ih}\sum_{q\geq 0}a_{\epsilon,q}
    (\X_{\epsilon,t-q})_{ih}.
\end{equation*}
Consider first $\sum_{h=1}^{\dmk}\sum_{t\in\psi_{ij,h}}(\sum_{q\geq 0}a_{e,q}\A_{e,k,i\cdot}'\X_{e,t-q}\bf{A}_{e,\text{-}k,h\cdot}) (\A_{\text{-}k})_{h\cdot}'\F_t' \A_{j\cdot}$. By Assumption (O1), (E1), (E2) and (F1), we have
\begin{equation}
\begin{split}
\label{eqn: II_H_to_prove_step1}
    & \hspace{5mm}
    \b{E}\Big[\big(
    \sum_{h=1}^{\dmk}
    \sum_{t\in\psi_{ij,h}}
    (\sum_{q\geq 0} a_{e,q}\A_{e,k,i\cdot}'\X_{e,t-q}\bf{A}_{e,\text{-}k,h\cdot}) (\A_{\text{-}k})_{h\cdot}'
    \F_t' \A_{j\cdot}\big)^2\Big] \\
    &=
    \text{cov}\Big(
    \sum_{h=1}^{\dmk}
    \sum_{t\in\psi_{ij,h}}
    (\A_{\text{-}k})_{h\cdot}'
    \big(\sum_{q\geq 0}a_{f,q}\X_{f,t-q}'\big)\A_{j\cdot}
    \big(\sum_{q\geq 0}a_{e,q}
    \A_{e,k,i\cdot}'\X_{e,t-q}\bf{A}_{e,\text{-}k,h\cdot}\big)
    , \\
    &
    \sum_{h=1}^{\dmk}
    \sum_{t\in\psi_{ij,h}}
    (\A_{\text{-}k})_{h\cdot}'
    \big(\sum_{q\geq 0}a_{f,q}\X_{f,t-q}'\big)\A_{j\cdot}
    \big(\sum_{q\geq 0}a_{e,q}
    \A_{e,k,i\cdot}'\X_{e,t-q}\bf{A}_{e,\text{-}k,h\cdot}\big) \Big) \\
    & =
    \sum_{h=1}^{\dmk}\sum_{l=1}^{\dmk}
    \sum_{t\in\psi_{ij,h}\cap \psi_{ij,l}}
    \sum_{q\geq 0}
    a_{f,q}^2a_{e,q}^2 \cdot
    \|\A_{j\cdot}\|^2 \cdot
    \| (\A_{\text{-}k})_{h\cdot}\|\cdot
    \| (\A_{\text{-}k})_{l\cdot}\|\cdot
    \|\A_{e,\text{-}k,h\cdot}\|\cdot
    \|\A_{e,\text{-}k,l\cdot}\|\cdot
    \|\A_{e,k,i\cdot}\|^2 \\
    &
    =
    O(T)\cdot \|\A_{j\cdot}\|^2 \cdot
    \|\A_{e,k,i\cdot}\|^2 .
\end{split}
\end{equation}
Consider also $\sum_{i=1}^{d_k}\sum_{h=1}^{\dmk}\sum_{t\in\psi_{ij,h}}\Q_{i\cdot} ([\mat{k}{\bSigma_\epsilon}]_{ih}\sum_{q\geq 0}a_{\epsilon,q} (\X_{\epsilon,t-q})_{ih}) (\A_{\text{-}k})_{h\cdot}'\F_t' \A_{j\cdot}$. Similarly, by Assumption (O1), (E1), (E2) and (F1), we have
\begin{equation}
\begin{split}
\label{eqn: II_H_to_prove_step2}
    & \hspace{5mm}
    \b{E}\Big[\big\|
    \sum_{i=1}^{d_k}
    \sum_{h=1}^{\dmk}
    \sum_{t\in\psi_{ij,h}}
    \Q_{i\cdot}([\mat{k}{\bSigma_\epsilon}]_{ih}\sum_{q\geq 0}a_{\epsilon,q} (\X_{\epsilon,t-q})_{ih})(\A_{\text{-}k})_{h\cdot}'
    \F_t' \A_{j\cdot}\big\|^2\Big] \\
    &=
    \text{cov}\Big(
    \sum_{i=1}^{d_k}
    \sum_{h=1}^{\dmk}
    \sum_{t\in\psi_{ij,h}}
    \Q_{i\cdot}
    (\A_{\text{-}k})_{h\cdot}'
    \big(\sum_{q\geq 0}a_{f,q}\X_{f,t-q}'\big)\A_{j\cdot}
    \big( [\mat{k}{\bSigma_\epsilon}]_{ih}\sum_{q\geq 0}a_{\epsilon,q}
    (\X_{\epsilon,t-q})_{ih}\big)
    , \\
    &
    \sum_{i=1}^{d_k}
    \sum_{h=1}^{\dmk}
    \sum_{t\in\psi_{ij,h}}
    \Q_{i\cdot}
    (\A_{\text{-}k})_{h\cdot}'
    \big(\sum_{q\geq 0}a_{f,q}\X_{f,t-q}'\big)\A_{j\cdot}
    \big([\mat{k}{\bSigma_\epsilon}]_{ih}\sum_{q\geq 0}a_{\epsilon,q}
    (\X_{\epsilon,t-q})_{ih}\big) \Big) \\
    & =
    \sum_{i=1}^{d_k}
    \sum_{h=1}^{\dmk}
    \sum_{t\in\psi_{ij,h}}\sum_{q\geq 0}
    a_{f,q}^2a_{\epsilon,q}^2
    \cdot \|\A_{j\cdot}\|^2 \cdot \|
    (\A_{\text{-}k})_{h\cdot}\|^2
    \cdot \Sigma_{\epsilon,k,h,ii}
    \cdot \|\Q_{i\cdot}\|^2\\
    & =
    O(T)\cdot \|\A_{j\cdot}\|^2 \cdot \|
    \A_{\text{-}k}\|^2 \cdot
    \|\Q\|^2 .
\end{split}
\end{equation}
Hence it holds that
\begin{equation*}
\begin{split}
    &\hspace{5mm} \bigg\|
    \sum_{i=1}^{d_k}\H_j\Q_{i\cdot}
    \sum_{h=1}^{d_{\text{-}k}}
    \frac{1}{|\psi_{ij,h}|}
    \sum_{t\in\psi_{ij,h}}
    E_{t,ih} (\A_{\text{-}k})_{h\cdot}'
    \F_t' \A_{j\cdot}\bigg\|^2 \\
    &\leq
    \|\H_j\|_F^2 \cdot \bigg\|
    \sum_{i=1}^{d_k}
    \sum_{h=1}^{\dmk}
    \frac{1}{|\psi_{ij,h}|}
    \sum_{t\in\psi_{ij,h}}
    \Q_{i\cdot}([\mat{k}{\bSigma_\epsilon}]_{ih}\sum_{q\geq 0}a_{\epsilon,q} (\X_{\epsilon,t-q})_{ih})(\A_{\text{-}k})_{h\cdot}'
    \F_t' \A_{j\cdot}\bigg\|^2
    \\
    &+
    \|\H_j\|_F^2 \cdot
    \bigg( \sum_{i=1}^{d_k}
    \|\Q_{i\cdot}\|^2 \bigg)\cdot
    \sum_{i=1}^{d_k} \bigg\|
    \sum_{h=1}^{\dmk}
    \frac{1}{|\psi_{ij,h}|}
    \sum_{t\in\psi_{ij,h}}(\sum_{q\geq 0}a_{e,q}\A_{e,k,i\cdot}'\X_{e,t-q}\bf{A}_{e,\text{-}k,h\cdot}) (\A_{\text{-}k})_{h\cdot}'\F_t' \A_{j\cdot}
    \bigg\|^2 \\
    &=
    O_P\Bigg(
    d_k^{-\alpha_{k,1}}
    \bigg( \frac{1}{T} \bigg)
    \prod_{j=1}^Kd_j^{\alpha_{j,1}}
    \Bigg),
\end{split}
\end{equation*}
where we used Assumption (L1), (\ref{eqn: II_H_to_prove_step1}) and (\ref{eqn: II_H_to_prove_step2}) in the last equality.

For (\ref{eqn: III_H_to_prove}), by Assumption (O1), (E1) and (E2), we have from the proof of Proposition \ref{Prop:assumption_implications} that
\begin{equation*}
\begin{split}
    &\hspace{5mm}
    \text{Var}\Big(
    \sum_{i=1}^{d_k}
    \sum_{h=1}^{\dmk}
    \sum_{t\in\psi_{ij,h}}
    \Q_{i\cdot}
    E_{t,ih} E_{t,jh}\Big) \\
    & =
    O(1) \cdot
    \sum_{i=1}^{d_k}\sum_{w=1}^{d_k}
    \sum_{h=1}^{\dmk}\sum_{l=1}^{\dmk}
    \sum_{t\in\psi_{ij,h}\cap\psi_{wj,l}}
    \sum_{n=1}^{r_{e,k}}
    \sum_{m=1}^{r_{e,\text{-}k}}
    \sum_{q\geq 0}
    a_{e,q}^4 A_{e,k,in}
    A_{e,k,wn} A_{e,k,jn}^2
    A_{e,\text{-}k,hm}^2
    A_{e,\text{-}k,lm}^2 \\
    &\cdot
    \|\Q_{i\cdot}\|
    \cdot \|\Q_{w\cdot}\|
    \cdot \text{Var}((\X_{e,t-q})_{nm}^2)\\
    & +
    O(1) \cdot
    \sum_{i=1}^{d_k}
    \sum_{h=1}^{\dmk}
    \sum_{t\in\psi_{ij,h}\cap\psi_{wj,l}}
    \sum_{q\geq 0} a_{\epsilon,q}^4
    \Sigma_{\epsilon,k,h,ii}
    \Sigma_{\epsilon,k,h,jj}
    \cdot \|\Q_{i\cdot}\|^2
    \cdot \text{Var}(
    (\X_{\epsilon,t-q})_{ih}
    (\X_{\epsilon,t-q})_{jh}) \\
    &
    = O(T + T \dmk)=O(T \dmk).
\end{split}
\end{equation*}
Moreover, it holds that
\begin{equation*}
\begin{split}
    \b{E}\Big[
    \sum_{i=1}^{d_k}
    \sum_{h=1}^{\dmk}
    \sum_{t\in\psi_{ij,h}}
    E_{t,ih} E_{t,jh}\Big] &=
    \sum_{i=1}^{d_k}\sum_{h=1}^{\dmk}
    \sum_{t\in\psi_{ij,h}}\bigg(
    \|\A_{e,\text{-}k,h\cdot}\|^2
    \cdot \|\A_{e,k,i\cdot}\|
    \cdot \|\A_{e,k,j\cdot}\|
    + \Sigma_{\epsilon,k,h,ij}\bigg)
    =O(T\dmk),
\end{split}
\end{equation*}
and with $\max_i\|\Q_{i\cdot}\|^2\leq\|\A_{j\cdot}\|^2\cdot \|\Z_k^{-1/2}\|^2=O_P\big(d_k^{-\alpha_{k,r_k}}\big)$, we thus have
\begin{equation*}
\begin{split}
    \bigg\|\sum_{i=1}^{d_k}\H_j\Q_{i\cdot}
    \sum_{h=1}^{\dmk}
    \frac{1}{|\psi_{ij,h}|}
    \sum_{t\in\psi_{ij,h}}
    E_{t,ih} E_{t,jh}\bigg\|^2
    &\leq
    \|\H_j\|_F^2 \cdot
    \bigg\|\sum_{i=1}^{d_k}\Q_{i\cdot}
    \sum_{h=1}^{\dmk}
    \frac{1}{|\psi_{ij,h}|}
    \sum_{t\in\psi_{ij,h}}
    E_{t,ih} E_{t,jh}\bigg\|^2 \\
    &=
    O\bigg(\frac{\dmk}{T}+\dmk^2
    d_k^{-\alpha_{k,r_k}}\bigg).
\end{split}
\end{equation*}

Now consider (\ref{eqn: I_HQ_diff_to_prove}). Similar to (\ref{eqn: I_H_diff_to_prove}), we have
\begin{equation*}
\begin{split}
    &\hspace{5mm} \bigg\|
    \sum_{i=1}^{d_k}(\wh{\bf{Q}}_{i\cdot}-\H_j\Q_{i\cdot})
    \sum_{h=1}^{d_{\text{-}k}}
    \frac{1}{|\psi_{ij,h}|}
    \sum_{t\in\psi_{ij,h}}
    E_{t,jh} (\A_{\text{-}k})_{h\cdot}'
    \F_t' \A_{i\cdot}\bigg\|^2 \\
    &\leq
    \bigg( \sum_{i=1}^{d_k}
    \|\wh{\bf{Q}}_{i\cdot}-
    \H_j\Q_{i\cdot}\|_F^2 \bigg)
    \cdot \sum_{i=1}^{d_k} \bigg(
    \sum_{h=1}^{d_{\text{-}k}}
        \frac{1}{|\psi_{ij,h}|}
        \sum_{t\in\psi_{ij,h}}
        E_{t,jh} (\A_{\text{-}k})_{h\cdot}'
    \F_t' \A_{i\cdot} \bigg)^2 \\
    &=
    O(\dmk^2) \cdot
    \bigg( \sum_{i=1}^{d_k}
    \|\wh{\bf{Q}}_{i\cdot}-
    \H_j\Q_{i\cdot}\|_F^2 \bigg)
    \cdot \sum_{i=1}^{d_k} \|\bf{u}_i\|^2
    \bigg(\sum_{h=1}^{\dmk}\sum_{t\in\psi_{ij,h}}
    \frac{1}{\dmk \cdot |\psi_{ij,h}|}
    E_{t,jh} \bf{v}_h'
    \F_t' \frac{1}{\|\bf{u}_i\|}\bf{u}_i\bigg)^2 \\
    &=
    O_P\Bigg(d_k^{3\alpha_{k,1} - 2\alpha_{k,r_k}}
    \bigg( \frac{\dmk}{T} \bigg)
    \bigg(\frac{1}{T\dmk}+\frac{1}{d_k}
    \bigg)\prod_{j=1}^Kd_j^{2(1-\alpha_{j,1})}
    \Bigg),
\end{split}
\end{equation*}
where we used Lemma \ref{lemma:2}, Proposition \ref{Prop:assumption_implications}.2 and (L1) in the last equality.

Finally, we consider the following ratios with $d_1,\dots,d_K,T \to\infty$:
\begin{equation*}
\begin{split}
    &\hspace{5mm}
    \bigg\|
    \sum_{i=1}^{d_k}(\H_j-\H^a)\Q_{i\cdot}
    \sum_{h=1}^{d_{\text{-}k}}
    \frac{1}{|\psi_{ij,h}|}
    \sum_{t\in\psi_{ij,h}}
    E_{t,jh} (\A_{\text{-}k})_{h\cdot}'
    \F_t' \A_{i\cdot}\bigg\|^2
    d_k^{2(\alpha_{k,1} -\alpha_{k,r_k})} g_s^{-2}\bigg/ (T\omega_B)^{-1} \\
    &=
    O_P\Bigg(
    d_k^{2(\alpha_{k,1}-\alpha_{k,r_k})}
    \cdot \frac{1}{T}\Bigg) = o_P(1), \\
    &\hspace{5mm}
    \bigg\|
    \sum_{i=1}^{d_k}\H_j\Q_{i\cdot}
    \sum_{h=1}^{d_{\text{-}k}}
    \frac{1}{|\psi_{ij,h}|}
    \sum_{t\in\psi_{ij,h}}
    E_{t,ih} (\A_{\text{-}k})_{h\cdot}'
    \F_t' \A_{j\cdot}\bigg\|^2
    d_k^{2(\alpha_{k,1} -\alpha_{k,r_k})} g_s^{-2}\bigg/ (T\omega_B)^{-1} \\
    &=
    O_P\Big(
    d_k^{-\alpha_{k,1}}
    \prod_{j\in[K]\setminus\{k\}} d_j^{\alpha_{j,1}-1}
    \Big)
    = o_P(1), \\
    &\hspace{5mm}
    \bigg\|\sum_{i=1}^{d_k}\H_j\Q_{i\cdot}
    \sum_{h=1}^{\dmk}
    \frac{1}{|\psi_{ij,h}|}
    \sum_{t\in\psi_{ij,h}}
    E_{t,ih} E_{t,jh}\bigg\|^2
    d_k^{2(\alpha_{k,1} -\alpha_{k,r_k})} g_s^{-2}\bigg/ (T\omega_B)^{-1} \\
    &=
    O_P\Big(
    T\dmk d_k^{-\alpha_{k,r_k}-\alpha_{k,1}}
    \Big)
    = o_P(1), \\
    &\hspace{5mm}
    \bigg\|
    \sum_{i=1}^{d_k}(\wh{\bf{Q}}_{i\cdot}-\H_j\Q_{i\cdot})
    \sum_{h=1}^{d_{\text{-}k}}
    \frac{1}{|\psi_{ij,h}|}
    \sum_{t\in\psi_{ij,h}}
    E_{t,jh} (\A_{\text{-}k})_{h\cdot}'
    \F_t' \A_{i\cdot}\bigg\|^2
    d_k^{2(\alpha_{k,1} -\alpha_{k,r_k})} g_s^{-2}\bigg/ (T\omega_B)^{-1} \\
    &=
    O_P\Bigg(d_k^{2(\alpha_{k,1} - \alpha_{k,r_k})}
    \bigg(\frac{1}{T\dmk}+\frac{1}{d_k}
    \bigg)\prod_{j=1}^Kd_j^{2(1-\alpha_{j,1})}
    \Bigg)
    = o_P(1),
\end{split}
\end{equation*}
by Assumptions (R1) and the rate assumptions $d_k^{2\alpha_{k,1} - 3\alpha_{k,r_k}} = o(\dmk)$ (from the statement of Theorem \ref{thm:imputation_consistency}) and $T\dmk = o(d_k^{\alpha_{k,1} + \alpha_{k,r_k}})$.
Hence $\c{I}_{H,j}$ is indeed the dominating term in (\ref{eqn: Q_expansion}). In other words, we have
\begin{equation}
\label{eqn: I_H_j_dominance}
    \wh\Q_{j\cdot}-\H_j\Q_{j\cdot}
    = \c{I}_{H,j} + o_P(1).
\end{equation}

Then we want to show that
\begin{equation}
\begin{split}
\label{eqn: lim_I_H_j}
    \sqrt{T\omega_B} \cdot
    \c{I}_{H,j} &=
    \sqrt{T\omega_B} \cdot
    \wh\D^{-1}\H^a
    \sum_{i=1}^{d_k} \Q_{i\cdot}
    \sum_{h=1}^{d_{\text{-}k}}
    \frac{1}{|\psi_{ij,h}|}
    \sum_{t\in\psi_{ij,h}}
        E_{t,jh} (\A_{\text{-}k})_{h\cdot}'
        \F_t' \A_{i\cdot} \\
    & \xrightarrow{p}
    \sqrt{T\omega_B} \cdot
    \D^{-1}\H_k^{a,\ast}
    \sum_{i=1}^{d_k} \Q_{i\cdot}
    \sum_{h=1}^{d_{\text{-}k}}
    \frac{1}{|\psi_{ij,h}|}
    \sum_{t\in\psi_{ij,h}}
        E_{t,jh} (\A_{\text{-}k})_{h\cdot}'
        \F_t' \A_{i\cdot} \\
    & \xrightarrow{\c{D}}
    \cN(\0, T\omega_B \cdot
    \D^{-1}\H_k^{a,\ast}
    \bf\Xi_{k,j}(\H_k^{a,\ast})'\D^{-1}),
\end{split}
\end{equation}
where $\D$ and $\H_k^{a,\ast}$ are from Lemma \ref{lemma:limit}, and we require Assumption (AD1) for the covariance matrix to have constant rate. In fact, using Lemma \ref{lemma:limit} and Proposition \ref{Prop:assumption_implications}, the upper bound is of constant order by
\begin{equation*}
\begin{split}
    &\hspace{5mm} \bigg\| T\omega_B \cdot
    \D^{-1}\H_k^{a,\ast}
    \bf\Xi_{k,j}(\H_k^{a,\ast})'\D^{-1} \bigg\|_F \\
    &=
    O\big(T\omega_B \big) \cdot \big\|\D^{-1} \big\|_F^2 \cdot \bigg\| \sum_{i=1}^{d_k}\Q_{i\cdot} \sum_{h=1}^{d_{\text{-}k}} \frac{1}{|\psi_{ij,h}|} \sum_{t\in\psi_{ij,h}} E_{t,jh} (\A_{\text{-}k})_{h\cdot}' \F_t' \A_{i\cdot}\bigg\|^2 \\
    &=
    O\Big(T (\dmk d_k^{\alpha_{k,1}})^{-1} \Big) \cdot \bigg( \sum_{i=1}^{d_k}
    \|\Q_{i\cdot}\|^2 \bigg)\cdot
    \sum_{i=1}^{d_k} \bigg(
    \sum_{h=1}^{d_{\text{-}k}}
        \frac{1}{|\psi_{ij,h}|}
        \sum_{t\in\psi_{ij,h}}
        E_{t,jh} (\A_{\text{-}k})_{h\cdot}'
    \F_t' \A_{i\cdot} \bigg)^2 =
    O_P(1).
\end{split}
\end{equation*}

We will adapt the central limit theorem for $\alpha$-mixing processes (\cite{Fan_Yao_2003}, Theorem 2.21). Due to the existence of missing data and the general missing patterns that we allow, we construct an auxiliary time series to facilitate the proof. Formally, define $\{\B_{j,t}\}_{t\in[T]}$ as
\begin{equation}
\begin{split}
\label{eqn: B_j_t_def}
    \B_{j,t} &:= \sqrt{\omega_B} \cdot
    \D^{-1} \H_k^{a,\ast}
    \sum_{i=1}^{d_k}\sum_{h=1}^{\dmk}
    \Q_{i\cdot}\frac{T}{|\psi_{ij,h}|} \cdot E_{t,jh} (\Amk)_{h\cdot}'
    \F_t' \A_{i\cdot} \cdot \b{1}\{t\in\psi_{ij,h}\}.
\end{split}
\end{equation}
Hence we have the following,
\begin{equation*}
\begin{split}
    \sqrt{T\omega_B} \cdot
    \c{I}_{H,j} &\xrightarrow{p}
    \frac{1}{\sqrt{T}}
    \sum_{t=1}^T \B_{j,t}.
\end{split}
\end{equation*}
It is easy to see that $\b{E}[\B_{j,t}]=\0$ by Assumption (E1), (E2) and (F1). Moreover, $\B_{j,t}$ is also $\alpha$-mixing over $t$. To see this, consider
\begin{equation*}
\begin{split}
    E_{t,jh} (\Amk)_{h\cdot}'
    \F_t' \A_{i\cdot} &=
    \big\{ \sum_{q\geq 0}a_{e,q}
    \A_{e,k,i\cdot}'\X_{e,t-q}\A_{e,\text{-}k,h\cdot}+
    \Sigma_{\epsilon,k,h,ii}^{1/2}
    \sum_{q\geq 0} a_{\epsilon,q}
    (\X_{\epsilon,t-q})_{ih} \big\}
    \A_{\text{-}k, h\cdot}'
    \big( \sum_{q\geq 0} a_{f,q}\X_{f,t-q}'
    \big) \A_{i\cdot} .
\end{split}
\end{equation*}
If we define $\bf{b}_{e,t} :=\sum_{q\geq 0}a_{e,q}\X_{e,t-q}$, $\bf{b}_{\epsilon,ih,t} :=\sum_{q\geq 0}a_{\epsilon,q}(\X_{\epsilon,t-q})_{ih}$ and $\bf{b}_{f,t} :=\sum_{q\geq 0}a_{f,q}\X_{f,t-q}$ which are independent of each other by Assumption (E2), we can then rewrite
\begin{equation*}
    \B_{j,t} =h\Big(\bf{b}_{e,t},
    (\bf{b}_{\epsilon,ih,t})_{i\in
    [d_k],h\in[\dmk]}, \bf{b}_{f,t}\Big),
\end{equation*}
for some function $h$, and hence Theorem 5.2 in \cite{Bradley_2005} implies the $\alpha$-mixing property. Then similar to \cite{ChenFan2023}, it is left to show that there exists an $m>2$ such that $\b{E}[\|\B_{j,t}\|^m]\leq C$ for some constant $C$. With Assumption (E1), (E2) and (F1) and similar to the proof of Proposition \ref{Prop:assumption_implications}, we have
\begin{equation*}
\begin{split}
    \b{E}\bigg[\sum_{h=1}^{\dmk}
    \frac{T}{|\psi_{ij,h}|} \cdot E_{t,jh} \bf{u}'
    \F_t \bf{v} \cdot \b{1}\{t\in\psi_{ij,h}\} \bigg]^2
    = O(\dmk),
\end{split}
\end{equation*}
where $\bf{u}\in\b{R}^{r_k}$ and $\bf{v}\in\b{R}^{\rmk}$ are any deterministic vectors of constant order. Hence
\begin{equation*}
\begin{split}
    &\hspace{5mm}
    \b{E}[\|\B_{j,t}\|^m] \\
    &\leq
    \omega_B^{m/2} \, \|\D^{-1}\|_F^m
    \, \|\H_k^{a,\ast}\|_F^m
    \, \bigg(\sum_{i=1}^{d_k}
    \|\Q_{i\cdot}\|^2\bigg)^{m/2}
    \b{E}\Bigg\{\bigg[\sum_{i=1}^{d_k}
    \bigg(\sum_{h=1}^{\dmk}
    \frac{T}{|\psi_{ij,h}|} \, E_{t,jh} (\Amk)_{h\cdot}'
    \F_t' \A_{i\cdot} \cdot \b{1}\{t\in\psi_{ij,h}\}
    \bigg)^2 \bigg]^{m/2}\Bigg\} \\
    &=
    O\Big( (\omega_B \dmk
    d_k^{\alpha_{k,1}})^{m/2}\Big)
    \cdot\|\D^{-1}\|_F^m
    =
    O_P\Bigg(\bigg(
    \omega_B \dmk
    d_k^{3\alpha_{k,1}-2\alpha_{k,r_k}}
    \prod_{j=1}^Kd_j^{-2\alpha_{j,1}}
    \bigg)^{m/2} \Bigg)
    = O_P(1),
\end{split}
\end{equation*}
where we used Lemma \ref{lemma:2} and the definition of $\omega_B$ in the last line. Theorem 2.21 in \cite{Fan_Yao_2003} then applies. With (\ref{eqn: I_H_j_dominance}), (\ref{eqn: lim_I_H_j}) and Lemma \ref{lemma:limit}, we can directly establish that
\begin{equation}
\label{eqn: asymp_loading_part1}
    \sqrt{T\omega_B}
    \cdot
    (\wh\Q_{j\cdot}-\H_j\Q_{j\cdot})
    \xrightarrow{\c{D}}
    \cN(\0, T\omega_B \cdot
    \D^{-1}\H_k^{a,\ast}
    \bf\Xi_{k,j}(\H_k^{a,\ast})'\D^{-1}) .
\end{equation}

Consider now the second term in (\ref{eqn: asymp_loading_expansion}). By Lemma \ref{lemma:3} and \ref{lemma:limit}, we have
\begin{equation*}
    \bigg\|(\H_j-\H^a)\Q_{j\cdot}\bigg\|^2
    \leq
    \|\H_j-\H^a\|_F^2
    \cdot \|\A_{j\cdot}\|^2
    \cdot \|\Z_k^{-1/2}\|^2
    = O_P\Bigg(
    \min\bigg(\frac{1}{T}, \frac{\eta^2}{(1-\eta)^2}\bigg)
    d_k^{2\alpha_{k,1}-3\alpha_{k,r_k}}
    \Bigg),
\end{equation*}
implying
\begin{align*}
    &\hspace{5mm}
    \bigg\|(\H_j-\H^a)\Q_{j\cdot}\bigg\|^2
    \bigg/ \bigg\|\wh{\D}^{-1}
    \sum_{i=1}^{d_k}\H^a \Q_{i\cdot}
    \sum_{h=1}^{d_{\text{-}k}}
    \frac{1}{|\psi_{ij,h}|}
    \sum_{t\in\psi_{ij,h}}
    E_{t,jh} (\A_{\text{-}k})_{h\cdot}'
    \F_t' \A_{i\cdot}\bigg\|^2 \\
    &= O_P\Bigg(
    \min\bigg(1, \frac{T\eta^2}{(1-\eta)^2}\bigg)
    d_k^{3(\alpha_{k,1}-\alpha_{k,r_k})}
    \prod_{j\in[K]\setminus\{k\} }d_j^{2\alpha_{j,1}-1}
    \Bigg),
\end{align*}
which is unrealistic to be $o_P(1)$ in the presence of missing data in general. Thus $(\H_j-\H^a)\Q_{j\cdot}$ contributes to the asymptotic distribution of $(\wh\Q_{j\cdot}-\H^a\Q_{j\cdot})$. Rewrite
\begin{equation*}
\begin{split}
    (\H_j-\H^a)\Q_{j\cdot} &=
    \wh\D^{-1}\sum_{i=1}^{d_k} \wh\Q_{i\cdot}
    \sum_{h=1}^{\dmk}
    \Big(
    \frac{1}{|\psi_{ij,h}|}
        \sum_{t\in\psi_{ij,h}}
        \bf{Q}_{i\cdot}'
        \bf{F}_{Z,t}
        \bf{\Lambda}_{h\cdot}
        \bf{\Lambda}_{h\cdot}'
        \bf{F}_{Z,t}'
        -
        \frac{1}{T}
        \sum_{t=1}^T
        \bf{Q}_{i\cdot}'
        \bf{F}_{Z,t}
        \bf{\Lambda}_{h\cdot}
        \bf{\Lambda}_{h\cdot}'
        \bf{F}_{Z,t}'
    \Big) \Q_{j\cdot} \\
    &=
    \wh{\bf{D}}^{-1}\sum_{i=1}^{d_k} \wh{\Q}_{i\cdot}\A_{i\cdot}'
    \sum_{h=1}^{\dmk}
    \Delta_{F,k,ij,h}\Z_k^{1/2}\Q_{j\cdot}
    \\
    &=
    \wh\D^{-1}\sum_{i=1}^{d_k} (\wh{\Q}_{i\cdot} - \H^a \Q_{i\cdot}) \A_{i\cdot}'
    \sum_{h=1}^{\dmk}
    \Delta_{F,k,ij,h}\Z_k^{1/2}\Q_{j\cdot}
    + \wh\D^{-1}\H^a \sum_{i=1}^{d_k} \Q_{i\cdot}\A_{i\cdot}' \sum_{h=1}^{\dmk}
    \Delta_{F,k,ij,h}\Z_k^{1/2}\Q_{j\cdot} .
\end{split}
\end{equation*}
Note the first term is dominated by the second term due to Theorem \ref{thm:factor_loading_consistency}. Using Assumption (AD2) and the Slutsky's theorem, we have
\begin{equation}
\label{eqn: asymp_loading_part2}
\sqrt{T d_k^{\alpha_{k,r_k}} }\cdot \wh\D^{-1}\H^a
\sum_{i=1}^{d_k} \Q_{i\cdot}\A_{i\cdot}' \sum_{h=1}^{\dmk}\Delta_{F,k,ij,h}\Z_k^{1/2}\Q_{j\cdot}
\to \c{N}(\0, \D^{-1}\H_k^{a,\ast} h_{k,j}(\A_{j\cdot})(\H_k^{a,\ast})' \D^{-1} )
\;\;\;
\text{$\c{G}^T$-stably}.
\end{equation}
Furthermore, $\c{I}_{H,j}$ and $(\H_j-\H^a)\Q_{j\cdot}$ are asymptotically independent since the randomness of $\c{I}_{H,j}$ comes from $E_{t,jh} (\Amk)_{h\cdot}'\F_t'$ while that of $(\H_j-\H^a)\Q_{j\cdot}$ comes from $\Delta_{F,k,ij,h}$. From (\ref{eqn: asymp_loading_part1}) and (\ref{eqn: asymp_loading_part2}), we conclude that
\[
\sqrt{T d_k^{\alpha_{k,r_k}} }\cdot
(\wh\Q_{j\cdot}-\H^a\Q_{j\cdot})
\xrightarrow{\c{D}}
\cN(\0, \D^{-1}\H_k^{a,\ast}
(T d_k^{\alpha_{k,r_k}} \cdot\bf\Xi_{k,j} + h_{k,j}(\A_{j\cdot}))(\H_k^{a,\ast})'\D^{-1}).
\]
On the other hand, if we have finite missingness or asymptotically vanishing missingness such that
\[
\min\bigg(1, \frac{T\eta^2}{(1-\eta)^2}\bigg)
    d_k^{3(\alpha_{k,1}-\alpha_{k,r_k})}
    \prod_{j\in[K]\setminus\{k\} }d_j^{2\alpha_{j,1}-1} = Td^{-1}g_s^2g_{\eta} d_k^{1+\alpha_{k,1} -3\alpha_{k,r_k}}
    = o(1),
\]
then (\ref{eqn: asymp_loading_part2}) is dominated by (\ref{eqn: asymp_loading_part1}), and hence it holds at the absence of (AD2) that
\[
\sqrt{T\omega_B}\cdot
(\wh\Q_{j\cdot}-\H^a \Q_{j\cdot})
\xrightarrow{\c{D}}
\cN(\0, T\omega_B \cdot
\D^{-1}\H_k^{a,\ast}
\bf\Xi_{k,j}(\H_k^{a,\ast})'\D^{-1}).
\]
This completes the proof of Theorem \ref{thm:asymp_normality_loadings}.
$\square$

\textbf{\textit{Proof of Theorem \ref{thm:covariance_estimator}.}}
By Lemma \ref{lemma:limit}, $\wh\D_k$ is consistent for $\D_k$, and $\H_k^a$ is consistent for $\H_k^{a,\ast}$. Similar to the proof of Theorem 5 in \cite{ChenFan2023}, it suffice to prove that the HAC estimator $\wh\bSigma_{HAC}$ based on $\{\wh\Q_k, \mat{k}{\wh\cC_t}, \mat{k}{\wh\cE_t}\}_{t\in[T]}$ is a consistent estimator for $\H_k^a\bf\Xi_{k,j}(\H_k^a)'$. Recall that
\begin{equation*}
\begin{split}
    \H_k^a\bf\Xi_{k,j}(\H_k^a)'
    &=
    \text{Var}\bigg(
    \sum_{i=1}^{d_k} \H_k^a\Q_{k,i\cdot}
    \sum_{h=1}^{\dmk}
    \frac{1}{|\psi_{k,ij,h}|}
    \sum_{t\in\psi_{k,ij,h}}
        E_{t,jh} (\Amk)_{h\cdot}'
        \F_t' \A_{k,i\cdot}
    \bigg) \\
    &=
    \text{Var}\bigg(
    \sum_{i=1}^{d_k}
    \Big(\wh\D_k^{-1}
    \frac{1}{T}\sum_{t=1}^T
    \wh\Q_k'\Q_k
    \F_{Z,t}\bf{\Lambda}'
    \bf{\Lambda}\bf{F}_{Z,t}'\Big)
    \Q_{k,i\cdot}
    \sum_{h=1}^{\dmk}
    \frac{1}{|\psi_{k,ij,h}|}
    \sum_{t\in\psi_{k,ij,h}}
        E_{t,jh} (\Amk)_{h\cdot}'
        \F_t' \A_{k,i\cdot}
    \bigg) \\
    &=
    \text{Var}\bigg(
    \sum_{i=1}^{d_k}
    \Big(\frac{1}{T}
    \sum_{t=1}^T \wh\D_k^{-1}
    \wh\Q_k' \cdot
    \mat{k}{\cC_t}\cdot
    \mat{k}{\cC_t}_{i\cdot}\Big)\\
    &\;\;\;\cdot
    \sum_{h=1}^{\dmk}
    \frac{1}{|\psi_{k,ij,h}|}
    \sum_{t\in\psi_{k,ij,h}}
    \mat{k}{\cE_t}_{jh}\cdot \mat{k}{\cC_t}_{ih} \bigg) .
\end{split}
\end{equation*}
By Theorem \ref{thm:factor_loading_consistency}, and the rate assumption in the statement of Theorem \ref{thm:covariance_estimator}, we have $\wh\Q_k$ being consistent for a version of $\Q_k$ (in Frobenius norm) for any $k\in[K]$. By Theorem \ref{thm:imputation_consistency} and the assumption that the rate for individual common component imputation error is going to 0, $\wh\cC_{t,i_1,\dots,i_K}$ is consistent for $\cC_{t,i_1,\dots,i_K}$ for any $k\in[K],i_k\in[d_k],t\in[T]$. Hence, it also holds that $\wh\cE_{t,i_1,\dots,i_K}$ is consistent for $\cE_{t,i_1,\dots,i_K}$ for any $k\in[K],i_k\in[d_k],t\in[T]$. We can finally conclude that $\wh\bSigma_{HAC}$ is estimating $\H_k^a\bf\Xi_{k,j}(\H_k^a)'$ consistently \citep{HAC_1987}, which is result 1. We can also show a similar result for $\wh\bSigma_{HAC}^{\Delta}$, which is result 2 (details omitted). Combining both results, and consider the general statement of Theorem \ref{thm:asymp_normality_loadings}, we can easily conclude result 3. This completes the proof of the theorem.
$\square$

We will present two other lemmas before proving Theorem \ref{thm:imputation_consistency}. While we stick with the notations in (\ref{eqn: notation_short}), we use the following also hereafter:
\begin{equation}
\label{eqn: notation_short_factor}
    \begin{split}
    \bf{y}_t &:= \vec{\cY_t}, \;  \bf{m}_t := \vec{\cM_t}, \;
    \bf{f}_{Z,t} := \vec{\cF_{Z,t}}, \;
    \bm\varepsilon_t := \vec{\cE_t},
    \bf{c}_t := \vec{\cC_t}, \;
    \bf{f}_{t} := \vec{\cF_{t}}, \\
    \H_\otimes &:=
    \H_K^a\otimes\dots\otimes \H_1^a, \;
    \A_\otimes :=
    \A_K\otimes\dots\otimes \A_1, \;
    \Z_\otimes :=
    \Z_K\otimes\dots\otimes \Z_1,
    \end{split}
\end{equation}
where the hat versions (if any) of the above are defined similarly.

\begin{lemma}\label{lemma:5_additional}
Under the assumptions in Theorem \ref{thm:imputation_consistency}, for any $k\in[K]$ and $j\in[d_k]$,
\begin{align}
  \|\wh{\bf{Q}}_{k,j\cdot}-\bf{H}_k^a
    \bf{Q}_{k,j\cdot}\|_F^2 &=
    O_P\Big(T^{-1}\dmk
    d_k^{3\alpha_{k,1} - 2\alpha_{k,r_k}} g_s^{-2} + d^2 g_s^{-2}d_k^{2\alpha_{k,1}- 3\alpha_{k,r_k}-2} + g_{\eta} d_k^{2\alpha_{k,1}- 3\alpha_{k,r_k}}\Big). \label{eqn:lemm5rate}
\end{align}
\end{lemma}

\textbf{\textit{Proof of Lemma \ref{lemma:5_additional}.}}
First, consider the case when $T\dmk = o\big(d_k^{ \alpha_{k,r_k} + \alpha_{k,1}}\big)$. From (\ref{eqn: asymp_loading_part1}) in the proof of Theorem 3, we have $\|\wh\Q_{k,j\cdot}-\H_{k,j}\Q_{k,j\cdot}\|_F^2 = O_P(T^{-1}\omega_B^{-1})$. Hence it follows that
\begin{align}
    &\hspace{5mm}
    \|\wh{\bf{Q}}_{k,j\cdot}-\bf{H}_k^a
    \bf{Q}_{k,j\cdot}\|_F^2 = O_P\Big(
    \|\wh\Q_{k,j\cdot}-\H_{k,j}\Q_{k,j\cdot}\|_F^2
    + \|(\H_{k,j}-\H_k^a)\Q_{k,j\cdot} \|_F^2
    \Big) \notag\\
    &=
    O_P\Big( (T\omega_B)^{-1}+ g_{\eta} d_k^{2\alpha_{k,1}- 3\alpha_{k,r_k}}\Big)
    = O_P\Big(T^{-1}\dmk
    d_k^{3\alpha_{k,1} - 2\alpha_{k,r_k}} g_s^{-2} + g_{\eta} d_k^{2\alpha_{k,1}- 3\alpha_{k,r_k}}\Big),
    \label{eqn:lemma5_eq1}
\end{align}
where we used Lemma \ref{lemma:3} in the second equality, and
\begin{align*}
  \norm{\Q_{k,j\cdot}}^2 &= \|\Z_{k}^{-1/2}\A_{k,j\cdot}\|^2
  = O_P( d_k^{-\alpha_{k,r_{k}}}).
\end{align*}

Now suppose $T\dmk = o\big(d_k^{ \alpha_{k,r_k} + \alpha_{k,1}}\big)$ fails to hold. From the decomposition of $\wh\Q_{k,j\cdot}-\H_{k,j}\Q_{k,j\cdot}$ in (\ref{eqn: Q_expansion}), $\c{I}_{H,j}$ is not the leading term anymore, and the leading term among the expressions from (\ref{eqn: I_H_diff_to_prove}) to (\ref{eqn: III_HQ_diff_to_prove}) will be (\ref{eqn: III_H_to_prove}). It has rate
\[
O_P\bigg(\frac{\dmk}{T}+\dmk^2 d_k^{-\alpha_{k,r_k}}\bigg) = O_P(\dmk^2 d_k^{-\alpha_{k,r_k}}),
\]
where the above equality used the fact that $T\dmk = o\big(d_k^{ \alpha_{k,r_k} + \alpha_{k,1}}\big)$ does not hold. Together with the bound on $\|\wh\D_k^{-1}\|_F$ from Lemma \ref{lemma:2} , we have
\begin{equation}\label{eqn:lemma5_eq2}
\|\wh{\bf{Q}}_{k,j\cdot}-\H_{k,j}
\bf{Q}_{k,j\cdot}\|_F^2 = O_P\Big( d^2 g_s^{-2}
d_k^{2\alpha_{k,1}- 3\alpha_{k,r_k}-2}\Big).
\end{equation}
Combining (\ref{eqn:lemma5_eq1}) and (\ref{eqn:lemma5_eq2}), we arrive at the statement of the lemma.
$\square$

\begin{lemma}\label{lemma:5}
Under the Assumptions in Theorem \ref{thm:imputation_consistency}, with the notations in (\ref{eqn: notation_short}) and (\ref{eqn: notation_short_factor}), we have the following for any $j\in[d]$:
\begin{align}
    &\|\Q_{\otimes}\H_\otimes'\|_F^2 = O_P(1),
    \label{eqn: lemma5.2}\\
    &\|\wh\Q_{\otimes, j\cdot} - \H_\otimes\Q_{\otimes, j\cdot}\|^2 \notag\\
    &=
    O_P\Bigg(\max_{k\in[K]}\Bigg\{ T^{-1}\dmk
    d_k^{3\alpha_{k,1} - \alpha_{k,r_k}} g_s^{-2}g_w^{-1} + d^2 g_s^{-2}g_w^{-1}d_k^{2\alpha_{k,1}- 2\alpha_{k,r_k}-2} + g_{\eta}g_w^{-1} d_k^{2\alpha_{k,1}- 2\alpha_{k,r_k}} \Bigg\}\Bigg),
    \label{eqn: lemma5.3}\\
    &\|\wh\Q_{\otimes} - \Q_{\otimes}\H_\otimes'\|_F^2 \notag\\
    &=
    O_P\Bigg(\max_{k\in[K]}\Bigg\{
     T^{-1}d
    d_k^{3\alpha_{k,1} - 2\alpha_{k,r_k}} g_s^{-2} + d^2 g_s^{-2}d_k^{2\alpha_{k,1}- 3\alpha_{k,r_k}-1} + g_{\eta} d_k^{2\alpha_{k,1}- 3\alpha_{k,r_k}+1}
    \Bigg\}\Bigg).
    \label{eqn: lemma5.4}
\end{align}
\end{lemma}

\textbf{\textit{Proof of Lemma \ref{lemma:5}.}}
For (\ref{eqn: lemma5.2}), we have
\begin{equation*}
\begin{split}
    \|\Q_{\otimes}\H_\otimes'\|_F^2 &\leq \norm{\H_{\otimes}}_F^2\cdot \prod_{k=1}^K\norm{\Q_k}_F^2 = O_P(1),
\end{split}
\end{equation*}
where we used Assumption (L1) in the last equality.

To show (\ref{eqn: lemma5.3}), for any $j\in[d_k]$, by a simple induction argument (omitted),
\begin{equation*}
\begin{split}
    &\|\wh\Q_{\otimes, j\cdot}
    -\H_\otimes\Q_{\otimes, j\cdot}\|^2
    =
    \|(\wh\Q_{\otimes}-
    \Q_{\otimes}\H_\otimes')_{j\cdot}
    \|^2
    =
    \Big\|[(\wh\Q_K\otimes\dots\otimes \wh\Q_1) -
    (\Q_K\H_K' \otimes\dots\otimes \Q_1\H_1')]_{j\cdot} \Big\|^2 \\
    &\leq \sum_{k=1}^K\bigg\{\max_{j\in[d_k]} \norm{\wh\Q_{k,j\cdot} - \H_k^a\Q_{k,j\cdot}}^2\prod_{\ell\in[K]\setminus\{k\}}
    \max_{j\in[d_\ell]}\norm{\wh\Q_{\ell,j\cdot}}^2\bigg\}  \\
    &=
    O_P\Bigg(\max_{k\in[K]}\Bigg\{\Bigg[T^{-1}\dmk
    d_k^{3\alpha_{k,1} - 2\alpha_{k,r_k}} g_s^{-2} + d^2 g_s^{-2}d_k^{2\alpha_{k,1}- 3\alpha_{k,r_k}-2} + g_{\eta} d_k^{2\alpha_{k,1}- 3\alpha_{k,r_k}}\Bigg]
    \cdot \prod_{\ell\in[K]\setminus\{k\}}d_\ell^{-\alpha_{\ell,r_{\ell}}}\Bigg\}\Bigg)\\
    &= O_P\Bigg(\max_{k\in[K]}\Bigg\{ T^{-1}\dmk
    d_k^{3\alpha_{k,1} - \alpha_{k,r_k}} g_s^{-2}g_w^{-1} + d^2 g_s^{-2}g_w^{-1}d_k^{2\alpha_{k,1}- 2\alpha_{k,r_k}-2} + g_{\eta}g_w^{-1} d_k^{2\alpha_{k,1}- 2\alpha_{k,r_k}} \Bigg\}\Bigg),
\end{split}
\end{equation*}
where the second last equality used (\ref{eqn:lemm5rate}) and
\begin{align*}
  \norm{\wh\Q_{\ell,j\cdot}}^2 &\leq 2(\norm{\wh\Q_{\ell,j\cdot} - \H_k^a\Q_{\ell,j\cdot}}^2 + \norm{\H_\ell^a\Q_{\ell,j\cdot}}^2) = O_P(\norm{\wh\Q_{\ell,j\cdot} - \H_\ell^a\Q_{\ell,j\cdot}}^2 + \norm{\H_\ell^a\Z_{\ell}^{-1/2}\A_{\ell,j\cdot}}^2)\\
  &= O_P(\norm{\wh\Q_{\ell,j\cdot} - \H_\ell^a\Q_{\ell,j\cdot}}^2 + \norm{\H_\ell^a}_F^2\cdot (d_\ell^{-\alpha_{\ell,r_{\ell}}/2})^2\cdot 1) = O_P(d_\ell^{-\alpha_{\ell,r_{\ell}}}).
\end{align*}
Finally it also holds that
\begin{align*}
     \|\wh\Q_{\otimes}-
   &\Q_{\otimes}\H_\otimes'\|_F^2 =
    \|(\wh\Q_K\otimes\dots\otimes \wh\Q_1) -
    (\Q_K\H_K^{a'} \otimes\dots\otimes \Q_1\H_1^{a'})
    \|_F^2 =
    O(1)\cdot
    \sum_{k=1}^K
    \|\wh\Q_k - \Q_k\H_k^{a'}\|_F^2 \\
    &= O\Bigg( \max_{k\in[K]}\sum_{j=1}^{d_k}\|\wh\Q_{k,j\cdot} - \H_k^a\Q_{k,j\cdot}\|^2 \Bigg)\\
    &=
    O_P\Bigg(\max_{k\in[K]}\Bigg\{
     T^{-1}d
    d_k^{3\alpha_{k,1} - 2\alpha_{k,r_k}} g_s^{-2} + d^2 g_s^{-2}d_k^{2\alpha_{k,1}- 3\alpha_{k,r_k}-1} + g_{\eta} d_k^{2\alpha_{k,1}- 3\alpha_{k,r_k}+1}
    \Bigg\}\Bigg),
\end{align*}
where the second equality could be shown by a simple induction argument using $\norm{\Q_k} = O(1)$ (omitted), and the last equality is from (\ref{eqn:lemm5rate}).
$\square$

\vspace{12pt}

\textbf{\textit{Proof of Theorem \ref{thm:imputation_consistency}.}}
The equation (\ref{eqn: Ft-regression}) is essentially
\begin{align}
    \wh{\bf{f}}_{Z,t} &=
     \bigg(\sum_{j=1}^d m_{t,j}
     \wh\Q_{\otimes, j\cdot}
     \wh\Q_{\otimes, j\cdot}' \bigg)^{-1}
     \bigg(\sum_{j=1}^d m_{t,j}
     \wh\Q_{\otimes, j\cdot}
     y_{t,j} \bigg)\notag\\
     &=
     \bigg(\sum_{j=1}^d m_{t,j}
     \wh\Q_{\otimes, j\cdot}
     \wh\Q_{\otimes, j\cdot}' \bigg)^{-1}
     \bigg(\sum_{j=1}^d m_{t,j}
     \wh\Q_{\otimes, j\cdot}
     (\Q_{\otimes, j\cdot}'\bf{f}_{Z,t}
     + \varepsilon_{t,j} ) \bigg)\notag\\
     &=
     \bigg(\sum_{j=1}^d m_{t,j}
     \wh\Q_{\otimes, j\cdot}
     \wh\Q_{\otimes, j\cdot}' \bigg)^{-1}
     \bigg(\sum_{j=1}^d m_{t,j}
     \wh\Q_{\otimes, j\cdot}
     \Q_{\otimes, j\cdot}'\bf{f}_{Z,t}
     \bigg)
     +
     \bigg(\sum_{j=1}^d m_{t,j}
     \wh\Q_{\otimes, j\cdot}
     \wh\Q_{\otimes, j\cdot}' \bigg)^{-1}
     \bigg(\sum_{j=1}^d m_{t,j}
     \wh\Q_{\otimes, j\cdot}
     \varepsilon_{t,j} \bigg) \notag\\
     &=
     \bigg(\sum_{j=1}^d m_{t,j}
     \wh\Q_{\otimes, j\cdot}
     \wh\Q_{\otimes, j\cdot}' \bigg)^{-1}
     \bigg(\sum_{j=1}^d m_{t,j}
     \wh\Q_{\otimes, j\cdot}
     \wh\Q_{\otimes, j\cdot}'
     \bigg)(\H_\otimes')^{-1}
     \bf{f}_{Z,t} \notag\\
     &\;\;\;+
     \bigg(\sum_{j=1}^d m_{t,j}
     \wh\Q_{\otimes, j\cdot}
     \wh\Q_{\otimes, j\cdot}' \bigg)^{-1}
     \bigg(\sum_{j=1}^d m_{t,j}
     \wh\Q_{\otimes, j\cdot}
     \varepsilon_{t,j} \bigg) \notag\\
     &\;\;\;+
     \bigg(\sum_{j=1}^d m_{t,j}
     \wh\Q_{\otimes, j\cdot}
     \wh\Q_{\otimes, j\cdot}' \bigg)^{-1}
     \bigg(\sum_{j=1}^d m_{t,j}
     \wh\Q_{\otimes, j\cdot}
     (\H_\otimes\Q_{\otimes, j\cdot}-
     \wh\Q_{\otimes, j\cdot})'
     \bigg)(\H_\otimes')^{-1}
     \bf{f}_{Z,t} \notag\\
     &=:
     (\H_\otimes')^{-1}\bf{f}_{Z,t}
     + \wt{\bm\varepsilon}_{H,t}
     + \wt{\bm\varepsilon}_t
     + \wt{\bf{f}}_{Z,t}
     ,\;\text{ where } \label{eqn: F_expansion}
\end{align}
\begin{align*}
    \wt{\bm\varepsilon}_{H,t} &:=
    \bigg(\sum_{j=1}^d m_{t,j}
     \wh\Q_{\otimes, j\cdot}
     \wh\Q_{\otimes, j\cdot}' \bigg)^{-1}
     \bigg(\sum_{j=1}^d m_{t,j}
     \H_\otimes \Q_{\otimes, j\cdot}
     \varepsilon_{t,j} \bigg)
    , \\
    \wt{\bm\varepsilon}_t &:=
    \bigg(\sum_{j=1}^d m_{t,j}
     \wh\Q_{\otimes, j\cdot}
     \wh\Q_{\otimes, j\cdot}' \bigg)^{-1}
     \bigg(\sum_{j=1}^d m_{t,j}
     (\wh\Q_{\otimes, j\cdot}
     -\H_\otimes\Q_{\otimes, j\cdot})
     \varepsilon_{t,j} \bigg)
    , \\
    \wt{\bf{f}}_{Z,t} &:=
    \bigg(\sum_{j=1}^d m_{t,j}
     \wh\Q_{\otimes, j\cdot}
     \wh\Q_{\otimes, j\cdot}' \bigg)^{-1}
     \bigg(\sum_{j=1}^d m_{t,j}
     \wh\Q_{\otimes, j\cdot}
     (\H_\otimes\Q_{\otimes, j\cdot}-
     \wh\Q_{\otimes, j\cdot})'
     \bigg)(\H_\otimes')^{-1}
     \bf{f}_{Z,t} .
\end{align*}
Then we have
\begin{equation*}
\begin{split}
    & \hspace{5mm}
    \bigg\| \sum_{j=1}^d m_{t,j}
    \wh\Q_{\otimes, j\cdot}
    \wh\Q_{\otimes, j\cdot}'
    - \sum_{j=1}^d m_{t,j}
    \H_\otimes\Q_{\otimes, j\cdot}
    \Q_{\otimes, j\cdot}'\H_\otimes'
    \bigg\|_F \\
    &\leq
    \sum_{j=1}^d \bigg\|
    \wh\Q_{\otimes, j\cdot}
    \wh\Q_{\otimes, j\cdot}'
    - \H_\otimes\Q_{\otimes, j\cdot}
    \Q_{\otimes, j\cdot}'\H_\otimes'
    \bigg\|_F \\
    &\leq
    \sum_{j=1}^d
    \|\wh\Q_{\otimes, j\cdot}-\H_\otimes
    \Q_{\otimes, j\cdot}\|^2 +
    2 \sum_{j=1}^d
    \|\wh\Q_{\otimes, j\cdot}-\H_\otimes
    \Q_{\otimes, j\cdot}\| \cdot
    \|\H_\otimes\Q_{\otimes, j\cdot}\| \\
    &=O_P(\norm{\wh\Q_{\otimes} - \Q_{\otimes}\H_{\otimes}'}_F^2 + \norm{\wh\Q_{\otimes} - \Q_{\otimes}\H_{\otimes}'}_F)\\
    &=O_P\Bigg(\max_{k\in[K]}\Bigg\{
     T^{-1}d
    d_k^{3\alpha_{k,1} - 2\alpha_{k,r_k}} g_s^{-2} + d^2 g_s^{-2}d_k^{2\alpha_{k,1}- 3\alpha_{k,r_k}-1} + g_{\eta} d_k^{2\alpha_{k,1}- 3\alpha_{k,r_k}+1}
    \Bigg\}^{1/2}\Bigg) = o_P(1),
\end{split}
\end{equation*}
where we used the Cauchy--Schwarz inequality, (\ref{eqn: lemma5.2}) and (\ref{eqn: lemma5.4}) in the last equality, and Assumption (R1). Hence we have
\begin{equation}
\begin{split}
\label{eqn: QQ_hat_otimes_lim}
    \sum_{j=1}^d m_{t,j}
    \wh\Q_{\otimes, j\cdot}
    \wh\Q_{\otimes, j\cdot}'
    &\xrightarrow{p}
    \sum_{j=1}^d m_{t,j}
    \H_\otimes\Q_{\otimes, j\cdot}
    \Q_{\otimes, j\cdot}'\H_\otimes'.
\end{split}
\end{equation}
Note that by (\ref{eqn: H_result}) we have $\|\H_\otimes^{-1}\|_F=O_P(1)\cdot \|\Q_K'\wh\Q_K\otimes\dots\otimes\Q_1'\wh\Q_1\|=O_P(1)$, which will be used later in the proof. To bound $\bf{f}_{Z,t}$, from Assumption (F1), we have
\begin{equation*}
\begin{split}
    \b{E}\|\bf{f}_{t}\|^2
    &= r
    = O(1),
\end{split}
\end{equation*}
and hence with Assumption (L1),
\begin{equation}
\begin{split}
\label{eqn: f_Zt_bound}
    \|\bf{f}_{Z,t}\|^2 &\leq
    \|\Z_\otimes^{1/2}\|_F^2
    \cdot \|\bf{f}_{t}\|^2
    = O_P\Bigg(
    \prod_{j=1}^Kd_j^{\alpha_{j,1}}\Bigg) = O_P(g_s).
\end{split}
\end{equation}

With (\ref{eqn: lemma5.4}), (\ref{eqn: QQ_hat_otimes_lim}) and the above result,
\begin{equation}
\begin{split}
\label{eqn: tilde_f_Zt_bound}
    \|\wt{\bf{f}}_{Z,t}\|^2
    &= O_P(1)\cdot
    \bigg(\sum_{j=1}^d
    \|\wh\Q_{\otimes, j\cdot}\|^2
    \bigg)\bigg(\sum_{j=1}^d
    \|\H_\otimes\Q_{\otimes, j\cdot}-
    \wh\Q_{\otimes, j\cdot}\|^2
    \bigg) \\
    &\cdot
    \|(\bSigma_{A,K}\otimes
    \dots\otimes\bSigma_{A,1})^{-1}\|_F^2
    \cdot\|\H_\otimes^{-1}\|_F^6
    \cdot \|\bf{f}_{Z,t}\|^2 \\
    &=
    O_P\Bigg(\max_{k\in[K]}\Bigg\{
     T^{-1}d
    d_k^{3\alpha_{k,1} - 2\alpha_{k,r_k}} g_s^{-1} + d^2 g_s^{-1}d_k^{2\alpha_{k,1}- 3\alpha_{k,r_k}-1} + g_{\eta}g_s d_k^{2\alpha_{k,1}- 3\alpha_{k,r_k}+1}
    \Bigg\}\Bigg),
\end{split}
\end{equation}
where we also used Assumption (L1) in the last equality. Similarly, by (\ref{eqn: lemma5.3}),
\begin{equation}
\begin{split}
\label{eqn: tilde_varepsilon_t_bound}
    \|\wt{\bm\varepsilon}_{t}\|^2
    &=
    O_P(1)\cdot
    \|(\bSigma_{A,K}\otimes
    \dots\otimes\bSigma_{A,1})^{-1}\|_F^2
    \cdot\|\H_\otimes^{-1}\|_F^4\cdot
    O_P\big(\max_{j\in[d]}
    \|\H_\otimes\Q_{\otimes, j\cdot}-
    \wh\Q_{\otimes, j\cdot}\|^2
    \big)\bigg(
    \sum_{j,\ell=1}^d
    |\b{E}\varepsilon_{t,j}\varepsilon_{t,\ell}|
    \bigg)\\
    &=
        O_P\Bigg(\max_{k\in[K]}\Bigg\{ T^{-1}d\dmk
    d_k^{3\alpha_{k,1} - \alpha_{k,r_k}} g_s^{-2}g_w^{-1} + d^3 g_s^{-2}g_w^{-1}d_k^{2\alpha_{k,1}- 2\alpha_{k,r_k}-2} + dg_{\eta}g_w^{-1} d_k^{2\alpha_{k,1}- 2\alpha_{k,r_k}} \Bigg\}\Bigg),
\end{split}
\end{equation}
since $\sum_{j,\ell=1}^d|\b{E}\varepsilon_{t,j}\varepsilon_{t,\ell}| = O(d)$ by Assumption (E1). By the same token,
\begin{equation}
\begin{split}
\label{eqn: tilde_varepsilon_Ht_bound}
    \|\wt{\bm\varepsilon}_{H,t}\|^2
    &=
    O_P(1)\cdot
    \|(\bSigma_{A,K}\otimes
    \dots\otimes\bSigma_{A,1})^{-1}\|_F^2
    \cdot \|\H_\otimes^{-1}\|_F^2 \cdot \|\Z_{\otimes}^{-1/2}\|_F^2
    \cdot \bigg\|
    \sum_{j=1}^d m_{t,j}
    \A_{\otimes, j\cdot}
    \varepsilon_{t,j}
    \bigg\|^2
    =
    O_P(d/g_w),
\end{split}
\end{equation}
where we used
\begin{align*}
  \b{E}\bigg\|\sum_{j=1}^dm_{t,j}\A_{\otimes,j\cdot}\epsilon_{t,j}\bigg\|^2 \leq \max_{j\in[d]} \norm{\A_{\otimes,j\cdot}}_F^2\sum_{j,\ell=1}^d|\b{E}\epsilon_{t,j}\epsilon_{t,\ell}| = O(d).
\end{align*}
Therefore, we have from (\ref{eqn: tilde_f_Zt_bound}), (\ref{eqn: tilde_varepsilon_Ht_bound}) and (\ref{eqn: tilde_varepsilon_t_bound}),
\begin{align*}
    &\|\wh{\bf{f}}_{Z,t}-
    (\H_\otimes')^{-1}\bf{f}_{Z,t}\|^2
    \leq
    \|\wt{\bm\varepsilon}_{H,t}\|^2
    + \|\wt{\bm\varepsilon}_t\|^2
    + \|\wt{\bf{f}}_{Z,t}\|^2 \\
    &=
    O_P\Bigg(\max_{k\in[K]}\Bigg\{
     T^{-1}d
    d_k^{3\alpha_{k,1} - 2\alpha_{k,r_k}} g_s^{-1} + d^2 g_s^{-1}d_k^{2\alpha_{k,1}- 3\alpha_{k,r_k}-1} + g_{\eta}g_s d_k^{2\alpha_{k,1}- 3\alpha_{k,r_k}+1}
    \Bigg\} + \frac{d}{g_w}\Bigg),
\end{align*}
where we also used $1/2<\alpha_{k,r_k}\leq \alpha_{k,1} \leq 1$ from Assumption (L1) to conclude that $dg_s^{-1}g_w^{-1} = o(1)$, so that in fact $\|\wt{\bm\varepsilon}_t\|^2 = o_P(\|\wt{\bf{f}}_{Z,t}\|^2)$.

Now from (\ref{eqn: F_expansion}) and using the notations in (\ref{eqn: notation_short_factor}), we can obtain the vectorized imputed values, which are the vectorized estimated common components, as $\wh{\bf{c}}_{t}=\wh\Q_{\otimes}\wh{\bf{f}}_{Z,t}$ for any $t\in[T]$. Then for $j\in[d]$, we have the squared individual imputation error as
\begin{equation*}
\begin{split}
    &(\wh{\c{C}}_{t,i_1,\dots,i_K}-
    \c{C}_{t,i_1,\dots,i_K})^2
    =
    (\wh{\bf{c}}_{t} - \bf{c}_{t})_j^2
    =
    (\wh\Q_{\otimes,j\cdot}'
    \wh{\bf{f}}_{Z,t} -
    \Q_{\otimes,j\cdot}'\bf{f}_{Z,t})^2 \\
    &=
    \bigg[
    \bigg(\wh\Q_{\otimes,j\cdot}
    - \H_\otimes
    \Q_{\otimes,j\cdot}\bigg)'
    \bigg((\H_\otimes')^{-1}\bf{f}_{Z,t}
     + \wt{\bm\varepsilon}_{H,t}
     + \wt{\bm\varepsilon}_t
     + \wt{\bf{f}}_{Z,t}\bigg)
     +
     \A_{\otimes,j\cdot}'
     \Z_\otimes^{-1/2}\H_\otimes'
     \bigg(\wt{\bm\varepsilon}_{H,t}
     + \wt{\bm\varepsilon}_t
     + \wt{\bf{f}}_{Z,t}\bigg)
     \bigg]^2 \\
    &=
    O_P\Bigg(\max_{k\in[K]}\Bigg\{
     T^{-1}d
    d_k^{3\alpha_{k,1} - 2\alpha_{k,r_k}} g_s^{-1}g_w^{-1} + d^2 g_s^{-1}g_w^{-1}d_k^{2\alpha_{k,1}- 3\alpha_{k,r_k}-1} + g_{\eta}g_sg_w^{-1} d_k^{2\alpha_{k,1}- 3\alpha_{k,r_k}+1}
    \Bigg\} + \frac{d}{g_w^2}\Bigg),
\end{split}
\end{equation*}
where we used (\ref{eqn: lemma5.3}), (\ref{eqn: lemma5.4}) and Assumption (R1) in the last equality.

Lastly, we have the average imputation error as the following,
\begin{equation}
\begin{split}
    &\hspace{5mm}
    \frac{1}{Td} \sum_{t=1}^T
    \sum_{i_1,\dots,i_K=1}^{d_1,\dots,d_K}
    (\wh{\c{C}}_{t,i_1,\dots,i_K}-
    \c{C}_{t,i_1,\dots,i_K})^2 \\
    &=
    \frac{1}{Td}
    \sum_{t=1}^T\|\wh{\bf{c}}_{t} - \bf{c}_{t}\|^2
    =
    \frac{1}{Td}
    \sum_{t=1}^T\|\wh\Q_{\otimes}
    \wh{\bf{f}}_{Z,t} -
    \Q_{\otimes}\bf{f}_{Z,t}\|^2 \\
    &=
    \frac{1}{Td}
    \sum_{t=1}^T\bigg\|
    \big(\wh\Q_{\otimes}
    -\Q_{\otimes}\H_\otimes'\big)
    (\H_\otimes')^{-1}\bf{f}_{Z,t}
     +
     \wh\Q_{\otimes}
     \bigg(\wt{\bm\varepsilon}_{H,t}
     + \wt{\bm\varepsilon}_t
     + \wt{\bf{f}}_{Z,t}\bigg)
     \bigg\|^2 \\
    &=
   O_P\Bigg(\max_{k\in[K]}\Bigg\{
     T^{-1}
    d_k^{3\alpha_{k,1} - 2\alpha_{k,r_k}} g_s^{-1} + dg_s^{-1}d_k^{2\alpha_{k,1}- 3\alpha_{k,r_k}-1} + d^{-1}g_{\eta}g_s d_k^{2\alpha_{k,1}- 3\alpha_{k,r_k}+1}
    \Bigg\} + \frac{1}{g_w}\Bigg),
\end{split}
\end{equation}
where we used (\ref{eqn: tilde_f_Zt_bound}), (\ref{eqn: tilde_varepsilon_Ht_bound}), (\ref{eqn: tilde_varepsilon_t_bound}) and Lemma \ref{lemma:5} in the last equality, and the fact that $\|\wh\Q_{\otimes}\|_F^2 = r = O(1)$. This completes the proof of Theorem \ref{thm:imputation_consistency}.
$\square$

\textbf{\textit{Proof of Theorem \ref{thm:number_of_factors}.}} Firstly, we use the notations in (\ref{eqn: notation_short}), and define also $\Z := \Z_k$ and $\R^* := \R_k^*$, which coincides with the $\R^*$ defined in (\ref{eqn:S_decomposition}). Then for $j\in[r_k]$,
\begin{align}
   \lambda_{j}(\R^*) &= \lambda_{j}\Bigg(\frac{1}{T}\sum_{t=1}^T\Q \F_{Z,t} \bf{\Lambda}'\bf{\Lambda} \F_{Z,t}' \Q'\Bigg) = \lambda_{j}\Bigg(\frac{1}{T}\sum_{t=1}^T\A_k\F_{t}
\big[\otimes_{j\in[K]\setminus\{k\}}\A_j\big]'\big[\otimes_{j\in[K]\setminus\{k\}}\A_j\big]\F_{t}'\A_k'\Bigg)\notag\\
    &= \lambda_j\Bigg(\A_k'\A_k \cdot \frac{1}{T}\sum_{t=1}^T\F_t[\otimes_{\ell\in[K]\setminus\{k\}}\A_\ell]'[\otimes_{\ell\in[K]\setminus\{k\}}
\A_\ell]\F_t'\Bigg)\notag\\
&\asymp_P \lambda_j(\A_k'\A_k\cdot \tr(\otimes_{\ell\in[K]\setminus\{k\}}\A_\ell'\A_\ell)) = \lambda_j(\A_k'\A_k)\prod_{\ell\in[K]\setminus\{k\}}\tr(\A_\ell'\A_\ell)\notag\\
&\asymp \lambda_j(\Z\Q'\Q)\prod_{\ell\in[K]\setminus\{k\}}\sum_{i=1}^{r_k}d_\ell^{\alpha_{\ell,i}}
\asymp \lambda_j(\bSigma_{A,k}^{1/2}\Z\bSigma_{A,k}^{1/2})\cdot \prod_{\ell\in[K]\setminus\{k\}}d_\ell^{\alpha_{\ell,1}}\notag\\
&\asymp \lambda_j(\Z)d_k^{-\alpha_{k,1}}g_s \asymp g_sd_k^{\alpha_{k,j}-\alpha_{k,1}}, \label{eqn:lambda_j(R*)}
\end{align}
where the third line uses Assumption (F1), and Assumption (L1) in the second last line. The last line uses Theorem 1 of \cite{Ostrowsk1959} on the eigenvalues of a congruent transformation $\bSigma_{A,k}^{1/2}\Z\bSigma_{A,k}^{1/2}$ of $\Z$, and from Assumption (L1) that $\bSigma_{A,k}$ has eigenvalues uniformly bounded away from 0 and infinity.

Since $\wh\S = \R^* + (\wh\S - \R^*)$, for $j\in[r_k]$, we have by Weyl's inequality that
\begin{align}
  |\lambda_j(\wh\S) - \lambda_j(\R^*)| &\leq \norm{\wh\S - \R^*}
\leq \norm{\wt\R - \R^*} + \norm{\R_1} + \norm{\R_2} + \norm{\R_3} \notag\\
&\leq \omega_k\bigg\{\sup_{\|\bgamma\|=1}|\wt R(\bgamma) - R^*(\bgamma)| + \sup_{\|\bgamma\|=1}R_1
+\sup_{\|\bgamma\|=1}R_2 +\sup_{\|\bgamma\|=1}R_3\bigg\} = o_P(\omega_k), \label{eqn:difference_in_evalues}
\end{align}
where we use the decomposition in (\ref{eqn:S_decomposition}) in the first line, and $\omega_k := g_sd_k^{\alpha_{k,r_k}-\alpha_{k,1}}$ is defined at the beginning of the proof of Lemma \ref{lemma:2}.
The second line uses $\wt R(\bgamma), R^*(\bgamma), R_1, R_2$ and $R_3$ defined in (\ref{eqn:different_R}), and the convergence in probability in (\ref{eqn:RminusRtildeto0}) and (\ref{eqn:RtildeminusR*to0}).

Secondly, with Assumption (R1) and our choice of $\xi$ (see also (\ref{eqn:lem2_withoutE1})),
\begin{align}
  \xi/\omega_k \asymp dg_s^{-1}d_k^{\alpha_{k,1} - \alpha_{k,r_k}}[(T\dmk)^{-1/2} + d_k^{-1/2}] = o(1).
\label{eqn:xi_over_omegak}
\end{align}

For $r_k > 1$, if $j\in[r_k-1]$, using (\ref{eqn:difference_in_evalues}) and (\ref{eqn:xi_over_omegak}), consider
\begin{align}
  \frac{\lambda_{j+1}(\wh\S) + \xi}{\lambda_j(\wh\S) + \xi} &\leq \frac{\lambda_{j+1}(\R^*) + \xi + |\lambda_{j+1}(\wh\S) - \lambda_{j+1}(\R^*)|}{\lambda_j(\R^*) + \xi - |\lambda_j(\wh\S)-\lambda_j(\R^*)|}
= \frac{\lambda_{j+1}(\R^*) + o_P(\omega_k)}{\lambda_j(\R^*) + o_P(\omega_k)}\notag\\
&= \frac{\lambda_{j+1}(\R^*)}{\lambda_j(\R^*)}(1+o_P(1)) \asymp_P d_k^{\alpha_{k,j+1} - \alpha_{k,j}},
\label{eqn:perturbed_ratio1}
\end{align}
where the last line uses (\ref{eqn:lambda_j(R*)}). Also, for $j\in[r_k-1]$,
\begin{align}
  \frac{\lambda_{r_k+1}(\wh\S) + \xi}{\lambda_{r_k}(\wh\S) + \xi} &= \frac{\lambda_{r_k+1}(\wh\S) + \xi}{\omega_k(1+o_P(1))} = O_P\bigg(\frac{\lambda_{r_k+1}(\wh\S)}{\omega_k} + \frac{\xi}{\omega_k}\bigg)\\
&= O_P\Big(\sup_{\|\bgamma\|=1}(\wt{R}(\bgamma) - R^*(\bgamma) + R_1+R_2+R_3) + \xi/\omega_k\Big)\notag\\
&=O_P(\xi/\omega_k) = o_P(d_k^{\alpha_{k,j+1} - \alpha_{k,j}}), \label{eqn:perturbed_ratio2}
\end{align}
where the second last equality uses (\ref{eqn:lem2_withoutE1}), (\ref{eqn:R2rate}) and (\ref{eqn:RtildeminusR*to0}) together with our choice of $\xi$,
and the last equality uses the extra rate assumption in the statement of the theorem.
In the third equality, we assume the following is true (to be shown at the end of this proof):
\begin{equation}\label{eqn:r_k+1_evalue_of_Shat}
\lambda_{j}(\wh\S) =  \lambda_j((\wt\R - \R^*) + \R_1 + \R_2 + \R_3), \;\;\; j = r_k+1,\ldots,d_k,
\end{equation}
so that
\begin{align*}
\frac{\lambda_j(\wh\S)}{\omega_k} =\lambda_j\bigg(\frac{1}{\omega_k}((\wt\R - \R^*) + \R_1 + \R_2 + \R_3)\bigg)  \leq \sup_{\|\bgamma\|=1}((\wt{R}(\bgamma) - R^*(\bgamma)) + R_1+R_2+R_3).
\end{align*}
Hence for $j=r_k+1,\ldots,\lfloor d_k/2 \rfloor$ (true also for $r_k=1$),
\begin{align}
  \frac{\lambda_{j+1}(\wh\S) + \xi}{\lambda_j(\wh\S) + \xi} \geq \frac{\xi/\omega_k}{\sup_{\|\bgamma\|=1}((\wt{R}(\bgamma) - R^*(\bgamma))+R_1+R_2+R_3) + \xi/\omega_k} \geq \frac{1}{C} \label{eqn:perturbed_ratio3}
\end{align}
in probability for some generic constant $C>0$, where the last inequality uses (\ref{eqn:lem2_withoutE1}), (\ref{eqn:R2rate}) and (\ref{eqn:RtildeminusR*to0}) together with our choice of $\xi$.
Combining (\ref{eqn:perturbed_ratio1}), (\ref{eqn:perturbed_ratio2}) and (\ref{eqn:perturbed_ratio3}), we can easily see that our proposed $\wh{r}_k$ is a consistent estimator for $r_k$.

If $r_k=1$, then (\ref{eqn:perturbed_ratio2}) becomes
\[\frac{\lambda_{r_k+1}(\wh\S) + \xi}{\lambda_{r_k}(\wh\S) + \xi} = O_P(\xi/\omega_k) = o_P(1).\]
When combined with (\ref{eqn:perturbed_ratio3}) which is true also for $r_k=1$, we can see that $\wh{r}_k = 1$ in probability, showing that $\wh{r}_k$ is a consistent estimator of $r_k$.

It remains to show (\ref{eqn:r_k+1_evalue_of_Shat}). To this end, from (\ref{eqn:lambda_j(R*)}) and (\ref{eqn:difference_in_evalues}), the first $r_k$ eigenvalues of $\wh\S$ coincides with those of $\R^*$ asymptotically, so that the first $r_k$ eigenvectors corresponding to $\wh\S$ coincides with those for $\R^*$ asymptotically as $T,d_k\rightarrow \infty$, which are necessarily in $\mathcal{N}^\perp:= \text{Span}(\Q)$, the linear span of the columns of $\Q$ (see (\ref{eqn:R*}), where $\R^*$ is sandwiched by $\Q$ and $\Q'$). This means that the $(r_k+1)$-th largest eigenvalue of $\wh\S$ and beyond will asymptotically have eigenvectors in $\mathcal{N}$, the orthogonal complement of $\mathcal{N}^\perp$. Then for any unit vectors $\bgamma \in \mathcal{N}$, we have from the definitions of $\R^*$, $\wt\R$. $\R_1$, $\R_2$ and $\R_3$ in (\ref{eqn:different_R}) that
\begin{align*}
  \bgamma'\wh\S\bgamma = \bgamma'(\R^* + (\wt\R-\R^*)+\R_1+\R_2+\R_3)\bgamma = \bgamma'((\wt\R-\R^*)+\R_1+\R_2+\R_3)\bgamma,
\end{align*}
which is equivalent to (\ref{eqn:r_k+1_evalue_of_Shat}). This completes the proof of the theorem. $\square$

\textbf{\textit{Proof of Proposition \ref{Prop:example_AD2}.}}
We can show stable convergence in law similar to Proposition 3.1 in \cite{Xiong_Pelger}. First, using Assumption (F1) we can write
\begin{equation*}
\begin{split}
    &\hspace{5mm}
    \sqrt{T d_k^{\alpha_{k,r_k}} }\cdot \D^{-1}\H_k^{a,\ast}
    \sum_{i=1}^{d_k} \Q_{k,i\cdot} \A_{k,i\cdot}' \sum_{h=1}^{\dmk}
    \Delta_{F,k,ij,h} \A_{k,j\cdot} \\
    &=
    \sum_{t=1}^T \sqrt{\frac{{d_k^{\alpha_{k,r_k}}}}{T}} \cdot \sum_{i=1}^{d_k}  \sum_{h=1}^{\dmk}
    \Bigg(\frac{T \cdot \b{1}\{t\in\psi_{k,ij,h}\} }{|\psi_{k,ij,h}|} -1 \Bigg) \\
    &\cdot
    \D^{-1}\H_k^{a,\ast}
    \Q_{k,i\cdot} \A_{k,i\cdot}'
    (\textnormal{mat}_k(\cF_t)\v_{k,h}
    \v_{k,h}'\textnormal{mat}_k(\cF_t)'
    - \v_{k,h}'\v_{k,h}\bSigma_k)
    \A_{k,j\cdot} .
\end{split}
\end{equation*}
Define the filtration $\c{G}^T:=\sigma(\cup_{s=1}^T\c{G}_s)$ where the sigma-algebra $\c{G}_s:=\sigma(\{\c{M}_{t,i_1, \dots,i_K}\mid t\leq s\}, \A_1, \dots,\A_K)$. Let $\bf{u}\in\b{R}^{r_k}$ be a non-random unit vector. For a given $k\in[K], j\in[d_k]$, define also the random variable
\begin{equation*}
\begin{split}
g_{k,j,t} &:= \bf{u}'\sqrt{\frac{{d_k^{\alpha_{k,r_k}}}}{T}} \cdot
    \sum_{i=1}^{d_k}  \sum_{h=1}^{\dmk}
    \Bigg(\frac{T \cdot \b{1}\{t\in\psi_{k,ij,h}\} }{|\psi_{k,ij,h}|} -1 \Bigg) \\
    &\cdot
    \D^{-1}\H_k^{a,\ast}
    \Q_{k,i\cdot} \A_{k,i\cdot}'
    (\textnormal{mat}_k(\cF_t)\v_{k,h}
    \v_{k,h}'\textnormal{mat}_k(\cF_t)'
    - \v_{k,h}'\v_{k,h}\bSigma_k)
    \A_{k,j\cdot} .
\end{split}
\end{equation*}
Since each entry in $\cF_t$ is i.i.d. by Assumption (F1) and is independent of $(\cM_t,\A_1,\dots,\A_K)$ by Assumptions (O1) and (L1), we have $\b{E}[g_{k,j,t}\mid \c{G}_{t-1}]= 0$. Define $\bf\Xi_{F,k}:= \text{var}\Big[\Vec{\textnormal{mat}_k(\cF_t)\Amk'\Amk\textnormal{mat}_k(\cF_t)' - \tr(\Amk\Amk')\bSigma_k}\Big]$ and $\bf{x}_{F,k,j,il}:= \Vec{[\A_{k,j\cdot}'\otimes (\D^{-1}\H_k^{a,\ast} \Q_{k,i\cdot} \A_{k,i\cdot}')] \bf\Xi_{F,k} [\A_{k,j\cdot}'\otimes (\D^{-1}\H_k^{a,\ast} \Q_{k,l\cdot} \A_{k,l\cdot}')]' }$, so that we have
\begin{equation*}
    \|\bf{x}_{F,k,j,il}\|^2 \leq
    \|\bf\Xi_{F,k}\|_F^2 \cdot \|\D^{-1}\|_F^4
    \cdot \|\Z_k^{-1/2}\|_F^4
    =
    O_P\Big(
    d_k^{4\alpha_{k,1}-6\alpha_{k,r_k}}
    \prod_{j=1}^Kd_j^{-4\alpha_{j,1}}
    \prod_{j\in[K]\setminus\{k\} }d_j^{4\alpha_{j,1}}
    \Big)
    = O_P(d_k^{-6\alpha_{k,r_k}}) ,
\end{equation*}
leading to
\begin{align*}
    &d_k^{\alpha_{k,r_k}} \cdot\b{E}\Bigg\|
    \sum_{i=1}^{d_k}\sum_{l=1}^{d_k}
    \Bigg(\frac{T \cdot \b{1}\{t\in\psi_{k,ij}\} }{|\psi_{k,ij}|} -1 \Bigg)
    \Bigg(\frac{T \cdot \b{1}\{t\in\psi_{k,lj}\} }{|\psi_{k,lj}|} -1 \Bigg)
    (\bf{x}_{F,k,j,il} - \b{E}[\bf{x}_{F,k,j,il}])
    \Bigg\|^2\\
    &= O_P(d_k^{4-5\alpha_{k,r_k}}) = o_P(1).
\end{align*}
Hence, it holds that
\begin{equation*}
\begin{split}
    \sum_{t=1}^T\b{E}[g_{k,j,t}^2\mid \c{G}_{t-1}] &=
    \frac{d_k^{\alpha_{k,r_k}}}{T} \cdot
    \sum_{t=1}^T \b{E}\Bigg\{\bf{u}'
    \sum_{i=1}^{d_k}\sum_{l=1}^{d_k}
    \Bigg(\frac{T \cdot \b{1}\{t\in\psi_{k,ij}\} }{|\psi_{k,ij}|} -1 \Bigg)
    \Bigg(\frac{T \cdot \b{1}\{t\in\psi_{k,lj}\} }{|\psi_{k,lj}|} -1 \Bigg) \\
    &\cdot
    \D^{-1}\H_k^{a,\ast}[\A_{k,j\cdot}'\otimes (
    \Q_{k,i\cdot} \A_{k,i\cdot}')]
    \Vec{\textnormal{mat}_k(\cF_t)\Amk'
    \Amk\textnormal{mat}_k(\cF_t)'
    - \tr(\Amk\Amk')\bSigma_k} \\
    & \cdot \Vec{\textnormal{mat}_k(\cF_t)\Amk'
    \Amk\textnormal{mat}_k(\cF_t)'
    - \tr(\Amk\Amk')\bSigma_k}' \\
    &\cdot
    [\A_{k,j\cdot}'\otimes (
    \Q_{k,l\cdot} \A_{k,l\cdot}')]'
    (\H_k^{a,\ast})'\D^{-1} \bf{u}
    \mid \c{G}_{t-1}\Bigg\} \\
    & \xrightarrow{p}
    d_k^{2+\alpha_{k,r_k}}\cdot
    \lim_{d_k\to\infty}\frac{1}{d_k^2}
    \sum_{i=1}^{d_k}\sum_{l=1}^{d_k}
    \Bigg(\frac{T \cdot \b{1}\{t\in\psi_{k,ij}\} }{|\psi_{k,ij}|} -1 \Bigg)
    \Bigg(\frac{T \cdot \b{1}\{t\in\psi_{k,lj}\} }{|\psi_{k,lj}|} -1 \Bigg) \\
    &\cdot
    \bf{u}'\D^{-1}\H_k^{a,\ast}[\A_{k,j\cdot}'\otimes (
    \Q_k' \A_k)] \bf\Xi_{F,k}
    [\A_{k,j\cdot}'\otimes (\Q_k' \A_k)]'
    (\H_k^{a,\ast})'\D^{-1} \bf{u} \\
    & \xrightarrow{p}
    d_k^{2+\alpha_{k,r_k}} \omega_{\psi,k,j} \cdot
    \bf{u}'\D^{-1}\H_k^{a,\ast}[\A_{k,j\cdot}'\otimes (
    \Q_k' \A_k)] \bf\Xi_{F,k}
    [\A_{k,j\cdot}'\otimes (\Q_k' \A_k)]'
    (\H_k^{a,\ast})'\D^{-1}\bf{u} ,
\end{split}
\end{equation*}
which satisfies the nesting condition of Theorem 6.1 in \cite{HauslerLuschgy2015}. From Assumption (O1), we have $|(T\cdot\b{1}\{t\in\psi_{k,ij,h}\}/|\psi_{k,ij,h}|)-1|\leq \max(\psi_0^{-1}-1, 1)$. Hence with $\epsilon$ from Proposition \ref{Prop:example_AD2}, we have
\begin{equation*}
\begin{split}
    &\hspace{5mm}
    \sum_{t=1}^T\b{E}[g_{k,j,t}^{2+\epsilon}\mid
    \c{G}_{t-1}]
    \leq
    \|\A_{k,j\cdot}\|^{2+\epsilon} \cdot
    \frac{d_k^{\alpha_{k,r_k}(1+\epsilon/2)}}{T^{1+\epsilon/2}}
    \cdot \|\D^{-1}\H_k^{a,\ast}\|^{2+\epsilon} \\
    &\cdot
    \sum_{t=1}^T
    \Bigg\| \sum_{i=1}^{d_k}  \sum_{h=1}^{\dmk}
    \Bigg(\frac{T \cdot \b{1}\{t\in\psi_{k,ij,h}\} }{|\psi_{k,ij,h}|} -1 \Bigg)
    \Q_{k,i\cdot} \A_{k,i\cdot}'
    (\textnormal{mat}_k(\cF_t)\v_{k,h}
    \v_{k,h}'\textnormal{mat}_k(\cF_t)'
    - \v_{k,h}'\v_{k,h}\bSigma_k )
    \Bigg\|^{2+\epsilon}\\
    &=
    O_P(1)\cdot
    (d_k^{\alpha_{k,1}-\alpha_{k,r_k}} g_s^{-1}
    )^{2+\epsilon} \cdot
    \frac{d_k^{\alpha_{k,r_k}(1+\epsilon/2)}}{T^{\epsilon/2}} \cdot
    d_k^{-\alpha_{k,r_k}(1+\epsilon/2)}
    g_s^{2+\epsilon}
    =
    O_P\Big(\frac{d_k^{\alpha_{k,1}-\alpha_{k,r_k}}}{T^{\epsilon/2}}\Big) = o_P(1),
\end{split}
\end{equation*}
which is sufficient for the conditional Lindeberg condition in \cite{HauslerLuschgy2015} to hold. Then by the stable martingale central limit theorem (Theorem 6.1 in \cite{HauslerLuschgy2015}), we have
\[
\sum_{t=1}^T g_{k,j,t} \to \c{N}(0, \D^{-1}\H_k^{a,\ast} h_{k,j}(\A_{k,j\cdot})(\H_k^{a,\ast})'\D^{-1} )
\;\;\;
\text{$\c{G}^T$-stably as } T\to\infty,
\]
where $h_{k,j}(\A_{k,j\cdot}) = d_k^{2+\alpha_{k,r_k}} \omega_{\psi,k,j} \cdot [\A_{k,j\cdot}'\otimes (\Q_k' \A_k)] \bf\Xi_{F,k} [\A_{k,j\cdot}'\otimes (\Q_k' \A_k)]' $. This completes the proof of Proposition \ref{Prop:example_AD2}.
$\square$

\bibliographystyle{apalike}
\bibliography{tensorrank}

\end{document}